%% file: paper.tex
\documentclass{siamltex}
\usepackage[numbers]{natbib}
\usepackage{amsmath}
\usepackage{amscd}
\usepackage{amssymb}
\usepackage{graphicx}
\usepackage{caption}
\usepackage{url}
\usepackage{color}
\usepackage{listings}
\usepackage{booktabs}
\usepackage{caption}
\usepackage{textcomp}
\usepackage{mathtools}
\usepackage{mcode}
\usepackage{hyperref}
\usepackage[utf8]{inputenc}
\usepackage{tikz}
\usepackage{pgfplots}
\usepackage{environ}
\usepackage{algorithm}
\usepackage{enumerate}
\usepackage[noend]{algpseudocode}
\usepackage{subcaption}
\definecolor{blue}{rgb}{0.2980392156862745, 0.4470588235294118, 0.6901960784313725}
\definecolor{green}{rgb}{0.3333333333333333, 0.6588235294117647, 0.40784313725490196}
\definecolor{red}{rgb}{0.7686274509803922, 0.3058823529411765, 0.3215686274509804}

\newcommand*\Let[2]{\State #1 $\gets$ #2}
\algrenewcommand\algorithmicrequire{\textbf{Input:}}

\newcommand{\dx}{\ \mathrm{d}x}
\newcommand{\ds}{\ \mathrm{d}s}

\definecolor{darkblue}{rgb}{0.00,0.00,0.55}
\definecolor{black}{rgb}{0.00,0.00,0.00}
\hypersetup{
    linkcolor = darkblue,
    anchorcolor = darkblue,
    citecolor = darkblue,
    filecolor = darkblue,
    urlcolor = darkblue,
    colorlinks  = true,
}

\makeatletter
\providecommand{\env@tikzpicture@save@env}{}
\providecommand{\env@tikzpicture@process}{}

\makeatother

\title{The computation of disconnected bifurcation diagrams}

\author{
  P.~E.~Farrell\thanks{Mathematical Institute, University of Oxford, Oxford, UK.
    Center for Biomedical Computing, Simula Research Laboratory, Oslo, Norway
    (\texttt{patrick.farrell@maths.ox.ac.uk}).}
  \and
  C.~H.~L.~Beentjes\thanks{Mathematical Institute, University of Oxford, Oxford, UK (\texttt{beentjes@maths.ox.ac.uk}).}
  \and
  \'A.~Birkisson\thanks{Mathematical Institute, University of Oxford, Oxford, UK (\texttt{birkisson@maths.ox.ac.uk}).
    This research is funded by EPSRC grants EP/K030930/1 and EP/M019721/1, by a Clarendon Fund Scholarship, by a New College Graduate Scholarship,
    by the European Research Council under the European Union's Seventh Framework Programme (FP7/2007-2013)/ERC grant 291068,
    and by a Center of Excellence grant from the Research
    Council of Norway to the Center for Biomedical Computing at Simula
    Research Laboratory. The authors would like to thank L.~N.~Trefethen for useful discussions and G.~N.~Wells for furnishing the hyperelastic solver used in section \ref{sec:beam_pde}.}
  }

\begin{document}
\maketitle

\begin{abstract}
Arclength continuation and branch switching are enormously successful algorithms
for the computation of bifurcation diagrams. Nevertheless, their combination
suffers from three significant disadvantages. The first is that they attempt to
compute only the part of the diagram that is continuously connected to the
initial data; disconnected branches are overlooked. The second is that the
subproblems required (typically determinant calculation and nullspace
construction) are expensive and hard to scale to very large discretizations. The
third is that they can miss connected branches associated with nonsimple
bifurcations, such as when an eigenvalue of even multiplicity crosses the
origin. Without expert knowledge or lucky guesses, these techniques alone can
paint an incomplete picture of the dynamics of a system.

In this paper we propose a new algorithm for computing bifurcation diagrams,
called deflated continuation, that is capable of overcoming all three of these
disadvantages. The algorithm combines classical continuation with a deflation
technique that elegantly eliminates known branches from consideration, allowing
the discovery of disconnected branches with Newton's method. Deflated
continuation does not rely on any device for detecting bifurcations and does not
involve computing eigendecompositions; all subproblems required in deflated
continuation can be solved efficiently if a good preconditioner is available
for the underlying nonlinear problem.
We prove sufficient conditions for the convergence of Newton's method
to multiple solutions from the same initial guess, providing insight into which
unknown branches will be discovered. We illustrate the success of the method on
several examples where standard techniques fail.

\end{abstract}



\section{Introduction}
We consider numerical methods for computing the solutions of
\begin{equation} \label{eqn:fundament}
f(u, \lambda) = 0,
\end{equation}
where $f: U \times \mathbb{R} \to Y$ is the $C^1$ problem residual, $U$ and $Y$
are isomorphic Banach spaces, $u \in U$ is referred to as the solution, and
$\lambda \in \mathbb{R}$ is referred to as the parameter. In our applications,
\eqref{eqn:fundament} typically represents the residual of a stationary ordinary or partial
differential equation, along with boundary conditions. The associated bifurcation
diagram visualizes how the behaviour of a functional of the solutions changes as
$\lambda$ is varied over some interval of interest $[\lambda_{\textrm{min}},
\lambda_{\textrm{max}}]$.

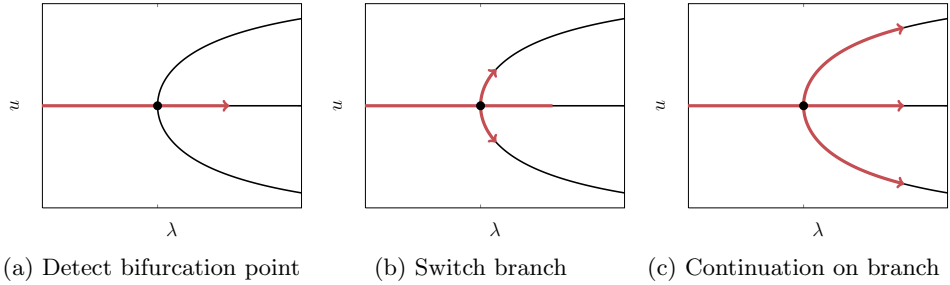
\begin{figure}
\centering
\subcaptionbox{Detect bifurcation point \label{fig:sbif1}}{
\scalebox{0.3}{\input{tikz/standardbif1.tikz}
}}
\
\subcaptionbox{Switch branch\label{fig:sbif2}}{
\scalebox{0.3}{\input{tikz/standardbif2.tikz}
}}
\
\subcaptionbox{Continuation on branch\label{fig:sbif3}}{
\scalebox{0.3}{\input{tikz/standardbif3.tikz}
}}
\caption{Sketch of switching continuation. Arclength (or pseudo-arclength) continuation is applied
on a branch as in Figure \ref{fig:sbif1}. When a bifurcation point is detected,
the nullspace of the Fr\'echet derivative there is computed and is used to
switch branch, depicted in Figure \ref{fig:sbif2}. Once
a point on the emanating branch is known, continuation traces out
the branch as in Figure \ref{fig:sbif3}.}
\label{fig:branchswitch}
\end{figure}

Arclength continuation and branch switching
\cite{keller1977,doedel1981,golubitsky1985,cliffe2008,seydel2010} are central techniques in the
computational analysis of \eqref{eqn:fundament} and are routinely used
throughout science and engineering. Given an initial point $(u_0, \lambda_0)$ on
a branch, arclength continuation (or its popular variant, pseudo-arclength
continuation) robustly traces out the remainder of that branch. It parameterizes
the solution and parameter $(u(s), \lambda(s))$ as a function of arclength along
the curve $s$ from the initial point, as this allows the method to continue
through fold bifurcations. It applies a predictor computed from previous
solutions to estimate the solution and parameter for $s + \Delta s$, and
corrects this guess with a solver such as Newton's method. 

Branch
switching algorithms attempt to detect bifurcation points along a branch and to
construct initial solutions on the branches emanating from it, Figure
\ref{fig:branchswitch}. The detection step typically relies on the computation
of the bifurcation test functional
\begin{equation}\label{eqn:tau_signdet}
\tau(u, \lambda)=\text{sign}(\det J(u,\lambda))=\text{sign}\left(\prod_i \mu_i(u,\lambda)\right),
\end{equation}
where $\mu_i$ are the eigenvalues of the (discretized) Jacobian $J$, including
multiplicities.  This test functional is cheaply computable if an LU
decomposition of $J$ is already available from a
continuation step \cite{seydel2010}, but is very difficult to estimate if a
Krylov method is used.
Once a bifurcation point has been identified, initial guesses for solutions on
the emanating branches are constructed from the nullspace of $J$ there. Once
one solution on each emanating branch is known, arclength continuation is used
to complete the branch. Henceforth, this combination of arclength continuation
and branch switching will be referred to as \emph{switching continuation}.

Switching continuation computes \emph{fragments} of bifurcation
diagrams: it attempts to compute the part of the bifurcation diagram that is
continuously connected to the initial point $(u_0, \lambda_0)$. However, it is
often the case that the bifurcation diagram is not connected, with multiple
branches that do not meet at bifurcation points. For example, pitchfork and
transcritical bifurcations are not generic; they are destroyed under
perturbation \cite[Chapter IV]{golubitsky1985}, and subsequently the complete bifurcation diagram cannot be
computed in one pass with the approach described above. Other examples will be given in
section \ref{sec:examples}. In these cases, the diagram returned by switching
continuation along a single path is incomplete, giving an unsatisfactory picture
of the solutions to \eqref{eqn:fundament}.

\begin{figure}
\centering
\subcaptionbox{Initial continuation \label{fig:dbif1}}{
\scalebox{0.3}{\input{tikz/deflationbif1.tikz}
}}
\
\subcaptionbox{Apply deflation\label{fig:dbif2}}{
\scalebox{0.3}{\input{tikz/deflationbif2.tikz}
}}
\
\subcaptionbox{Continuation of branches\label{fig:dbif3}}{
\scalebox{0.3}{\input{tikz/deflationbif3.tikz}
}}
\caption{Sketch of deflated continuation, the algorithm proposed in this work. Continuation is applied
on a branch as in Figure \ref{fig:dbif1}. Along the branch, we stop and fix the
parameter $\lambda$, then attempt a deflation step to find multiple solutions for
this parameter as in Figure \ref{fig:dbif2}. If the deflation step is
successful, we have points on multiple branches and can continue
these branches as in Figure \ref{fig:dbif3}.}
\label{fig:deflationbifurcation}
\end{figure}
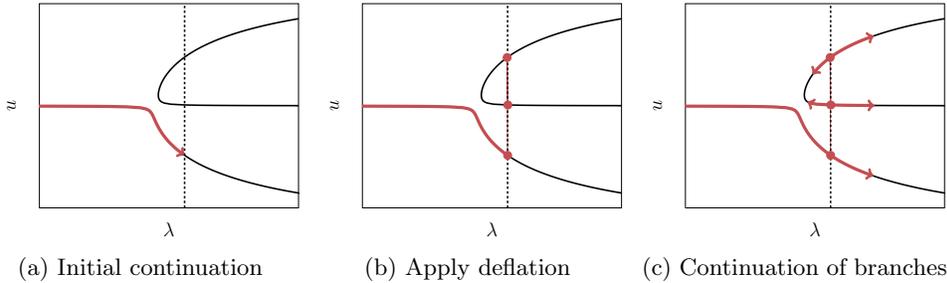

In this work we develop and analyze an alternative algorithm, \emph{deflated continuation}, that is capable
of discovering disconnected branches from known ones, without requiring that a
bifurcation point connect the two, Figure \ref{fig:deflationbifurcation}. At
the heart of the method is the deflation of known solutions
\cite{brown1971,farrell2014}. Deflation is a technique that systematically modifies a
nonlinear problem to guarantee that Newton's method will not converge to
a known root, thus enabling unknown roots to be discovered from the same initial
guess. Fix $\lambda$ in \eqref{eqn:fundament} to yield the nonlinear problem
\begin{equation} \label{eqn:lambda_fixed}
F(u) = 0,
\end{equation}
where $F: U \to Y$.  Suppose Newton's method is applied to
$F$ from initial guess $u_0$ to yield the solution $u_1^*$, with
the Fr\'echet derivative $F'(u_1^*)$ nonsingular. Suppose
further that we suspect that \eqref{eqn:lambda_fixed} permits solutions
other than $u_1^*$, but no additional initial guesses are available.
We thus construct a modified problem
\begin{equation} \label{eqn:deflatedproblem}
G(u) = M(u; u_1^*) F(u),
\end{equation}
via the application of a \emph{deflation operator} $M(u; u_1^*)$ to the 
residual $F$. By construction, this deflated residual satisfies two properties. The
first is the preservation of solutions of $F$, i.e.
for $u \neq u_1^*$, $G(u) = 0$ iff $F(u) = 0$. The second is that Newton's method
applied to $G$ will not discover $u_1^*$ again,
as
\begin{equation}
\liminf_{u \rightarrow u_1^*} \|G(u)\| > 0,
\end{equation}
i.e.~along any sequence converging to the known root, the deflated
residual does not converge to zero. Thus,
if Newton's method applied to $G$ converges from $u_0$,
it will converge to a distinct solution
$u_2^* \neq u_1^*$. The process can then be repeated until no more solutions
are found from $u_0$ in a specified number of Newton iterations.
In this
work, we use the shifted deflation operator
\begin{equation} \label{eqn:deflation}
M(u; u_1^*) = \left(\frac{1}{\|u - u_1^*\|^p} + \sigma\right)\mathbb{I},
\end{equation}
where $\mathbb{I}$ is the identity on $Y$, $p$ is the power, and $\sigma$ is the
shift. All of the examples below use $p = 2$ and $\sigma = 1$. Importantly,
it is possible to efficiently solve the Newton step for $G$ if a good preconditioner
is available for the Newton step of $F$. For
more details, see \cite{abbott1977,farrell2014,farrell2015a}.

We present the proposed bifurcation algorithm in section
\ref{sec:continuation_and_deflation}. To analyze its behaviour, in section
\ref{sec:convergence_analysis} we develop an initial theory of
\emph{multiconvergence} of Newton's method. We derive novel sufficient
conditions under which Newton's method will converge to two different solutions,
starting from the same initial guess. In section \ref{sec:examples}, the
method is applied to several problems on which switching continuation fails.

\section{Deflated continuation} \label{sec:continuation_and_deflation}
\begin{algorithm}
\caption{Deflated continuation.}
\label{alg:mainalg}

\begin{algorithmic}[1]
  \Require{Initial parameter value $\lambda_{\textrm{min}}$.}
  \Require{Final parameter value $\lambda_{\textrm{max}} > \lambda_{\textrm{min}}$.}
  \Require{Step size $\Delta \lambda > 0$.}
  \Require{Nonlinear residual $f(u, \lambda)$.}
  \Require{Deflation operator $M(u; u^*)$.}
  \Require{Initial solutions $S(\lambda_{\textrm{min}})$ to $f(\cdot, \lambda_{\textrm{min}})$.}
  \Statex
  \Let{$\lambda$}{$\lambda_{\textrm{min}}$}
  \While{$\lambda < \lambda_{\textrm{max}}$}
    \Let{$F(\cdot)$}{$f(\cdot, \lambda + \Delta \lambda)$} \Comment{Fix the value of $\lambda$ to solve for.}
    \Let{$S(\lambda + \Delta \lambda)$}{$\varnothing$}
    \For{$u_0 \in S(\lambda)$}                            \Comment{Continue known branches.}
      \State{apply Newton's method to $F$ from initial guess $u_0$.}
      \If{solution $u^*$ found}
        \Let{$S(\lambda + \Delta \lambda)$}{$S(\lambda + \Delta \lambda) \cup \{u^*\}$} \Comment{Record solution.}
        \Let{$F(\cdot)$}{$M(\cdot; u^*)F(\cdot)$} \Comment{Deflate solution.}
      \EndIf
    \EndFor
    \Statex
    \For{$u_0 \in S(\lambda)$}                            \Comment{Seek new branches.}
      \Let{success}{true}
      \While{success}
        \State{apply Newton's method to $F$ from initial guess $u_0$.}
        \If{solution $u^*$ found}                                   \Comment{New branch found.}
          \Let{$S(\lambda + \Delta \lambda)$}{$S(\lambda + \Delta \lambda) \cup \{u^*\}$} \Comment{Record solution.}
          \Let{$F(\cdot)$}{$M(\cdot; u^*)F(\cdot)$} \Comment{Deflate solution.}
        \Else{}
          \Let{success}{false}
        \EndIf
      \EndWhile
    \EndFor

    \Let{$\lambda$}{$\lambda + \Delta \lambda$}
  \EndWhile
  \State \Return{$S$}
\end{algorithmic}
\end{algorithm}

Let $[\lambda_{\textrm{min}}, \lambda_{\textrm{max}}]$ be the interval  of interest
for the parameter $\lambda$, and let $\Delta \lambda$ be the continuation step size. 
For a given $\lambda$, let $S(\lambda) \subset U$ denote the set of known solutions to
\eqref{eqn:lambda_fixed}. Given $S(\lambda)$, Algorithm
\ref{alg:mainalg}
constructs $S(\lambda + \Delta \lambda)$ as follows. In the first pass (lines
5--9), known solutions are continued with standard classical continuation; each
known solution $u_0 \in S(\lambda)$ is used as initial guess for $f(\cdot, \lambda + \Delta \lambda)$ in turn. If some
solutions are not successfully continued, the algorithm proceeds with the other
branches regardless. (This can happen at fold bifurcations, for example.) As
each solution is continued, it is recorded and deflated. In the second pass
(lines 10--18), each initial guess is again considered in turn. Deflation
guarantees that Newton's method will not return to the known branch, and hence if
Newton's method converges, it will converge to a previously unknown solution
(line 14). Each initial guess is attempted repeatedly until no solution is found
within a certain number of Newton iterations.  (Recall that Newton's method is
undecidable, i.e.~it is impossible to decide in general if Newton's method will eventually
converge for a given initial guess \cite{blum1998}.) Once all initial guesses
have been exhausted, the algorithm increments $\lambda$ and continues until the
end of the interval has been reached.

\subsection{Variants of the algorithm}

Various modifications to the basic algorithm are possible. The analyst may
decide to seek unknown branches with a step size larger than $\Delta \lambda$,
to reduce the effort spent on unsuccessful Newton iterations.  The continuation
and discovery stages are independent and may be executed in parallel; one group
of processors can continue known solutions forwards, while other groups follow
behind, seeking new solutions to continue.

If the system \eqref{eqn:fundament} has a finite symmetry group $\mathcal{G}$
such that for all $g \in \mathcal{G}$,
\begin{equation}
f(u, \lambda) = 0 \iff f(gu, \lambda) = 0,
\end{equation}
then when a solution $u$ is discovered, its actions $\mathcal{G}u$ should be
recorded and deflated as well, assuming that it is possible to represent each
$gu$ exactly with the discretization employed. If the discretization does not
respect this symmetry (e.g. a finite element discretization on an unstructured
mesh), then the projection of $gu$ should be used as initial guess for
Newton's method instead. Deflating infinite symmetry groups will be studied in
future research.

If the system \eqref{eqn:fundament} has a trivial branch $\bar{u}$ such that
$f(\bar{u}, \lambda) = 0$ for all $\lambda$, then this branch must be excluded
from the set of initial guesses to use in deflation. This is because the initial
residual of \eqref{eqn:deflatedproblem} will evaluate to $0/0$. As such trivial
branches are obvious from the equations, the simplest approach is just to
deflate any trivial solutions away before beginning Algorithm
\ref{alg:mainalg}.

In the discovery stage, problem-specific guesses other than the previous
solutions may be employed; for example, in nonlinear eigenproblems it may be
useful to use the eigenmodes of an associated linear problem. It may be
necessary to break the symmetry of the guesses: if the system
\eqref{eqn:fundament} has a $\mathbb{Z}_2$ symmetry $R$ such that $f(Ru,
\lambda) = Rf(u, \lambda)$, then if Newton's method is initialized with a
symmetric initial guess satisfying $Ru_0 = u_0$ then all subsequent iterates will
also remain symmetric. This will cause nonconvergence to nonsymmetric
solutions, such as those introduced at a symmetry-breaking bifurcation. In this
regard it may be advantageous to deliberately break the symmetry of the
discretization, or if this is not possible (such as when using a spectral
method), to deliberately break the symmetry of the initial guesses.

It is straightforward in principle to employ other continuation approaches in
Algorithm \ref{alg:mainalg}: if arclength continuation is used, then in the
first pass each branch is synchronized at $\lambda + \Delta \lambda$, deflation
is applied to seek new branches, and the process is repeated.

\section{Convergence analysis} \label{sec:convergence_analysis}
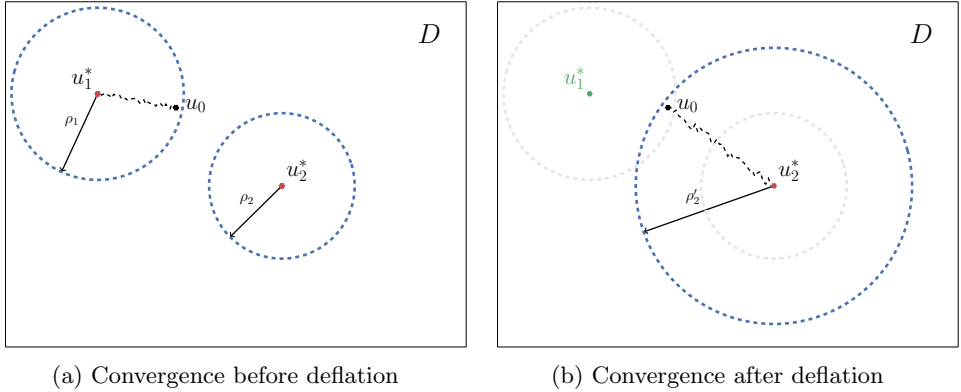
\begin{figure}
\centering
\subcaptionbox{Convergence before deflation \label{fig:sufficientconditions1}}{
\scalebox{0.5}{\input{tikz/sufficientconditions1.tikz}
}}
\
\subcaptionbox{Convergence after deflation \label{fig:sufficientconditions2}}{
\scalebox{0.5}{\input{tikz/sufficientconditions2.tikz}
}}
\caption{Sketch of the regions of convergence around solutions $u_1^*$ and $u_2^*$ before
and after deflation.
Before deflation, the initial guess $u_0$ converges to $u_1^*$; after deflating
this solution, the region of convergence around $u_2^*$ expands and $u_0$ now lies
within it.}
\label{fig:sufficientconditions}
\end{figure}

The central question in the analysis of Algorithm \ref{alg:mainalg} is: under
what circumstances will unknown branches be discovered, and under what
circumstances will they be missed? Given an initial guess $u_0$, we wish to
derive sufficient conditions that guarantee convergence to at least two
solutions $u_1^*$ and $u_2^*$ with Newton's method and deflation, Figure
\ref{fig:sufficientconditions}. In the context of Algorithm \ref{alg:mainalg},
$u_0$ is the known solution for $f(\cdot, \lambda)$, $u_1^*$ is the solution on
the same branch for $f(\cdot, \lambda + \Delta \lambda)$, and $u_2^*$ is another solution
to $f(\cdot, \lambda + \Delta \lambda)$ on a different branch.

The best-known theorem of convergence for Newton's method is the theorem of
Kantorovich \cite{kantorovich1948}, who first formulated and analyzed Newton's
method in Banach spaces. We state the theorem (and all subsequent results) in affine-covariant form
\cite{deuflhard2011}.
\begin{theorem}[Affine-covariant Newton--Kantorovich \cite{kantorovich1948}]\label{thm:kantorovich}
Let $F:D\to Y$ be a continuously Fr\'echet differentiable function on the open
convex subset $D\subseteq U$. Given $u_0\in D$, assume that \\
\begin{enumerate}[i)]
\item $F'(u_0)^{-1}$ exists; let $\alpha = \|F'(u_0)^{-1}F(u_0)\|$; \\
\item $\|F'(u_0)^{-1}\left(F'(u)-F'(v)\right)\|\leq \omega_0 \|u-v\|$ for all
$u,v\in D$; \\
\item $h_0 = \alpha  \omega_0 \leq \frac{1}{2}$; \\ 
\item $\mathcal{B}=\bar{B}(u_0,\rho_0)\subset D$ for
$\rho_0=({1-\sqrt{1-2h_0}})/{\omega_0}$, where $B$ defines an open ball. \\
\end{enumerate}
Then the Newton sequence from $u_0$ is well-defined and remains within the ball $\mathcal{B}$. 
A solution $u^*\in \mathcal{B}$ with $F(u^*)=0$ exists, and the
Newton sequence converges to it.
{Furthermore, if we define $\rho^+=({1+\sqrt{1-2h_0}})/{\omega_0}$, then $u^*$
is unique within $D\cap B(u_0,\rho^+)$}.
\end{theorem}

One of the main features of this theorem is that all of its assumptions except
for Lipschitz continuity are verified at the initial guess $u_0$. Convergence
can be assured \emph{a priori}, without needing to assume the existence of a
root beforehand.

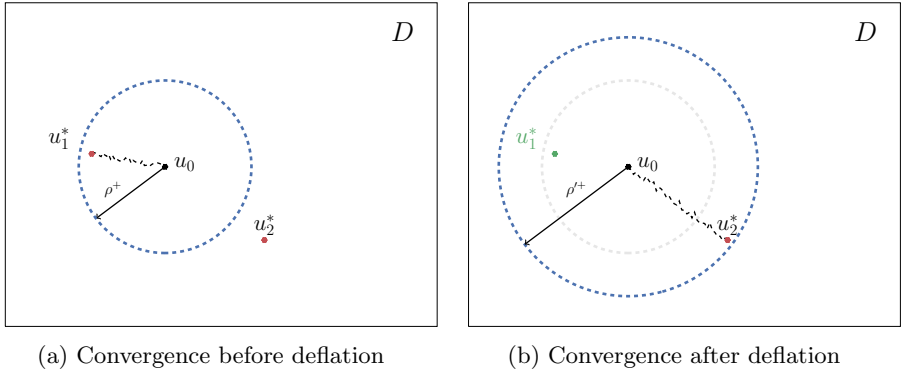
\begin{figure}
\centering
\subcaptionbox{Convergence before deflation \label{fig:centering1}}{
\scalebox{0.5}{\input{tikz/centering1.tikz}
}}
\
\subcaptionbox{Convergence after deflation \label{fig:centering2}}{
\scalebox{0.5}{\input{tikz/centering2.tikz}
}}
\caption{Sketch of why sufficient conditions for multiconvergence cannot be
based on the Newton--Kantorovich theorem. The Newton--Kantorovich theorem
describes a ball of convergence centred at the initial guess $u_0$.
In order to show convergence to
multiple solutions $u_1^*$ and $u_2^*$ we need the convergence region to grow, i.e. $\rho'^+ >
\rho^+$. This would imply that the deflated root lies within the new region of
convergence, which poses regularity problems on the Fr\'echet derivative of the
deflated function in the convergence region.}
\label{fig:centering}
\end{figure}

Nevertheless, this theorem is not a suitable foundation for the purpose at hand.
Suppose there exist $u_1^*$ and $u_2^*$ with $F(u_1^*) = F(u_2^*) = 0$ and
$u_1^* \neq u_2^*$. Now consider an initial guess $u_0$ that provably converges
to $u_1^*$ by the Newton--Kantorovich theorem. The result is a $\rho^+$ such
that $u_1^* \in B(u_0, \rho^+)$ and $u_2^* \notin B(u_0, \rho^+)$, Figure
\ref{fig:centering1}. In order to prove convergence of the deflated function
$M(u; u_1^*) F(u)$, we would need to establish a $\rho'^+ > \rho^+$ such that
$u_2^* \in B(u_0, \rho'^+)$. However, this would imply that $u_1^* \in B(u_0,
\rho'^+)$, Figure \ref{fig:centering2}. The assumptions of the
Newton--Kantorovich theorem imply that the Fr\'echet derivative is invertible
everywhere in the ball, but the Fr\'echet derivative of the deflated
function is not defined at $u_1^*$, and hence the assumptions cannot hold after
deflation.
The same argument holds for the Newton--Mysovskikh theorem \cite{mysovskikh1949}.

We therefore seek to base our analysis on results whose conditions are verified
at the roots themselves, instead of at the initial guess. The theorem we will
build upon is the Rall--Rheinboldt theorem \cite{rall1974,rheinboldt1978}, again
stated in affine-covariant form.
\begin{theorem}[Affine-covariant Rall-Rheinboldt \cite{rall1974,rheinboldt1978}]\label{thm:acRR}
Let $F:D\to Y$ be a continuously Fr\'echet differentiable function on the open
convex subset $D\subseteq U$. Suppose that there exists a $u^*\in D$ such that
$F(u^*)=0$, and suppose further that \\
\begin{enumerate}[i)]
\item $F'(u^*)^{-1}$ exists; \\
\item $\|F'(u^*)^{-1}\left(F'(u)-F'(v)\right)\|\leq \omega^* \|u-v\|$ for all
$u,v\in D$. \\
\end{enumerate}
Then any $\rho^*\leq 2/(3\omega^*)$ such that
$\mathcal{B}={B}(u^*,\rho^*)\subset D$ has the property that starting at
$u_0\in\mathcal{B}$, the Newton sequence is well-defined and remains within
$\mathcal{B}$. The Newton sequence converges to $u^*\in \mathcal{B}$.
{Furthermore, if we define $\rho^+={1}/{\omega^*}$, then $u^*$ is unique
within $D\cap B(u^*,\rho^+)$}.
\end{theorem}

A crucial ingredient of this theorem is the affine covariant Lipschitz
continuity of the Fr\'echet derivative $F'$. Before extending this theorem
to the deflated case, we first give a lemma regarding
the product of Lipschitz continuous functions.
\begin{lemma}[Product of Lipschitz continuous functions]\label{lem:product_Lipschitz}
Let $X, Y$ and $Z$ be Banach spaces and let $L(Y, Z)$ be the vector space of bounded linear
operators from $Y$ to $Z$ with induced operator norm. Let $F:X \to Y$ and $G: X \to L(Y, Z)$ be Lipschitz
continuous functions on the open subset $D\subseteq X$ with Lipschitz constants
$\omega_F$ and $\omega_G$ respectively. Assume further that $F$ and $G$ are
bounded on $D$, i.e.\ there exist $N_F, N_G \in \mathbb{R}$ such that
$\|F(x)\|<N_F$ and $\|G(x)\| < N_G$ for all $x \in D$.  Then the product
$GF:X\to Z$ is bounded and Lipschitz continuous on $D$ with Lipschitz constant
$(N_F\omega_G+N_G\omega_F)$.
\end{lemma}

\begin{proof}
Let $x,y\in D$. As both $F$ and $G$ are bounded on $D$ their product is bounded
as well:
\begin{equation}
\|G(x)F(x)\|\leq \|G(x)\|\|F(x)\|\leq N_G N_F< \infty.
\end{equation}
Furthermore,
\begin{align}
\|G(x)F(x)-G(y)F(y)\|&=\|G(x)F(x)-G(x)F(y)+G(x)F(y)-G(y)F(y\| \nonumber \\ \nonumber
&\leq \|G(x)F(x)-G(x)F(y)\|+\|G(x)F(y)-G(y)F(y)\|\\ \nonumber
&\leq \|G(x)\| \|F(x)-F(y)\| + \|G(x)-G(y)\|\|F(y)\|\\ \nonumber
&\leq N_G\|F(x)-F(y)\|+N_F\|G(x)-G(y)\|\\ \nonumber
&\leq N_G\omega_F \|x-y\| + N_F\omega_G \|x-y\|\\
&=(N_F\omega_G+N_G\omega_F)\|x-y\|,
\end{align}
which proves the claim.
\end{proof}

We now consider the situation where one solution $u_1^*$ is known and has
been deflated. We state sufficient conditions on the original residual and
deflation operator that guarantee convergence to another solution $u_2^*$.
\begin{theorem} \label{thm:known_u1}
Let $F: D\to Y$ be a continuously Fr\'echet differentiable function on the open
convex subset $D\subseteq U$. Suppose there exists $u_2^* \in D$ such that $F(u_2^*)=0$. 
Further assume there exists
$u_1^*\in D$, $u_1^*\neq u_2^*$, such that $F(u_1^*)=0$. This solution is
deflated with a deflation operator $M(\cdot;u_1^*):D\setminus\{u_1^*\}\to GL(Y,Y)$.
Suppose there exists an open bounded convex subset $E\subseteq
D\setminus\{u_1^*\}$ with $u_2^* \in E$ such that the following conditions
hold: \\
\begin{enumerate}[i)]
\item $F'(u_2^*)^{-1}$ exists; \\

\item $\|F'(u_2^*)^{-1}\left(F'(u)-F'(v)\right)\|\leq \omega^* \|u-v\|$ for all
$u,v\in E$; \\

\item $M(u; u_1^*)$ is continuously Fr\'echet differentiable for all $u\in E$; \\

\item $\|M'(u; u_1^*)-M'(v; u_1^*)\|\leq
\omega_{{M}'} \| u-v\|$ for all $u, v \in E$. \\
\end{enumerate}
Then there exists a $\rho>0$ such that the Newton sequence from $u_0 \in \mathcal{B} =
B(u_2^*, \rho)$ on the deflated function ${M}(u;
u_1^*)F(u)$ is well-defined, remains in $\mathcal{B}$ and converges to $u_2^*\in
\mathcal{B}$.
\end{theorem}

If the norm on $U$ is twice continuously differentiable on $E$, the deflation
operator \eqref{eqn:deflation} satisfies these conditions. In this case, we can
use the composition rule for differentiable functions to show that the deflation
operator is in turn
twice differentiable on $E$. This implies that the deflation operator is
Lipschitz continuous, as $E$ is bounded. For $2 \leq p < \infty$, the $L^p$ norm
is at least twice continuously differentiable on any open subset not containing
zero \cite[Theorem 8]{sundaresan1967}. More generally, if the Banach space $U$
is isomorphic to a Hilbert space, then it can be equipped with twice
differentiable norms \cite{leonard1973,fabian2011}. Thus, the conditions
demanded are satisfied in typical cases of interest.

\begin{proof}
As $F(u)$ and ${M}(u;u_1^*)$ are continuously Fr\'echet differentiable on $E$,
they are Lipschitz continuous as well by boundedness of $E$.
Lipschitz continuity implies that the operators are bounded on $E$ and thus
$F(u),F'(u),{M}(u;u_1^*)$ and ${M}'(u;u_1^*)$ are all bounded on $E$.
As a result the Fr\'echet derivative of the deflated operator
\begin{equation}
\left({M}(u;u_1^*)F(u)\right)'={M}(u;u_1^*)F'(u)+{M}'(u;u_1^*)F(u)
\end{equation}
is Lipschitz continuous by use of the triangle inequality and Lemma
\ref{lem:product_Lipschitz}.

Since $u_2^*$ is a root of $F$, the Fr\'echet derivative of the
deflated residual there is $M(u_2^*; u_1^*) F'(u_2^*)$. For any
$u \in D$ the deflation operator $M(u; u_1^*) \in GL(Y, Y)$ and is
thus invertible. The Fr\'echet derivative of the deflated residual
is thus invertible at $u_2^*$ with
\begin{equation}
\left.\left( \left[M(u; u_1^*) F(u)\right]' \right)^{-1}\right|_{u_2^*} = 
 F'(u_2^*)^{-1} M(u_2^*; u_1^*)^{-1}.
\end{equation}

Combining these facts, there exists an (affine covariant)
$\tilde{\omega}_2>0$ such that 
\begin{equation}
\left\| \left( \left[{M}(u_2^*;u_1^*) F(u_2^*) \right]' \right)^{-1}
\left[({M}(u;u_1^*)F(u))'-({M}(v;u_1^*)F(v))'\right]\right\|\leq
\tilde{\omega}_2\|u-v\|,
\end{equation}
for all $u,v\in E$. Hence the conditions of Theorem \ref{thm:acRR} are satisfied
for both $F(u)$ and ${M}(u;u_1^*)F(u)$, and it can be applied to prove the
claim.
\end{proof}

We are now in a position to state sufficient conditions for convergence
to two solutions with deflation and Newton's method. The proof applies
the previous theorem, Theorem \ref{thm:known_u1}, and the Rall--Rheinboldt
theorem, Theorem \ref{thm:acRR}.
\begin{theorem}[Deflated Rall-Rheinboldt \cite{beentjes2015}] \label{thm:deflated_RR}
Let $F: D\to Y$ be a continuously Fr\'echet differentiable function on an open
subset $D\subseteq U$. Suppose there exist $u_1^*,u_2^* \in D$ such that
$F(u_1^*)=F(u_2^*)=0$, $u_1^* \neq u_2^*$.
Let $E_1$ be an open bounded convex subset such that $E_1\subset D\setminus
\{u_2^*\}$ and $u_1^*\in E_1$. Furthermore let $E_2$ be an open bounded convex
subset such that $E_2\subset D\setminus \{u_1^*\}$ and $u_2^*\in E_2$.
Let $M(\cdot;u_1^*):D\setminus\{u_1^*\}\to GL(Y,Y)$ be a deflation
operator such that the following conditions hold: \\
\begin{enumerate}[i)]
\item $F'(u_1^*)^{-1}$ and $F'(u_2^*)^{-1}$ exist; \\

\item $\|F'(u_1^*)^{-1}\left(F'(u)-F'(v)\right)\|\leq \omega_1^* \|u-v\|$ for all $u,v\in E_1$; \\

\item $\|F'(u_2^*)^{-1}\left(F'(u)-F'(v)\right)\|\leq \omega_2^* \|u-v\|$ for all $u,v\in E_2$; \\

\item $M(u; u_1^*)$ is continuously Fr\'echet differentiable for all $u\in E_2$; \\

\item $\|{M}'(u; u_1^*)-{M}'(v;u_1^*)\|\leq \omega_{{M}'} \|u-v\|$ for all $x,y \in E_2$. \\

\end{enumerate}
Then there exists an $\tilde{\omega}_2>0$ such that for all $u, v\in E_2$ there holds
\begin{equation}
\left\| \left( \left[{M}(u_2^*;u_1^*) F(u_2^*) \right]' \right)^{-1}
\left[({M}(u ;u_1^*)F(u))'-({M}(v;u_1^*)F(v))'\right]\right\|\leq
\tilde{\omega}_2\|u-v\|.
\end{equation}
If $\|u_1^*-u_2^*\|<\rho_1+\rho_2$ for some $\rho_1\leq 2/(3\omega_1^*)$ and
$\rho_2\leq 2/(3\tilde{\omega}_2)$ such that we have
${B}_1=B(u_1^*,\rho_1)\subset E_1$ and ${B}_2=B(u_2^*,\rho_2)\subset E_2$, then
the intersection ${B}_1\cap {B}_2$ is nonempty. Starting from any $u_0 \in
{B}_1\cap {B}_2$, Newton's
method will first converge to $u_1^*\in E_1$ and then after deflation with
${M}(\cdot;u_1^*)$ will converge to $u_2^*\in E_2$.
\end{theorem}


The argument of Theorems \ref{thm:known_u1} and \ref{thm:deflated_RR} can
be applied again to derive sufficient conditions for a single initial guess
to converge to three or more solutions.

A natural question to ask is if Algorithm \ref{alg:mainalg} will recover the
behaviour of switching continuation, i.e.~if it will always discover connected
branches for sufficiently small $\Delta \lambda$. This is discussed in the
following corollary.
\begin{corollary}[Connected roots]\label{conj:connected_roots}
Let $f:D\times \mathbb{R}\to Y$ and suppose there exists a
$\lambda_c\in\mathbb{R}$ such that $f(\cdot,\lambda):D\to Y$ is a continuously
Fr\'echet differentiable function on the open subset $D\subseteq U$ for
$\lambda>\lambda_c$. Furthermore assume that there exists
$u_1^*(\lambda),u_2^*(\lambda) \in D$ such that
$f(u_1^*(\lambda),\lambda)=f(u_2^*(\lambda),\lambda)=0$ and $u_1^*(\lambda) \neq
u_2^*(\lambda)$ for $\lambda>\lambda_c$ and $u_1^*(\lambda_c)=u_2^*(\lambda_c)$.
Assume that for fixed $\lambda>\lambda_c$ all conditions from Theorem
\ref{thm:deflated_RR} hold for the function $f(\cdot,\lambda):D\to Y$ so that
$\rho_1(\lambda),\rho_2(\lambda)\in\mathbb{R}$ as in Theorem
\ref{thm:deflated_RR} are well defined. If
\begin{equation} \label{eqn:conjecture_limit}
\lim_{\lambda\downarrow \lambda_c} \frac{\|u_1^*(\lambda)-u_2^*(\lambda)\|}{\rho_1(\lambda)+\rho_2(\lambda)}<1,
\end{equation}
then an initial guess $u_0\in D$ exists which converges to both
$u_1^*$ and $u_2^*$ using Newton's method and deflation for $\lambda$
sufficiently close to $\lambda_c$.
\end{corollary}

By assumption, $\|u_1^*(\lambda)-u_2^*(\lambda)\| \to 0$ as $\lambda \downarrow
\lambda_c$. As $\rho_1(\lambda) < \|u_1^*(\lambda)-u_2^*(\lambda)\|$ (and
similarly for $\rho_2(\lambda)$), $\rho_1(\lambda) +
\rho_2(\lambda) \to 0$ also. Thus, the evaluation of the
left-hand side of \eqref{eqn:conjecture_limit} requires the application of
L'H\^{o}pital's rule.

A similar formulation applies to the case of branches meeting as $\lambda
\uparrow \lambda_c$, and to more than two roots.  Our practical experience does
indeed suggest that Algorithm \ref{alg:mainalg} is always able to find branches
connected via a bifurcation point; we conjecture that
\eqref{eqn:conjecture_limit} always holds for sufficiently regular functions.

Note that these results are nonconstructive, i.e.~the Lipschitz constants
arising and the resulting radii of convergence are not in general known.  Thus,
it could be the case that the region of multiconvergence is too small to be of
practical use in bifurcation analysis. We therefore apply Algorithm
\ref{alg:mainalg} to several problems of interest in the literature to
investigate the robustness and efficiency of deflated continuation.

\section{Examples} \label{sec:examples}
\subsection{Roots of unity}
\begin{figure}
\centering
\input{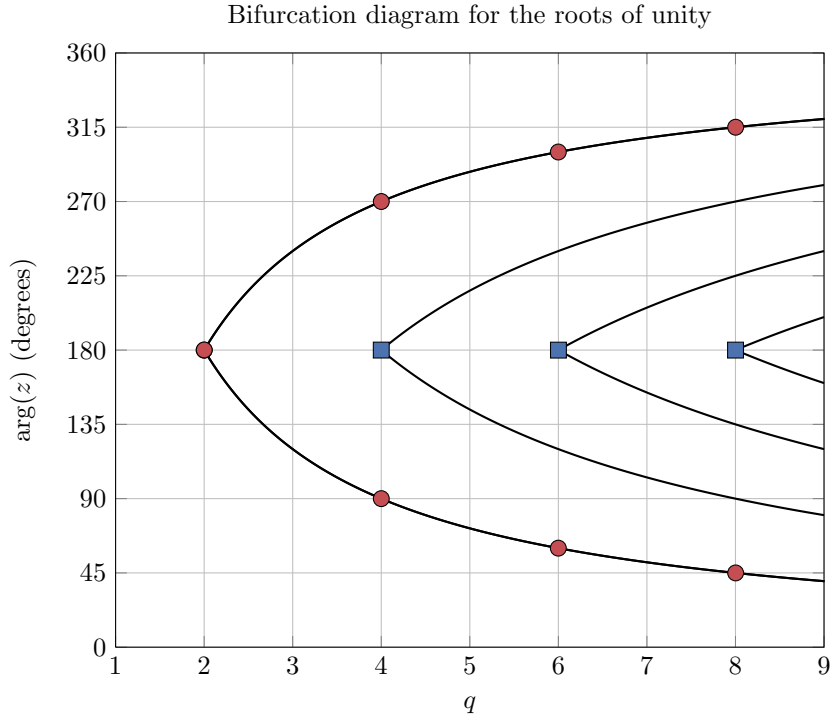}
\caption{Bifurcation diagram for the roots of unity \eqref{eqn:unity}
as a function of exponent $q$. The trivial branch $z = 1$ is not shown.
The bifurcation diagram is disconnected; switching continuation from $q = 2, z = -1$ would only
identify the branch marked with red circles. Blue squares denote the discovery
of disconnected branches with deflation.
}
\label{fig:unity}
\end{figure}

We consider the complex roots of unity
\begin{equation} \label{eqn:unity}
z^q - 1 = 0
\end{equation}
as the exponent $q$ is varied. For $q \in \mathbb{N}_+$, the solutions are
$\exp({2\pi i k}/{q})$ for $k = 1, \dots, q$; this example studies how these
solutions bifurcate for non-integer exponents.

Algorithm \ref{alg:mainalg} was applied to \eqref{eqn:unity} from $q = 2$ to $q
= 9$ with $\Delta q = 0.1$. Deflation was applied with power $p = 2$ and shift $\sigma = 1$. The
resulting bifurcation diagram is shown in Figure \ref{fig:unity}, where the
quantity plotted is the argument of the solution. For $q = 2$, the solutions
are $z = \pm 1$; the solution $z = -1$ bifurcates and the resulting solutions
approach $z = \pm i$ as $q \to 4$. At $q = 4$, $z = -1$ undergoes another
bifurcation, and the process repeats. In general there is a bifurcation at $z =
-1$ for $q \in 2\mathbb{N}_+$, and the resulting branches are mutually
disconnected from each other.

As the bifurcation diagram is disconnected, switching continuation would
identify at most one branch from any given initial solution. By contrast,
deflated continuation identifies the new solutions at $z = -1$
immediately and correctly computes the entire diagram.
\subsection{Deformation of a slender beam} \label{sec:beam_ode}
\begin{figure}
\centering
\subcaptionbox{$\mu = 0$ \label{fig:euler_sym}}{
\scalebox{0.5}{\input{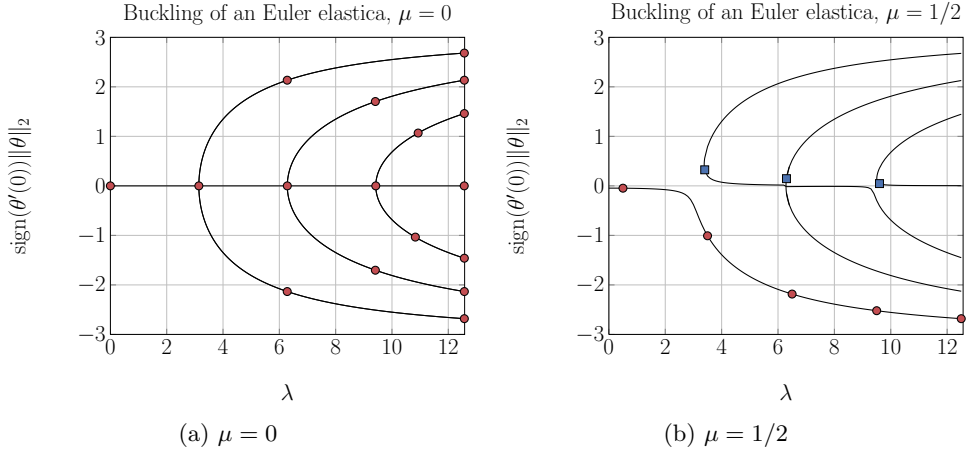}
}}
\
\subcaptionbox{$\mu = 1/2$ \label{fig:euler_asym}}{
\scalebox{0.5}{\input{tikz/eulerbeam-asym.tikz}
}}
\caption{Bifurcation diagrams for the Euler elastica equation \eqref{eqn:euler}
as a function of longitudinal loading $\lambda$, for $\mu = 0$ and $\mu = 1/2$.
For $\mu = 0$, the bifurcation diagram is continuous and both switching continuation
and deflated continuation discover the entire diagram. For $\mu =
1/2$, the bifurcation diagram is disconnected: switching continuation only discovers the part of the diagram labelled with red circles,
whereas deflated continuation correctly computes the entire diagram. Blue squares denote the discovery
of disconnected branches with deflation.}
\label{fig:euler}
\end{figure}

The deformation of a slender vertical beam under loading is governed by Euler's
elastica equation \cite{levien2008}
\begin{align} \label{eqn:euler}
\theta'' + \lambda^2 \sin{(\theta)} = \mu, \qquad \theta(0) = \theta(1) = 0,
\end{align}
where $s$ is the arclength along the beam, $\theta(s)$ is the angle relative to
the vertical axis, $\lambda$ is the longitudinal force and $\mu$ is the transversal
force. This system has long served as a model problem in bifurcation analysis \cite{reiss1969}.

Algorithm \ref{alg:mainalg} was applied to \eqref{eqn:euler}, from $\lambda = 0$
to $\lambda = 4\pi$, with continuation step $\Delta \lambda = 0.1$.
The equation was discretized with $10^4$ piecewise linear finite elements using
FEniCS \cite{logg2011} and PETSc \cite{balay2015}.
In the absence of a transversal force ($\mu = 0$), the initially straight
solution $\theta(s) = 0$ forms a trivial branch and was thus deflated
before beginning Algorithm \ref{alg:mainalg}. 
Newton's method was terminated with failure if convergence
did not occur within $10^2$ iterations.
Deflation was applied with power $p = 2$, shift $\sigma = 1$ and with distances measured
in the $H^1$ norm.
After the forward continuation
pass, arclength continuation backwards in $\lambda$ was performed (without deflation) to
complete the small sections of the bifurcation diagram where branches were not
immediately discovered (cf.~Figure \ref{fig:dbif2}). The functional used was the
$L^2$ norm, signed by $\textrm{sign}(\theta'(0))$.

A series of pitchfork bifurcations at $\lambda = n\pi$ for $n \in \mathbb{N}_+$
(corresponding to the eigenvalues of the associated linear problem) result in
the buckled modes emanating from the trivial branch. As all branches meet at
bifurcation points with the trivial branch, both switching continuation and
deflated continuation compute the entire bifurcation diagram, Figure
\ref{fig:euler_sym}.  However, if a transversal force is applied ($\mu = 1/2$),
the reflective symmetry is destroyed and the symmetric pitchfork bifurcations
degenerate. In this case, the initial branch disconnects from all other
branches, resulting in a disconnected bifurcation diagram, Figure
\ref{fig:euler_asym}. All of these other branches are missed with switching
continuation applied to this path, yielding an incomplete representation of the
dynamics of the system\footnote{It is possible to identify all of these solutions with
switching continuation as follows: set $\mu = 0$ and continue $\lambda$ from 0
to $4\pi$; set $\lambda = 4\pi$ and continue $\mu$ from 0 to $1/2$; set $\mu =
1/2$ and continue $\lambda$ from $4\pi$ to 0. However, this is laborious and
requires expert knowledge of the system; the right continuation strategy may not
be obvious in more complex cases.}. Deflated continuation correctly computes the
bifurcation diagram without continuation along multiple parameters.

\subsection{Nonlinear pendulum} \label{sec:pendulum}
In the previous example, an additional source term destroyed the symmetry of
the bifurcation diagram. This example serves to demonstrate that the same effect
can be achieved by inhomogeneous boundary conditions.
The angle of a pendulum to the vertical is described by the same equation,
\begin{equation} \label{eqn:pendulum}
\theta'' + \sin{\theta} = 0,
\end{equation}
but here we impose inhomogeneous Dirichlet conditions $\theta(0) = \theta(10) =
2$. It is well known that with these boundary conditions this equation permits
multiple solutions \cite{birkisson2013}. One possible way to compute these solutions is to attempt
a homotopy from the linear equation $\theta'' = 0$ via the addition of a
parameter $\varepsilon$ multiplying the nonlinear term:
\begin{equation} \label{eqn:pendulumh}
\theta'' + \varepsilon \sin{\theta} = 0, \qquad \theta(0) = \theta(10) = 2.
\end{equation}
For $\varepsilon = 0$, \eqref{eqn:pendulumh} reduces to a trivial linear
problem; for $\varepsilon = 1$, the problem of interest is recovered. It is
clear that homotopy methods based on switching continuation will identify a solution for $\varepsilon = 1$ only
if there is a branch that continuously connects it to the solution for
$\varepsilon = 0$ \cite[\S 11.3]{nocedal2006}.\footnote{Another approach would be
to consider the associated initial-value problem with boundary conditions $\theta(0) = 2,
\theta'(0) = \varepsilon$. The resulting IVP can be solved for varying values of $\varepsilon$
and the solutions with $\theta(10) = 2$ selected. This shooting approach does not generalize
to higher dimensions, and can be unstable.}

Algorithm \ref{alg:mainalg} was applied to \eqref{eqn:pendulumh}, from $\varepsilon
= 0$ to $\varepsilon = 1$, with continuation step $\Delta \varepsilon = 10^{-2}$.
The equation was discretized with $10^4$ standard piecewise linear finite elements
using FEniCS and PETSc.
The same deflation and solver settings were used as in the previous example.
The functional considered was the product of the derivative at the left endpoint
and the $H^1$ norm of the solution.

\begin{figure}
\centering
\input{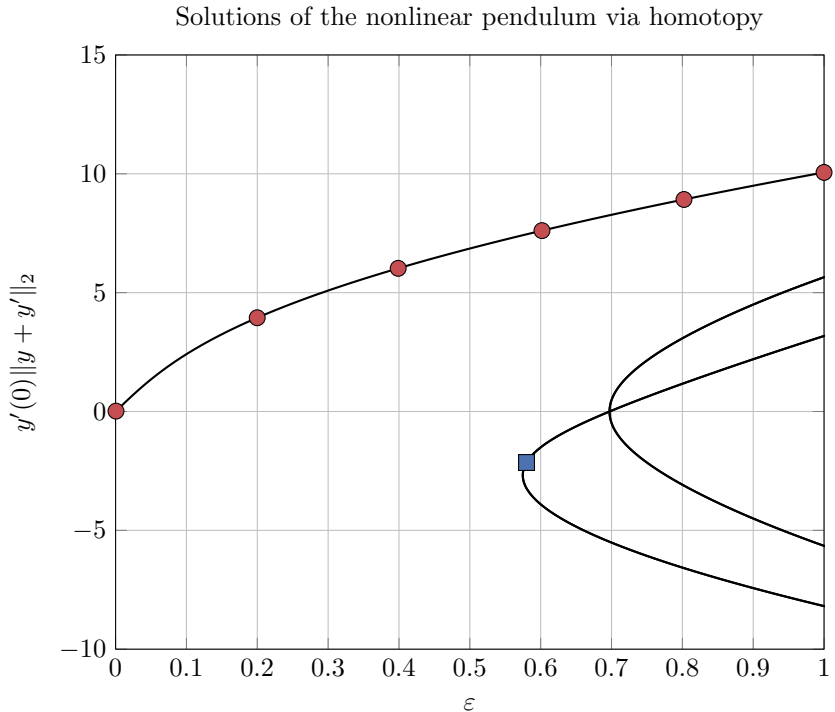}
\caption{Bifurcation diagram for the nonlinear pendulum \eqref{eqn:pendulumh}
as a function of homotopy parameter $\varepsilon$. As before, switching continuation
only discovers that part of the diagram labelled with
red circles. The blue square denotes the discovery
of a disconnected branch with deflation.}
\label{fig:pendulum}
\end{figure}

The resulting bifurcation diagram is shown in Figure \ref{fig:pendulum}. For
$\varepsilon = 0$, the problem has a unique solution; as continuation is applied
to this branch, no bifurcation points are encountered, and hence with switching
continuation only one solution would be identified for $\varepsilon = 1$.
As previously mentioned, a major difficulty with such homotopy methods is that the
resulting bifurcation diagram must continuously connect the solution for
$\varepsilon = 0$ to those of $\varepsilon = 1$;
homotopy works robustly if this property holds, and fails if it does not. With
deflated continuation, the requirements for success are weakened. Given a branch
$\{\left(u(t), \lambda(t)\right): t \in [0, 1]\}$, define its support to be $\{\lambda(t): t \in [0, 1]\} \subset \mathbb{R}$. Whereas
switching continuation homotopy finds a solution only if there exists a continuously connected
branch between it and the initial guess, deflated continuation homotopy only necessitates
that the union of the supports of the branches covers the interval $[\lambda_{\textrm{min}}, \lambda_{\textrm{max}}]$.
This is precisely the case in Figure \ref{fig:pendulum}.
As the supports of the branches intersect, Algorithm \ref{alg:mainalg} is able
to discover 
the disconnected branches that come into existence at $\varepsilon
\approx 0.575$ and $\varepsilon \approx 0.697$, and identifies four additional solutions that
switching continuation homotopy neglects along this path.

\subsection{Deformation of a hyperelastic beam} \label{sec:beam_pde}
A major strength of deflated continuation is that it scales to fine
discretizations of partial differential equations (PDEs). Unlike switching
continuation, deflated continuation does not demand the nonscalable computation
of determinants or difficult eigendecompositions to detect bifurcations or switch
branches. In fact, all of the subproblems arising in deflated continuation can
be solved efficiently if a good preconditioner is available for the underlying
forward problem.

The example of section \ref{sec:beam_ode} modelled the deformation of a beam
under compression with Euler's elastica equation. In this example, we model
the same physical phenomenon, but with a two-dimensional compressible
neo-Hookean hyperelastic PDE, solved with scalable Krylov methods and
preconditioners. The potential energy $\Pi$ is given by
\begin{equation} \label{eqn:potential_energy}
\Pi(u) = \int_\Omega \psi(u) \dx - \int_\Omega B\cdot u \dx - \int_{\partial \Omega} T \cdot u \ds,
\end{equation}
where $\Omega$ is the reference domain, $u: \Omega \to \mathbb{R}^2$ is the displacement,
$\psi$ is the elastic stored energy density, $B$ is the body force per unit reference area, 
and $T$ is the traction force per unit reference length. To define $\psi$, consider the
deformation gradient
\begin{equation}
F = I + \nabla u,
\end{equation}
the right Cauchy--Green tensor
\begin{equation}
C = F^T F,
\end{equation}
and its invariants $J = \mathrm{det}(C)$ and $I_c = \mathrm{tr}(C)$. The compressible
neo-Hookean stored energy density is given by
\begin{equation}
\psi = \frac{\mu}{2} (I_c - 2) - \mu \mathrm{log}(J) + \frac{\lambda}{2} \mathrm{log}(J)^2,
\end{equation}
where $\mu$ and $\lambda$ are the Lam\'e parameters, which are calculated from
the Young's modulus $E$ and Poisson ratio $\nu$. In this problem, we take
$\Omega = (0, 1) \times (0, 0.1)$, $B = (0, -1000)$, $T = 0$, $E = 10^6$, and
$\nu = 0.3$. In addition, Dirichlet conditions are imposed on the left and right
boundaries:
\begin{align} \label{eqn:hyperelastic_bcs}
u(0, \cdot) &= (0, 0), \\
u(1, \cdot) &= (0, -\varepsilon),
\end{align}
where $\varepsilon$ is the parameter to be continued. 

For a fixed $\varepsilon$, let $V_\varepsilon = \{u \in
H^1(\Omega; \mathbb{R}^2): u(0, \cdot) = (0, 0), \ u(1, \cdot) = (0,
-\varepsilon)\}$ be the function space of admissible displacements.
Minimizers of \eqref{eqn:potential_energy} are computed by seeking solutions of
the associated optimality condition: find $u \in V_\varepsilon$ such that
\begin{equation} \label{eqn:hyperelasticity}
\Pi'(u; v) = 0 \ \forall\ v \in V_0.
\end{equation}

\begin{figure} 
\centering 
\input{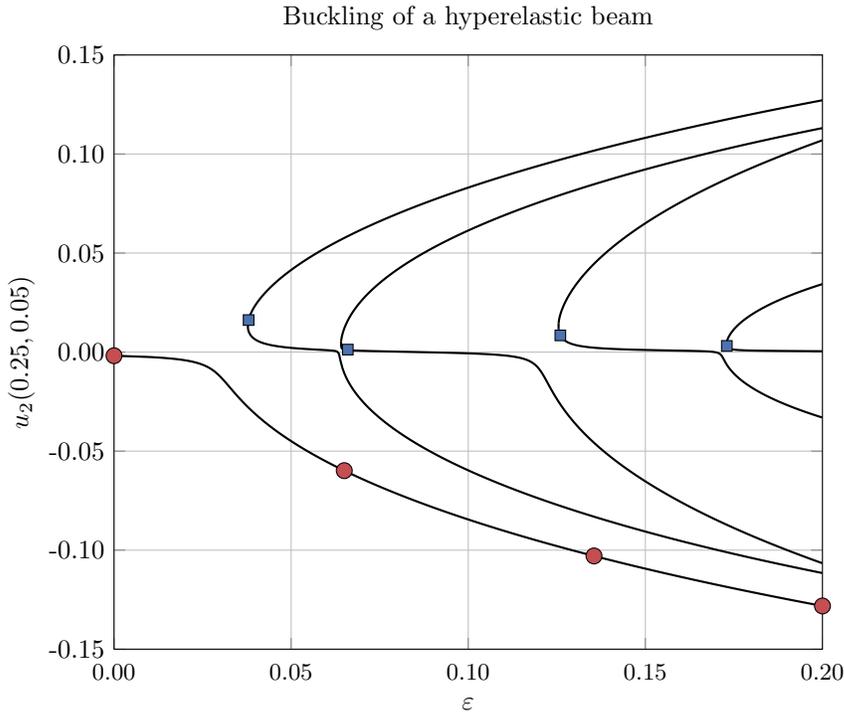} 
\caption{Bifurcation diagram for the hyperelastic PDE formulation for the deformation of a beam
\eqref{eqn:hyperelasticity} as a function of the displacement on the
right-hand boundary. As before, switching
continuation only discovers that part of the diagram labelled with red circles.
Blue squares denote the discovery of disconnected branches with deflation.}
\label{fig:hyperelasticity_diagram}
\end{figure}

\begin{figure}
\centering
\begin{tabular}{cc}
\includegraphics[height=4.9cm]{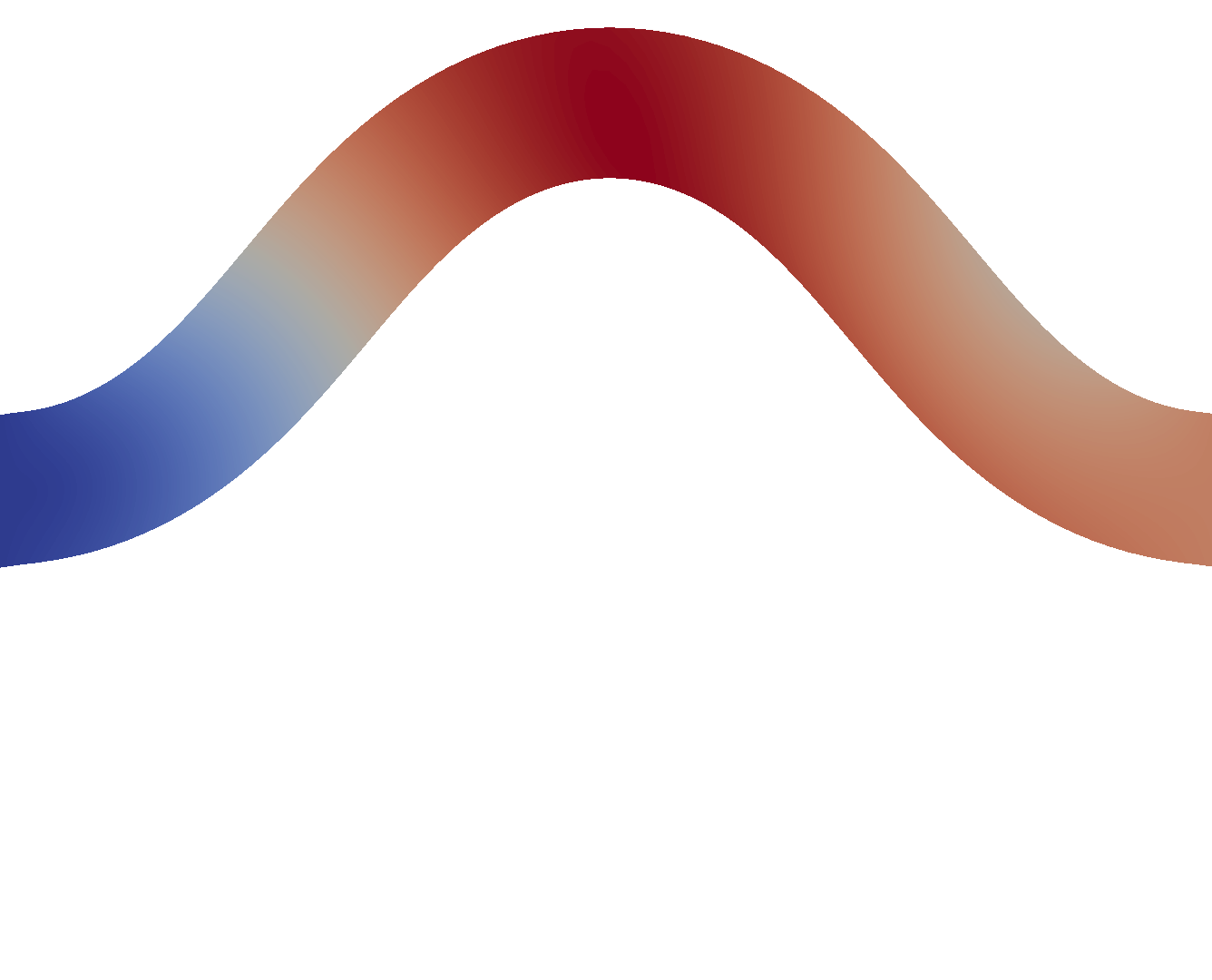} &
\includegraphics[height=4.9cm]{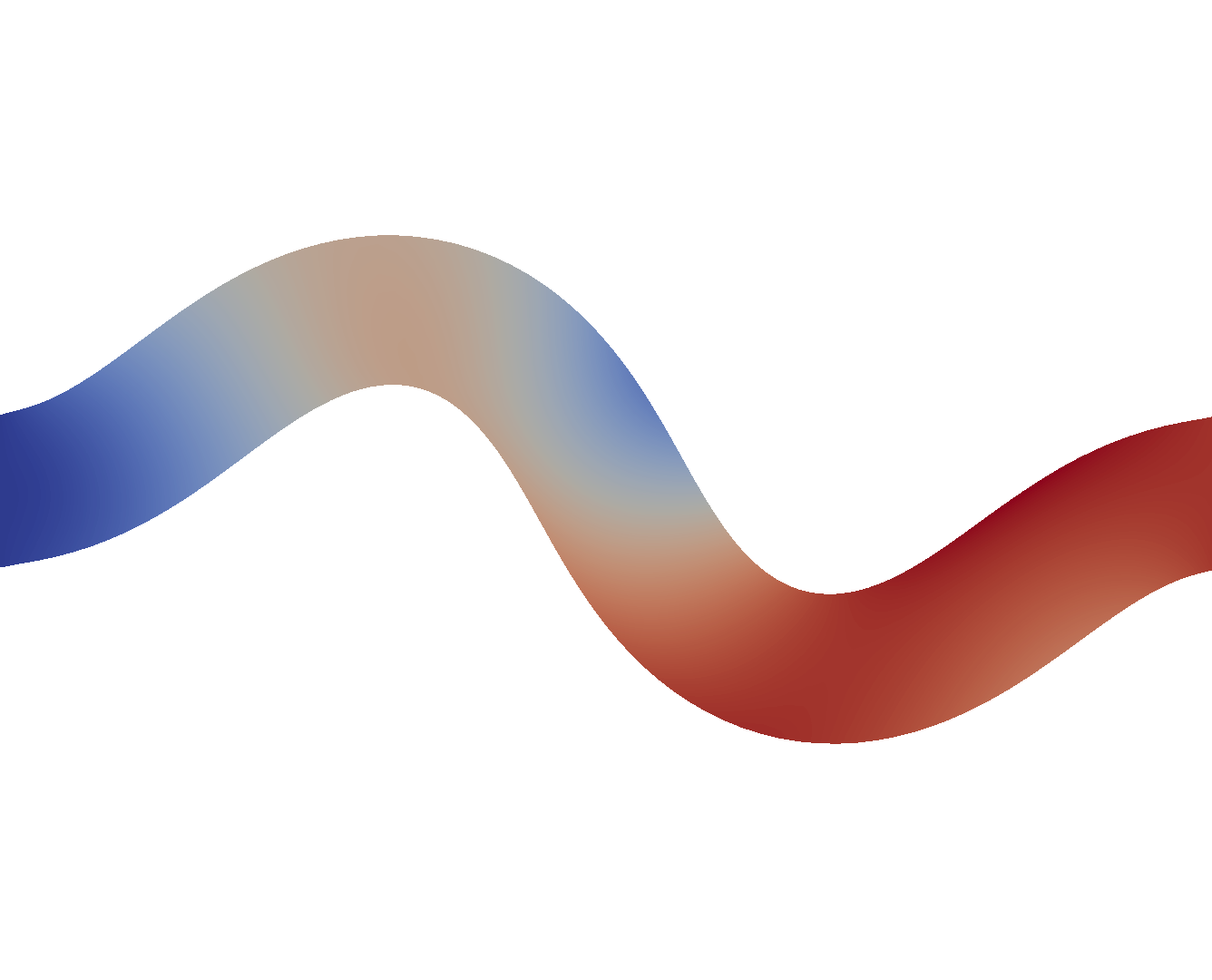} \\
\includegraphics[height=4.9cm]{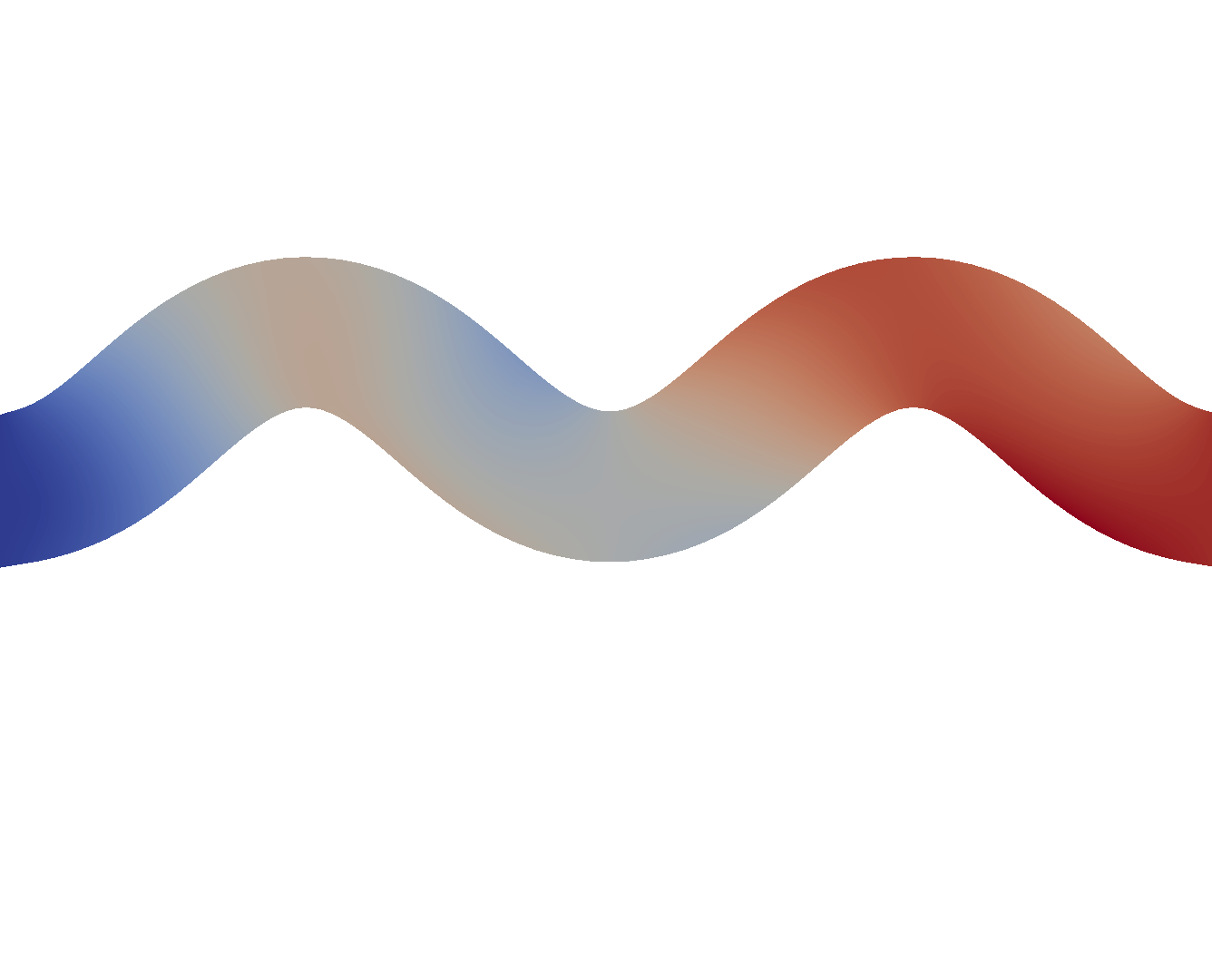} &
\includegraphics[height=4.9cm]{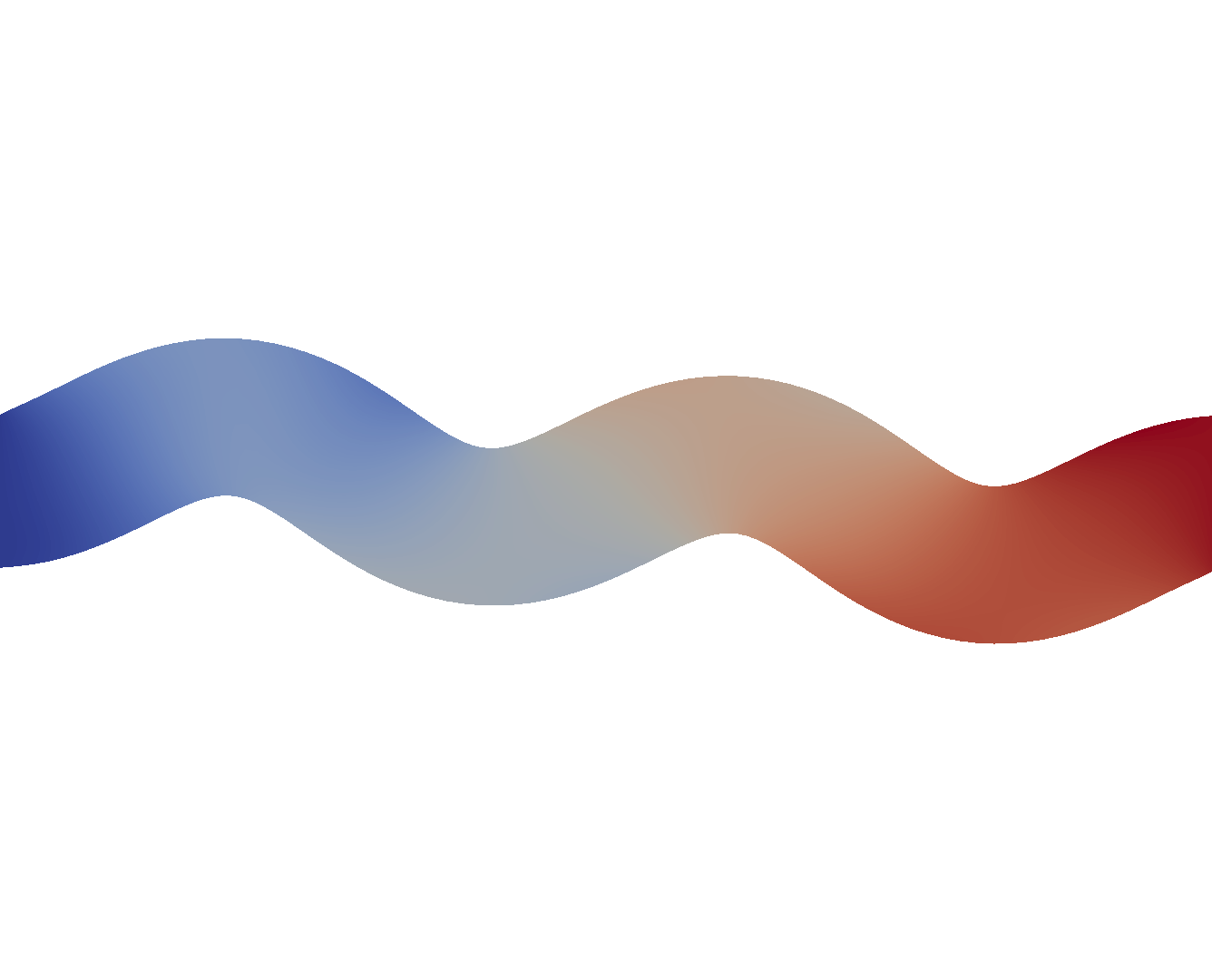} \\
\end{tabular}
\caption{Some of the solutions to the hyperelastic PDE
\eqref{eqn:hyperelasticity} for $\varepsilon = 0.2$, found with deflated
continuation. The color refers to the magnitude of the displacement from the
reference configuration.}
\label{fig:hyperelasticity_solutions}
\end{figure}

Algorithm \ref{alg:mainalg} was applied to \eqref{eqn:hyperelasticity}, from
$\varepsilon = 0$ to $\varepsilon = 0.2$, with continuation step $\Delta
\varepsilon = 0.0005$. The equation was discretized with $2 \times 10^4$
piecewise linear finite elements using FEniCS and PETSc. Newton's method was
terminated with failure if convergence did not occur within $10^2$ iterations.
Each Newton step was solved with the GAMG algebraic multigrid preconditioner
\cite{adams2004}, equipped with the near-nullspace of rigid body modes
\cite{falgout2006}. Deflation was applied with power $p = 2$, shift $\sigma = 1$
and with distances measured in the $H^1$ norm. The functional considered was the
vertical component of displacement evaluated at $(0.25, 0.05)$.

It is well-known that for $B = T = 0$, \eqref{eqn:hyperelasticity} enjoys a
$\mathbb{Z}_2$ reflective symmetry and its bifurcation diagram undergoes a
series of pitchfork bifurcations as $\varepsilon$ is increased, similar to
Figure \ref{fig:euler_sym}. However, in this configuration the reflective
symmetry has been broken by imposing a gravitational body force, causing the
bifurcation diagram to disconnect.  The resulting diagram is shown in Figure
\ref{fig:hyperelasticity_diagram}, and the computed solutions with positive
functional value are shown in Figure \ref{fig:hyperelasticity_solutions}.
The bifurcation diagram has been computed correctly, indicating that
Algorithm \ref{alg:mainalg} is robust to the use of indirect solvers,
and that it will scale to much finer discretizations of PDEs.

\subsection{A generalized Bratu--Gelfand problem in two dimensions}
Previous examples have demonstrated that deflated continuation is able to find
branches that switching continuation misses because they are disconnected.
Switching continuation can fail in other ways: for example, if a bifurcation
is caused by an eigenvalue of even multiplicity crossing the origin, the standard
bifurcation test functional \eqref{eqn:tau_signdet} will neglect it. This example exhibits such a bifurcation,
and demonstrates that deflated continuation is robust to this failure mode.

\begin{figure}
\centering
\input{tikz/mittelmann.tikz}
\caption{Bifurcation diagram for the Mittelmann problem in two dimensions
\eqref{eqn:mittelmann} as a function of $\lambda$. Standard switching continuation
approaches that rely on the sign of the determinant \eqref{eqn:tau_signdet} as a bifurcation test
functional only discover that part of
the diagram labelled with red circles. Blue squares denote the
branch overlooked by switching continuation, but found by deflated continuation.
Compare with \cite[Figure 1]{mittelmann1986}, \cite[Figure 2a]{uecker2014}.}
\label{fig:mittelmann}
\end{figure}

We consider the problem of Mittelmann \cite{mittelmann1986}:
\begin{gather} \label{eqn:mittelmann}
\begin{aligned}
-\nabla^2 y&=-10(y-\lambda e^y) \equiv \phi(y, \lambda), &\quad &\textrm{ in } \Omega=(-0.5,0.5)^2,\\
\nabla y \cdot \hat{n} &=0,&\quad & \textrm{ on } \partial\Omega.
\end{aligned}
\end{gather}
This is a generalization of the Bratu--Gelfand problem
to multiple dimensions with the addition of a linear term, and has
been used as a test problem for the \texttt{PLTMG} \cite{mittelmann1986} and
\texttt{pde2path} \cite{uecker2014} continuation codes.

As noted by Mittelmann, this equation has two spatially constant solutions that satisfy
$\phi(y, \lambda) = 0$, i.e.~$y(x) = \bar{y}$ with $\bar{y} = \lambda e^{\bar{y}}$. Linearising around $\bar{y}$ with 
$y = \bar{y} + w$ yields an eigenvalue problem
\begin{gather}
\begin{aligned}
-\nabla^2 w            &=\left.\frac{\partial \phi}{\partial y}\right\vert_{(\bar{y},\lambda)}w = -10(1 - \bar{y})w ,&\quad &\textrm{ in } \Omega, \\
\nabla w \cdot \hat{n} &=0,                                                                                       &\quad &\textrm{ on } \partial\Omega.
\end{aligned}
\end{gather}
Non-trivial perturbations of the constant solutions can be located by
examining the eigenvalues and corresponding eigenfunctions of the Laplacian on
$\Omega$. The bifurcation points are found by solving $-10(1 - \bar{y}) =
\mu_{m,n}$ where $\mu_{m,n}$ are the eigenvalues of the Laplacian with Neumann
boundary conditions. In this case, $\mu_{m, n}=(m^2+n^2)\pi^2$ for $m, n \in
\mathbb{N}_+,$ and thus the bifurcations occur when $\bar{y}_{m,n} =
1+\mu_{m,n}/10$ and $\lambda_{m,n} = \bar{y}_{m,n} e^{-\bar{y}_{m,n}}$. The
initial bifurcation points occur at $\lambda_{0,0} = e^{-1} \approx 0.3678$ (a
fold bifurcation), $\lambda_{0,1} = \lambda_{1,0} \approx 0.2724$ (a double
pitchfork bifurcation), and $\lambda_{1,1} \approx 0.1519$ (a simple pitchfork bifurcation).

Consider again the computation of the bifurcation test functional $\tau$,
defined in
\eqref{eqn:tau_signdet}.  In the case when a simple bifurcation point is crossed
$\tau$ will indicate this by negation, as one of the eigenvalues
will have passed through the origin. On the other hand, if
an eigenvalue of even multiplicity passes through the origin, $\tau$ remains
unchanged. In this case the bifurcation point is overlooked, and the machinery
of switching continuation is not activated.\footnote{Switching continuation can be
rescued by deliberately breaking the symmetry of the domain, to unfold the
double eigenvalues. Uecker et al.~\cite[Figure 4]{uecker2014} suggest breaking the rotational
symmetry of the domain by solving on $\tilde{\Omega} = (-0.5, 0.5) \times (-0.495, 0.495)$.
The solutions found can then be continued to the true $\Omega$.
While this strategy is successful, it is laborious and requires expert knowledge of
the system at hand.}

As the Mittelmann problem is two-dimensional, the Laplacian has degenerate
eigenvalues of even multiplicity. For example, its eigenvalues $\mu_{0, 1}$ and
$\mu_{1, 0}$ are identical but correspond to different eigenfunctions (related
by rotation). Thus, the associated bifurcation point at $\lambda \approx
0.2724$ is missed by switching continuation, even though the bifurcation diagram
is connected. This deficiency is not specific to this equation, and will manifest in any
situation where such degeneracy of eigenvalues occurs.

By contrast, the specific nature of the bifurcation is irrelevant to deflated
continuation; we expect the algorithm to find nearby solutions regardless of
the details of how the branches are connected (or not connected). To investigate this, Algorithm
\ref{alg:mainalg} was applied to \eqref{eqn:mittelmann}, for $\lambda \in
[0.3678, 0.05]$, with continuation step $\Delta \lambda = -0.0001$. The equation
was discretized with $1600$ piecewise linear finite elements using FEniCS and
PETSc. Newton's method was terminated with failure if convergence did not occur
within $10^2$ iterations. The same deflation and solver settings were used as in
all previous examples. Following Uecker et al.~\cite{uecker2014}, the functional
considered was the $L^2$ norm of the solution. 

The resulting diagram is shown in Figure \ref{fig:mittelmann}. The outer
branches (on the top and bottom) are the constant solutions, with branches
bifurcating from the upper branch at $\lambda_{0,1} = \lambda_{1,0}$ and $\lambda_{1,1}$ as
expected, and secondary bifurcations in turn emanating from these. Importantly,
all branches in this interval have been discovered, including the branch overlooked with
switching continuation (denoted with blue squares). Deflated continuation
applies to both \emph{connected} and disconnected bifurcation diagrams on which
switching continuation fails.

\section{Conclusion}

We have presented a new algorithm for bifurcation analysis that relies on the
elimination of known branches, rather than the detection and analysis of bifurcation
points. In this way, the algorithm applies equally to connected and disconnected
diagrams. We have developed an initial analysis of multiconvergence of Newton's
method, giving sufficient conditions for when convergence to two solutions is
guaranteed. In numerical experiments the algorithm is effective and succeeds where
switching continuation fails.

Unlike switching continuation, the algorithm relies only on the solution of the
original nonlinear problem with a fixed parameter value, and the solution of
deflations of that problem. The latter is straightforward to implement and solve
if a preconditioner for the former is available. There is no need to implement augmented
systems for different kinds of bifurcation points, or to compute expensive test
functionals, or to construct the nullspace of singular operators. Thus, if a
scalable preconditioner for the undeflated Jacobian is available, it will be
possible to apply the algorithm to massive discretizations of PDEs on supercomputers.

\bibliographystyle{siam}
\bibliography{literature}

\end{document}

%% file: tikz/standardbif1.tikz
%
%
\begin{tikzpicture}

\begin{axis}[%
width=4.520833in,
height=3.565625in,
at={(0.758333in,0.48125in)},
scale only axis,
xmin=0.2,
xmax=2,
ymin=-2.5,
ymax=2.5,
yticklabels={,,},
xticklabels={,,},
xtick = {1},
ytick = {0},
xlabel={\Huge $\lambda$},
ylabel={\Huge $u$}
]
\addplot [color=black, mark=*, mark options={solid}, mark size = 5.0pt, forget plot]
  table[row sep=crcr]{%
1.0	0\\
};

\addplot [color=black,solid,line width=2.0pt,forget plot]
  table[row sep=crcr]{%
0	0\\
0.10000001	0\\
0.20000002	0\\
0.30000003	0\\
0.40000004	0\\
0.50000005	0\\
0.60000006	0\\
0.70000007	0\\
0.80000008	0\\
0.90000009	0\\
1.0000001	0\\
1.10000011	0\\
1.20000012	0\\
1.30000013	0\\
1.40000014	0\\
1.50000015	0\\
1.60000016	0\\
1.70000017	0\\
1.80000018	0\\
1.90000019	0\\
2.0000002	0\\
};
\addplot [color=black,solid,line width=2.0pt,forget plot]
  table[row sep=crcr]{%
2.0252982505	2.148789793\\
1.9999999912	2.133854424\\
1.9914705034	2.1287172389\\
1.9666991515	2.1134968024\\
1.9423221535	2.098065624\\
1.918338977	2.0824246622\\
1.8947487901	2.066574951\\
1.8715504715	2.0505175997\\
1.8487425943	2.0342537682\\
1.8263234394	2.0177846474\\
1.8042910396	2.0011115109\\
1.7826431322	1.9842356034\\
1.7613772043	1.9671581895\\
1.7404905213	1.9498805533\\
1.7199801039	1.9324039459\\
1.6998427666	1.9147296013\\
1.6800751397	1.8968587335\\
1.6606736665	1.8787925069\\
1.6416346328	1.8605320417\\
1.6229541861	1.8420784107\\
1.6046283435	1.8234326213\\
1.586653015	1.8045956182\\
1.5690240205	1.7855682795\\
1.5517371024	1.7663514078\\
1.5347879447	1.7469457309\\
1.5181721867	1.7273518992\\
1.5018854386	1.7075704823\\
1.4859232948	1.68760197\\
1.4702813478	1.6674467704\\
1.4549552003	1.6471052109\\
1.4399404779	1.6265775387\\
1.4252328393	1.6058639226\\
1.4108279873	1.5849644539\\
1.396721678	1.5638791486\\
1.3829097315	1.542607953\\
1.3693880373	1.5211507425\\
1.3561525633	1.499507327\\
1.343199366	1.4776774597\\
1.330524591	1.4556608332\\
1.3181244826	1.4334570898\\
1.3059953892	1.4110658269\\
1.2941337664	1.3884866007\\
1.282536183	1.3657189338\\
1.2711993238	1.3427623216\\
1.2601199935	1.3196162397\\
1.2492951202	1.296280151\\
1.2387217572	1.2727535139\\
1.2283970861	1.2490357912\\
1.2183184189	1.2251264593\\
1.2084831991	1.2010250172\\
1.1988890038	1.1767309973\\
1.1895335447	1.1522439758\\
1.1804146687	1.1275635841\\
1.1715303587	1.1026895209\\
1.1628787341	1.0776215644\\
1.1544580507	1.0523595861\\
1.1462667005	1.0269035636\\
1.1383032109	1.0012535956\\
1.1305662446	0.97540991675\\
1.1230545979	0.94937291292\\
1.1157671995	0.92314313797\\
1.1087031084	0.89672132812\\
1.1018615123	0.87010842195\\
1.0952417249	0.84330557578\\
1.0888431827	0.81631418125\\
1.082665442	0.78913588376\\
1.0767081748	0.76177259973\\
1.0709711648	0.73422653418\\
1.0654543028	0.7065001983\\
1.060157581	0.67859642652\\
1.0550810875	0.65051839128\\
1.0502250003	0.62226962042\\
1.0455895801	0.59385401025\\
1.0411751635	0.56527583742\\
1.0369821547	0.53653977171\\
1.0330110179	0.50765088377\\
1.0292622682	0.47861465328\\
1.0257364627	0.44943697182\\
1.0224341913	0.42012414581\\
1.0193560672	0.39068289371\\
1.0165027168	0.36112034313\\
1.0138747702	0.33144402161\\
1.0114728513	0.3016618453\\
1.0092975683	0.27178210372\\
1.0073495042	0.24181344067\\
1.0056292079	0.21176483194\\
1.0041371857	0.18164555854\\
1.0028738935	0.15146517736\\
1.0018397298	0.12123348832\\
1.0010350294	0.090960498516\\
1.0004600583	0.060656384023\\
1.0001150095	0.03033145014\\
1.0000001	0\\
1.0001150095	-0.03033145014\\
1.0004600583	-0.060656384023\\
1.0010350294	-0.090960498516\\
1.0018397298	-0.12123348832\\
1.0028738935	-0.15146517736\\
1.0041371857	-0.18164555854\\
1.0056292079	-0.21176483194\\
1.0073495042	-0.24181344067\\
1.0092975683	-0.27178210372\\
1.0114728513	-0.3016618453\\
1.0138747702	-0.33144402161\\
1.0165027168	-0.36112034313\\
1.0193560672	-0.39068289371\\
1.0224341913	-0.42012414581\\
1.0257364627	-0.44943697182\\
1.0292622682	-0.47861465328\\
1.0330110179	-0.50765088377\\
1.0369821547	-0.53653977171\\
1.0411751635	-0.56527583742\\
1.0455895801	-0.59385401025\\
1.0502250003	-0.62226962042\\
1.0550810875	-0.65051839128\\
1.060157581	-0.67859642652\\
1.0654543028	-0.7065001983\\
1.0709711648	-0.73422653418\\
1.0767081748	-0.76177259973\\
1.082665442	-0.78913588376\\
1.0888431827	-0.81631418125\\
1.0952417249	-0.84330557578\\
1.1018615123	-0.87010842195\\
1.1087031084	-0.89672132812\\
1.1157671995	-0.92314313797\\
1.1230545979	-0.94937291292\\
1.1305662446	-0.97540991675\\
1.1383032109	-1.0012535956\\
1.1462667005	-1.0269035636\\
1.1544580507	-1.0523595861\\
1.1628787341	-1.0776215644\\
1.1715303587	-1.1026895209\\
1.1804146687	-1.1275635841\\
1.1895335447	-1.1522439758\\
1.1988890038	-1.1767309973\\
1.2084831991	-1.2010250172\\
1.2183184189	-1.2251264593\\
1.2283970861	-1.2490357912\\
1.2387217572	-1.2727535139\\
1.2492951202	-1.296280151\\
1.2601199935	-1.3196162397\\
1.2711993238	-1.3427623216\\
1.282536183	-1.3657189338\\
1.2941337664	-1.3884866007\\
1.3059953892	-1.4110658269\\
1.3181244826	-1.4334570898\\
1.330524591	-1.4556608332\\
1.343199366	-1.4776774597\\
1.3561525633	-1.499507327\\
1.3693880373	-1.5211507425\\
1.3829097315	-1.542607953\\
1.396721678	-1.5638791486\\
1.4108279873	-1.5849644539\\
1.4252328393	-1.6058639226\\
1.4399404779	-1.6265775387\\
1.4549552003	-1.6471052109\\
1.4702813478	-1.6674467704\\
1.4859232948	-1.68760197\\
1.5018854386	-1.7075704823\\
1.5181721867	-1.7273518992\\
1.5347879447	-1.7469457309\\
1.5517371024	-1.7663514078\\
1.5690240205	-1.7855682795\\
1.586653015	-1.8045956182\\
1.6046283435	-1.8234326213\\
1.6229541861	-1.8420784107\\
1.6416346328	-1.8605320417\\
1.6606736665	-1.8787925069\\
1.6800751397	-1.8968587335\\
1.6998427666	-1.9147296013\\
1.7199801039	-1.9324039459\\
1.7404905213	-1.9498805533\\
1.7613772043	-1.9671581895\\
1.7826431322	-1.9842356034\\
1.8042910396	-2.0011115109\\
1.8263234394	-2.0177846474\\
1.8487425943	-2.0342537682\\
1.8715504715	-2.0505175997\\
1.8947487901	-2.066574951\\
1.918338977	-2.0824246622\\
1.9423221535	-2.098065624\\
1.9666991515	-2.1134968024\\
1.9914705034	-2.1287172389\\
1.9999999912	-2.133854424\\
2.0252982505	-2.148789793\\
};
\addplot [->,color=red,solid,line width=4.0pt,forget plot]
  table[row sep=crcr]{%
0 0\\
1.5 0\\
};
\end{axis}

\end{tikzpicture}%

%% file: tikz/standardbif2.tikz
%
%
\begin{tikzpicture}

\begin{axis}[%
width=4.520833in,
height=3.565625in,
at={(0.758333in,0.48125in)},
scale only axis,
xmin=0.2,
xmax=2,
ymin=-2.5,
ymax=2.5,
yticklabels={,,},
xticklabels={,,},
xtick = {1},
ytick = {0},
xlabel={\Huge $\lambda$},
ylabel={\Huge $u$}
]
\addplot [color=black, mark=*, mark options={solid}, mark size = 5.0pt, forget plot]
  table[row sep=crcr]{%
1.0	0\\
};

\addplot [color=black,solid,line width=2.0pt,forget plot]
  table[row sep=crcr]{%
0	0\\
0.10000001	0\\
0.20000002	0\\
0.30000003	0\\
0.40000004	0\\
0.50000005	0\\
0.60000006	0\\
0.70000007	0\\
0.80000008	0\\
0.90000009	0\\
1.0000001	0\\
1.10000011	0\\
1.20000012	0\\
1.30000013	0\\
1.40000014	0\\
1.50000015	0\\
1.60000016	0\\
1.70000017	0\\
1.80000018	0\\
1.90000019	0\\
2.0000002	0\\
};
\addplot [color=black,solid,line width=2.0pt,forget plot]
  table[row sep=crcr]{%
2.0252982505	2.148789793\\
1.9999999912	2.133854424\\
1.9914705034	2.1287172389\\
1.9666991515	2.1134968024\\
1.9423221535	2.098065624\\
1.918338977	2.0824246622\\
1.8947487901	2.066574951\\
1.8715504715	2.0505175997\\
1.8487425943	2.0342537682\\
1.8263234394	2.0177846474\\
1.8042910396	2.0011115109\\
1.7826431322	1.9842356034\\
1.7613772043	1.9671581895\\
1.7404905213	1.9498805533\\
1.7199801039	1.9324039459\\
1.6998427666	1.9147296013\\
1.6800751397	1.8968587335\\
1.6606736665	1.8787925069\\
1.6416346328	1.8605320417\\
1.6229541861	1.8420784107\\
1.6046283435	1.8234326213\\
1.586653015	1.8045956182\\
1.5690240205	1.7855682795\\
1.5517371024	1.7663514078\\
1.5347879447	1.7469457309\\
1.5181721867	1.7273518992\\
1.5018854386	1.7075704823\\
1.4859232948	1.68760197\\
1.4702813478	1.6674467704\\
1.4549552003	1.6471052109\\
1.4399404779	1.6265775387\\
1.4252328393	1.6058639226\\
1.4108279873	1.5849644539\\
1.396721678	1.5638791486\\
1.3829097315	1.542607953\\
1.3693880373	1.5211507425\\
1.3561525633	1.499507327\\
1.343199366	1.4776774597\\
1.330524591	1.4556608332\\
1.3181244826	1.4334570898\\
1.3059953892	1.4110658269\\
1.2941337664	1.3884866007\\
1.282536183	1.3657189338\\
1.2711993238	1.3427623216\\
1.2601199935	1.3196162397\\
1.2492951202	1.296280151\\
1.2387217572	1.2727535139\\
1.2283970861	1.2490357912\\
1.2183184189	1.2251264593\\
1.2084831991	1.2010250172\\
1.1988890038	1.1767309973\\
1.1895335447	1.1522439758\\
1.1804146687	1.1275635841\\
1.1715303587	1.1026895209\\
1.1628787341	1.0776215644\\
1.1544580507	1.0523595861\\
1.1462667005	1.0269035636\\
1.1383032109	1.0012535956\\
1.1305662446	0.97540991675\\
1.1230545979	0.94937291292\\
1.1157671995	0.92314313797\\
1.1087031084	0.89672132812\\
1.1018615123	0.87010842195\\
1.0952417249	0.84330557578\\
1.0888431827	0.81631418125\\
1.082665442	0.78913588376\\
1.0767081748	0.76177259973\\
1.0709711648	0.73422653418\\
1.0654543028	0.7065001983\\
1.060157581	0.67859642652\\
1.0550810875	0.65051839128\\
1.0502250003	0.62226962042\\
1.0455895801	0.59385401025\\
1.0411751635	0.56527583742\\
1.0369821547	0.53653977171\\
1.0330110179	0.50765088377\\
1.0292622682	0.47861465328\\
1.0257364627	0.44943697182\\
1.0224341913	0.42012414581\\
1.0193560672	0.39068289371\\
1.0165027168	0.36112034313\\
1.0138747702	0.33144402161\\
1.0114728513	0.3016618453\\
1.0092975683	0.27178210372\\
1.0073495042	0.24181344067\\
1.0056292079	0.21176483194\\
1.0041371857	0.18164555854\\
1.0028738935	0.15146517736\\
1.0018397298	0.12123348832\\
1.0010350294	0.090960498516\\
1.0004600583	0.060656384023\\
1.0001150095	0.03033145014\\
1.0000001	0\\
1.0001150095	-0.03033145014\\
1.0004600583	-0.060656384023\\
1.0010350294	-0.090960498516\\
1.0018397298	-0.12123348832\\
1.0028738935	-0.15146517736\\
1.0041371857	-0.18164555854\\
1.0056292079	-0.21176483194\\
1.0073495042	-0.24181344067\\
1.0092975683	-0.27178210372\\
1.0114728513	-0.3016618453\\
1.0138747702	-0.33144402161\\
1.0165027168	-0.36112034313\\
1.0193560672	-0.39068289371\\
1.0224341913	-0.42012414581\\
1.0257364627	-0.44943697182\\
1.0292622682	-0.47861465328\\
1.0330110179	-0.50765088377\\
1.0369821547	-0.53653977171\\
1.0411751635	-0.56527583742\\
1.0455895801	-0.59385401025\\
1.0502250003	-0.62226962042\\
1.0550810875	-0.65051839128\\
1.060157581	-0.67859642652\\
1.0654543028	-0.7065001983\\
1.0709711648	-0.73422653418\\
1.0767081748	-0.76177259973\\
1.082665442	-0.78913588376\\
1.0888431827	-0.81631418125\\
1.0952417249	-0.84330557578\\
1.1018615123	-0.87010842195\\
1.1087031084	-0.89672132812\\
1.1157671995	-0.92314313797\\
1.1230545979	-0.94937291292\\
1.1305662446	-0.97540991675\\
1.1383032109	-1.0012535956\\
1.1462667005	-1.0269035636\\
1.1544580507	-1.0523595861\\
1.1628787341	-1.0776215644\\
1.1715303587	-1.1026895209\\
1.1804146687	-1.1275635841\\
1.1895335447	-1.1522439758\\
1.1988890038	-1.1767309973\\
1.2084831991	-1.2010250172\\
1.2183184189	-1.2251264593\\
1.2283970861	-1.2490357912\\
1.2387217572	-1.2727535139\\
1.2492951202	-1.296280151\\
1.2601199935	-1.3196162397\\
1.2711993238	-1.3427623216\\
1.282536183	-1.3657189338\\
1.2941337664	-1.3884866007\\
1.3059953892	-1.4110658269\\
1.3181244826	-1.4334570898\\
1.330524591	-1.4556608332\\
1.343199366	-1.4776774597\\
1.3561525633	-1.499507327\\
1.3693880373	-1.5211507425\\
1.3829097315	-1.542607953\\
1.396721678	-1.5638791486\\
1.4108279873	-1.5849644539\\
1.4252328393	-1.6058639226\\
1.4399404779	-1.6265775387\\
1.4549552003	-1.6471052109\\
1.4702813478	-1.6674467704\\
1.4859232948	-1.68760197\\
1.5018854386	-1.7075704823\\
1.5181721867	-1.7273518992\\
1.5347879447	-1.7469457309\\
1.5517371024	-1.7663514078\\
1.5690240205	-1.7855682795\\
1.586653015	-1.8045956182\\
1.6046283435	-1.8234326213\\
1.6229541861	-1.8420784107\\
1.6416346328	-1.8605320417\\
1.6606736665	-1.8787925069\\
1.6800751397	-1.8968587335\\
1.6998427666	-1.9147296013\\
1.7199801039	-1.9324039459\\
1.7404905213	-1.9498805533\\
1.7613772043	-1.9671581895\\
1.7826431322	-1.9842356034\\
1.8042910396	-2.0011115109\\
1.8263234394	-2.0177846474\\
1.8487425943	-2.0342537682\\
1.8715504715	-2.0505175997\\
1.8947487901	-2.066574951\\
1.918338977	-2.0824246622\\
1.9423221535	-2.098065624\\
1.9666991515	-2.1134968024\\
1.9914705034	-2.1287172389\\
1.9999999912	-2.133854424\\
2.0252982505	-2.148789793\\
};
\addplot [color=red,solid,line width=4.0pt,forget plot]
  table[row sep=crcr]{%
0 0\\
1.5 0\\
};
\addplot [<-,color=red,solid,line width=4.0pt,forget plot]
  table[row sep=crcr]{%
1.1087031084	0.89672132812\\
1.1018615123	0.87010842195\\
1.0952417249	0.84330557578\\
1.0888431827	0.81631418125\\
1.082665442	0.78913588376\\
1.0767081748	0.76177259973\\
1.0709711648	0.73422653418\\
1.0654543028	0.7065001983\\
1.060157581	0.67859642652\\
1.0550810875	0.65051839128\\
1.0502250003	0.62226962042\\
1.0455895801	0.59385401025\\
1.0411751635	0.56527583742\\
1.0369821547	0.53653977171\\
1.0330110179	0.50765088377\\
1.0292622682	0.47861465328\\
1.0257364627	0.44943697182\\
1.0224341913	0.42012414581\\
1.0193560672	0.39068289371\\
1.0165027168	0.36112034313\\
1.0138747702	0.33144402161\\
1.0114728513	0.3016618453\\
1.0092975683	0.27178210372\\
1.0073495042	0.24181344067\\
1.0056292079	0.21176483194\\
1.0041371857	0.18164555854\\
1.0028738935	0.15146517736\\
1.0018397298	0.12123348832\\
1.0010350294	0.090960498516\\
1.0004600583	0.060656384023\\
1.0001150095	0.03033145014\\
1.0000001	0\\
};
\addplot [<-,color=red,solid,line width=4.0pt,forget plot]
  table[row sep=crcr]{%
1.1087031084	-0.89672132812\\
1.1018615123	-0.87010842195\\
1.0952417249	-0.84330557578\\
1.0888431827	-0.81631418125\\
1.082665442	-0.78913588376\\
1.0767081748	-0.76177259973\\
1.0709711648	-0.73422653418\\
1.0654543028	-0.7065001983\\
1.060157581	-0.67859642652\\
1.0550810875	-0.65051839128\\
1.0502250003	-0.62226962042\\
1.0455895801	-0.59385401025\\
1.0411751635	-0.56527583742\\
1.0369821547	-0.53653977171\\
1.0330110179	-0.50765088377\\
1.0292622682	-0.47861465328\\
1.0257364627	-0.44943697182\\
1.0224341913	-0.42012414581\\
1.0193560672	-0.39068289371\\
1.0165027168	-0.36112034313\\
1.0138747702	-0.33144402161\\
1.0114728513	-0.3016618453\\
1.0092975683	-0.27178210372\\
1.0073495042	-0.24181344067\\
1.0056292079	-0.21176483194\\
1.0041371857	-0.18164555854\\
1.0028738935	-0.15146517736\\
1.0018397298	-0.12123348832\\
1.0010350294	-0.090960498516\\
1.0004600583	-0.060656384023\\
1.0001150095	-0.03033145014\\
1.0000001	-0\\
};\end{axis}
\end{tikzpicture}%

%% file: tikz/standardbif3.tikz
%
%
\begin{tikzpicture}

\begin{axis}[%
width=4.520833in,
height=3.565625in,
at={(0.758333in,0.48125in)},
scale only axis,
xmin=0.2,
xmax=2,
ymin=-2.5,
ymax=2.5,
yticklabels={,,},
xticklabels={,,},
xtick = {1},
ytick = {0},
xlabel={\Huge $\lambda$},
ylabel={\Huge $u$}
]
\addplot [color=black, mark=*, mark options={solid}, mark size = 5.0pt, forget plot]
  table[row sep=crcr]{%
1.0	0\\
};

\addplot [color=black,solid,line width=2.0pt,forget plot]
  table[row sep=crcr]{%
0	0\\
0.10000001	0\\
0.20000002	0\\
0.30000003	0\\
0.40000004	0\\
0.50000005	0\\
0.60000006	0\\
0.70000007	0\\
0.80000008	0\\
0.90000009	0\\
1.0000001	0\\
1.10000011	0\\
1.20000012	0\\
1.30000013	0\\
1.40000014	0\\
1.50000015	0\\
1.60000016	0\\
1.70000017	0\\
1.80000018	0\\
1.90000019	0\\
2.0000002	0\\
};
\addplot [color=black,solid,line width=2.0pt,forget plot]
  table[row sep=crcr]{%
2.0252982505	2.148789793\\
1.9999999912	2.133854424\\
1.9914705034	2.1287172389\\
1.9666991515	2.1134968024\\
1.9423221535	2.098065624\\
1.918338977	2.0824246622\\
1.8947487901	2.066574951\\
1.8715504715	2.0505175997\\
1.8487425943	2.0342537682\\
1.8263234394	2.0177846474\\
1.8042910396	2.0011115109\\
1.7826431322	1.9842356034\\
1.7613772043	1.9671581895\\
1.7404905213	1.9498805533\\
1.7199801039	1.9324039459\\
1.6998427666	1.9147296013\\
1.6800751397	1.8968587335\\
1.6606736665	1.8787925069\\
1.6416346328	1.8605320417\\
1.6229541861	1.8420784107\\
1.6046283435	1.8234326213\\
1.586653015	1.8045956182\\
1.5690240205	1.7855682795\\
1.5517371024	1.7663514078\\
1.5347879447	1.7469457309\\
1.5181721867	1.7273518992\\
1.5018854386	1.7075704823\\
1.4859232948	1.68760197\\
1.4702813478	1.6674467704\\
1.4549552003	1.6471052109\\
1.4399404779	1.6265775387\\
1.4252328393	1.6058639226\\
1.4108279873	1.5849644539\\
1.396721678	1.5638791486\\
1.3829097315	1.542607953\\
1.3693880373	1.5211507425\\
1.3561525633	1.499507327\\
1.343199366	1.4776774597\\
1.330524591	1.4556608332\\
1.3181244826	1.4334570898\\
1.3059953892	1.4110658269\\
1.2941337664	1.3884866007\\
1.282536183	1.3657189338\\
1.2711993238	1.3427623216\\
1.2601199935	1.3196162397\\
1.2492951202	1.296280151\\
1.2387217572	1.2727535139\\
1.2283970861	1.2490357912\\
1.2183184189	1.2251264593\\
1.2084831991	1.2010250172\\
1.1988890038	1.1767309973\\
1.1895335447	1.1522439758\\
1.1804146687	1.1275635841\\
1.1715303587	1.1026895209\\
1.1628787341	1.0776215644\\
1.1544580507	1.0523595861\\
1.1462667005	1.0269035636\\
1.1383032109	1.0012535956\\
1.1305662446	0.97540991675\\
1.1230545979	0.94937291292\\
1.1157671995	0.92314313797\\
1.1087031084	0.89672132812\\
1.1018615123	0.87010842195\\
1.0952417249	0.84330557578\\
1.0888431827	0.81631418125\\
1.082665442	0.78913588376\\
1.0767081748	0.76177259973\\
1.0709711648	0.73422653418\\
1.0654543028	0.7065001983\\
1.060157581	0.67859642652\\
1.0550810875	0.65051839128\\
1.0502250003	0.62226962042\\
1.0455895801	0.59385401025\\
1.0411751635	0.56527583742\\
1.0369821547	0.53653977171\\
1.0330110179	0.50765088377\\
1.0292622682	0.47861465328\\
1.0257364627	0.44943697182\\
1.0224341913	0.42012414581\\
1.0193560672	0.39068289371\\
1.0165027168	0.36112034313\\
1.0138747702	0.33144402161\\
1.0114728513	0.3016618453\\
1.0092975683	0.27178210372\\
1.0073495042	0.24181344067\\
1.0056292079	0.21176483194\\
1.0041371857	0.18164555854\\
1.0028738935	0.15146517736\\
1.0018397298	0.12123348832\\
1.0010350294	0.090960498516\\
1.0004600583	0.060656384023\\
1.0001150095	0.03033145014\\
1.0000001	0\\
1.0001150095	-0.03033145014\\
1.0004600583	-0.060656384023\\
1.0010350294	-0.090960498516\\
1.0018397298	-0.12123348832\\
1.0028738935	-0.15146517736\\
1.0041371857	-0.18164555854\\
1.0056292079	-0.21176483194\\
1.0073495042	-0.24181344067\\
1.0092975683	-0.27178210372\\
1.0114728513	-0.3016618453\\
1.0138747702	-0.33144402161\\
1.0165027168	-0.36112034313\\
1.0193560672	-0.39068289371\\
1.0224341913	-0.42012414581\\
1.0257364627	-0.44943697182\\
1.0292622682	-0.47861465328\\
1.0330110179	-0.50765088377\\
1.0369821547	-0.53653977171\\
1.0411751635	-0.56527583742\\
1.0455895801	-0.59385401025\\
1.0502250003	-0.62226962042\\
1.0550810875	-0.65051839128\\
1.060157581	-0.67859642652\\
1.0654543028	-0.7065001983\\
1.0709711648	-0.73422653418\\
1.0767081748	-0.76177259973\\
1.082665442	-0.78913588376\\
1.0888431827	-0.81631418125\\
1.0952417249	-0.84330557578\\
1.1018615123	-0.87010842195\\
1.1087031084	-0.89672132812\\
1.1157671995	-0.92314313797\\
1.1230545979	-0.94937291292\\
1.1305662446	-0.97540991675\\
1.1383032109	-1.0012535956\\
1.1462667005	-1.0269035636\\
1.1544580507	-1.0523595861\\
1.1628787341	-1.0776215644\\
1.1715303587	-1.1026895209\\
1.1804146687	-1.1275635841\\
1.1895335447	-1.1522439758\\
1.1988890038	-1.1767309973\\
1.2084831991	-1.2010250172\\
1.2183184189	-1.2251264593\\
1.2283970861	-1.2490357912\\
1.2387217572	-1.2727535139\\
1.2492951202	-1.296280151\\
1.2601199935	-1.3196162397\\
1.2711993238	-1.3427623216\\
1.282536183	-1.3657189338\\
1.2941337664	-1.3884866007\\
1.3059953892	-1.4110658269\\
1.3181244826	-1.4334570898\\
1.330524591	-1.4556608332\\
1.343199366	-1.4776774597\\
1.3561525633	-1.499507327\\
1.3693880373	-1.5211507425\\
1.3829097315	-1.542607953\\
1.396721678	-1.5638791486\\
1.4108279873	-1.5849644539\\
1.4252328393	-1.6058639226\\
1.4399404779	-1.6265775387\\
1.4549552003	-1.6471052109\\
1.4702813478	-1.6674467704\\
1.4859232948	-1.68760197\\
1.5018854386	-1.7075704823\\
1.5181721867	-1.7273518992\\
1.5347879447	-1.7469457309\\
1.5517371024	-1.7663514078\\
1.5690240205	-1.7855682795\\
1.586653015	-1.8045956182\\
1.6046283435	-1.8234326213\\
1.6229541861	-1.8420784107\\
1.6416346328	-1.8605320417\\
1.6606736665	-1.8787925069\\
1.6800751397	-1.8968587335\\
1.6998427666	-1.9147296013\\
1.7199801039	-1.9324039459\\
1.7404905213	-1.9498805533\\
1.7613772043	-1.9671581895\\
1.7826431322	-1.9842356034\\
1.8042910396	-2.0011115109\\
1.8263234394	-2.0177846474\\
1.8487425943	-2.0342537682\\
1.8715504715	-2.0505175997\\
1.8947487901	-2.066574951\\
1.918338977	-2.0824246622\\
1.9423221535	-2.098065624\\
1.9666991515	-2.1134968024\\
1.9914705034	-2.1287172389\\
1.9999999912	-2.133854424\\
2.0252982505	-2.148789793\\
};
\addplot [->,color=red,solid,line width=4.0pt,forget plot]
  table[row sep=crcr]{%
0 0\\
1.7 0\\
};
\addplot [<-,color=red,solid,line width=4.0pt,forget plot]
  table[row sep=crcr]{%
1.6998427666	1.9147296013\\
1.6800751397	1.8968587335\\
1.6606736665	1.8787925069\\
1.6416346328	1.8605320417\\
1.6229541861	1.8420784107\\
1.6046283435	1.8234326213\\
1.586653015	1.8045956182\\
1.5690240205	1.7855682795\\
1.5517371024	1.7663514078\\
1.5347879447	1.7469457309\\
1.5181721867	1.7273518992\\
1.5018854386	1.7075704823\\
1.4859232948	1.68760197\\
1.4702813478	1.6674467704\\
1.4549552003	1.6471052109\\
1.4399404779	1.6265775387\\
1.4252328393	1.6058639226\\
1.4108279873	1.5849644539\\
1.396721678	1.5638791486\\
1.3829097315	1.542607953\\
1.3693880373	1.5211507425\\
1.3561525633	1.499507327\\
1.343199366	1.4776774597\\
1.330524591	1.4556608332\\
1.3181244826	1.4334570898\\
1.3059953892	1.4110658269\\
1.2941337664	1.3884866007\\
1.282536183	1.3657189338\\
1.2711993238	1.3427623216\\
1.2601199935	1.3196162397\\
1.2492951202	1.296280151\\
1.2387217572	1.2727535139\\
1.2283970861	1.2490357912\\
1.2183184189	1.2251264593\\
1.2084831991	1.2010250172\\
1.1988890038	1.1767309973\\
1.1895335447	1.1522439758\\
1.1804146687	1.1275635841\\
1.1715303587	1.1026895209\\
1.1628787341	1.0776215644\\
1.1544580507	1.0523595861\\
1.1462667005	1.0269035636\\
1.1383032109	1.0012535956\\
1.1305662446	0.97540991675\\
1.1230545979	0.94937291292\\
1.1157671995	0.92314313797\\
1.1087031084	0.89672132812\\
1.1018615123	0.87010842195\\
1.0952417249	0.84330557578\\
1.0888431827	0.81631418125\\
1.082665442	0.78913588376\\
1.0767081748	0.76177259973\\
1.0709711648	0.73422653418\\
1.0654543028	0.7065001983\\
1.060157581	0.67859642652\\
1.0550810875	0.65051839128\\
1.0502250003	0.62226962042\\
1.0455895801	0.59385401025\\
1.0411751635	0.56527583742\\
1.0369821547	0.53653977171\\
1.0330110179	0.50765088377\\
1.0292622682	0.47861465328\\
1.0257364627	0.44943697182\\
1.0224341913	0.42012414581\\
1.0193560672	0.39068289371\\
1.0165027168	0.36112034313\\
1.0138747702	0.33144402161\\
1.0114728513	0.3016618453\\
1.0092975683	0.27178210372\\
1.0073495042	0.24181344067\\
1.0056292079	0.21176483194\\
1.0041371857	0.18164555854\\
1.0028738935	0.15146517736\\
1.0018397298	0.12123348832\\
1.0010350294	0.090960498516\\
1.0004600583	0.060656384023\\
1.0001150095	0.03033145014\\
1.0000001	0\\
};
\addplot [<-,color=red,solid,line width=4.0pt,forget plot]
  table[row sep=crcr]{%
1.6998427666	-1.9147296013\\
1.6800751397	-1.8968587335\\
1.6606736665	-1.8787925069\\
1.6416346328	-1.8605320417\\
1.6229541861	-1.8420784107\\
1.6046283435	-1.8234326213\\
1.586653015	-1.8045956182\\
1.5690240205	-1.7855682795\\
1.5517371024	-1.7663514078\\
1.5347879447	-1.7469457309\\
1.5181721867	-1.7273518992\\
1.5018854386	-1.7075704823\\
1.4859232948	-1.68760197\\
1.4702813478	-1.6674467704\\
1.4549552003	-1.6471052109\\
1.4399404779	-1.6265775387\\
1.4252328393	-1.6058639226\\
1.4108279873	-1.5849644539\\
1.396721678	-1.5638791486\\
1.3829097315	-1.542607953\\
1.3693880373	-1.5211507425\\
1.3561525633	-1.499507327\\
1.343199366	-1.4776774597\\
1.330524591	-1.4556608332\\
1.3181244826	-1.4334570898\\
1.3059953892	-1.4110658269\\
1.2941337664	-1.3884866007\\
1.282536183	-1.3657189338\\
1.2711993238	-1.3427623216\\
1.2601199935	-1.3196162397\\
1.2492951202	-1.296280151\\
1.2387217572	-1.2727535139\\
1.2283970861	-1.2490357912\\
1.2183184189	-1.2251264593\\
1.2084831991	-1.2010250172\\
1.1988890038	-1.1767309973\\
1.1895335447	-1.1522439758\\
1.1804146687	-1.1275635841\\
1.1715303587	-1.1026895209\\
1.1628787341	-1.0776215644\\
1.1544580507	-1.0523595861\\
1.1462667005	-1.0269035636\\
1.1383032109	-1.0012535956\\
1.1305662446	-0.97540991675\\
1.1230545979	-0.94937291292\\
1.1157671995	-0.92314313797\\
1.1087031084	-0.89672132812\\
1.1018615123	-0.87010842195\\
1.0952417249	-0.84330557578\\
1.0888431827	-0.81631418125\\
1.082665442	-0.78913588376\\
1.0767081748	-0.76177259973\\
1.0709711648	-0.73422653418\\
1.0654543028	-0.7065001983\\
1.060157581	-0.67859642652\\
1.0550810875	-0.65051839128\\
1.0502250003	-0.62226962042\\
1.0455895801	-0.59385401025\\
1.0411751635	-0.56527583742\\
1.0369821547	-0.53653977171\\
1.0330110179	-0.50765088377\\
1.0292622682	-0.47861465328\\
1.0257364627	-0.44943697182\\
1.0224341913	-0.42012414581\\
1.0193560672	-0.39068289371\\
1.0165027168	-0.36112034313\\
1.0138747702	-0.33144402161\\
1.0114728513	-0.3016618453\\
1.0092975683	-0.27178210372\\
1.0073495042	-0.24181344067\\
1.0056292079	-0.21176483194\\
1.0041371857	-0.18164555854\\
1.0028738935	-0.15146517736\\
1.0018397298	-0.12123348832\\
1.0010350294	-0.090960498516\\
1.0004600583	-0.060656384023\\
1.0001150095	-0.03033145014\\
1.0000001	0\\
};
\end{axis}

\end{tikzpicture}%

%% file: tikz/deflationbif1.tikz
%
%
\begin{tikzpicture}

\begin{axis}[%
width=4.520833in,
height=3.565625in,
at={(0.758333in,0.48125in)},
scale only axis,
xmin=0.2,
xmax=2,
ymin=-2.5,
ymax=2.5,
yticklabels={,,},
xticklabels={,,},
xtick = {0},
ytick = {0},
xlabel={\Huge $\lambda$},
ylabel={\Huge $u$}
]

\addplot [color=black,solid,line width=2.0pt,forget plot]
  table[row sep=crcr]{%
2.0252982505	-2.1511656348\\
1.9999309186	-2.1362543389\\
1.9749563582	-2.1211296905\\
1.9503766066	-2.1057937463\\
1.926191269	-2.0902474329\\
1.9023996495	-2.0744917679\\
1.8790007424	-2.0585278475\\
1.8559932461	-2.0423568163\\
1.8333755614	-2.0259798846\\
1.8111458073	-2.0093983028\\
1.7893018209	-1.9926133265\\
1.7678411841	-1.9756262475\\
1.7467612191	-1.9584383362\\
1.7260590322	-1.9410509061\\
1.7057314924	-1.9234652044\\
1.6857752645	-1.9056824608\\
1.6661868433	-1.8877038888\\
1.6469625361	-1.8695306389\\
1.6280984997	-1.8511638136\\
1.609590766	-1.8326044706\\
1.5914352373	-1.8138535919\\
1.5736277161	-1.7949120938\\
1.5561639264	-1.7757808286\\
1.5390395165	-1.7564605655\\
1.5222500824	-1.7369519971\\
1.5057911853	-1.7172557411\\
1.4896583586	-1.6973723292\\
1.4738471258	-1.6773022126\\
1.4583530139	-1.6570457631\\
1.4431715629	-1.6366032689\\
1.4282983382	-1.6159749386\\
1.4137289427	-1.5951609043\\
1.3994590235	-1.5741612199\\
1.3854842815	-1.5529758632\\
1.3718004856	-1.5316047503\\
1.3584034703	-1.5100477215\\
1.3452891468	-1.488304553\\
1.3324535167	-1.4663749736\\
1.3198926639	-1.4442586469\\
1.3076027684	-1.4219551908\\
1.2955801094	-1.399464182\\
1.2838210673	-1.3767851585\\
1.2723221278	-1.353917628\\
1.2610798849	-1.3308610749\\
1.2500910422	-1.3076149676\\
1.2393524145	-1.2841787652\\
1.22886093	-1.2605519292\\
1.2186136291	-1.2367339279\\
1.2086076653	-1.2127242489\\
1.1988403043	-1.1885224078\\
1.1893089229	-1.1641279594\\
1.1800110069	-1.1395405081\\
1.1709441484	-1.1147597201\\
1.1621060428	-1.0897853354\\
1.1534944836	-1.0646171804\\
1.1451073584	-1.0392551824\\
1.1369426409	-1.0136993826\\
1.1289983838	-0.98794995051\\
1.1212727096	-0.9620072015\\
1.1137637987	-0.93587160973\\
1.1064698758	-0.90954382591\\
1.0993891942	-0.88302469492\\
1.0925200155	-0.85631527127\\
1.0858605864	-0.82941683784\\
1.0794091099	-0.80233092322\\
1.0731637105	-0.77505931988\\
1.067122392	-0.7476041022\\
1.0612829839	-0.71996764303\\
1.0556430772	-0.69215263302\\
1.0501999416	-0.66416209651\\
1.0449504217	-0.63599940763\\
1.0398908039	-0.60766830659\\
1.0350166431	-0.57917291367\\
1.0303225356	-0.55051774264\\
1.0258018152	-0.52170771295\\
1.0214461414	-0.49274816079\\
1.0172449291	-0.46364484811\\
1.0131845461	-0.43440397316\\
1.0092471546	-0.40503217975\\
1.0054089979	-0.37553656972\\
1.0016377924	-0.34592472368\\
0.99788862615	-0.31620473681\\
0.99409725379	-0.28638529122\\
0.99016863314	-0.25647580595\\
0.98595623307	-0.22648677162\\
0.98122212496	-0.19643054715\\
0.97555334713	-0.16632345752\\
0.96816685946	-0.136192096\\
0.9573856425	-0.10609607997\\
0.93894272015	-0.076234226129\\
0.89941683675	-0.047635415889\\
0.80838945474	-0.026320112905\\
0.70311807578	-0.01804357165\\
0.60298032985	-0.014337731116\\
0.50309140379	-0.012218368249\\
0.40315157715	-0.010898140356\\
0.30317836462	-0.01005154642\\
0.20318820823	-0.0095213113758\\
0.10318989316	-0.0092268294122\\
0.0031883060395	-0.0091288003399\\
-1.8043484279e-06	-0.0091287076645\\
};
\addplot [color=black,solid,line width=2.0pt,forget plot]
  table[row sep=crcr]{%
2.0000689783	2.1314508024\\
1.975231865	2.1162785454\\
1.9507888522	2.1008954949\\
1.9267394877	2.085302569\\
1.9030830223	2.0695007707\\
1.8798184032	2.0534911832\\
1.8569442979	2.0372749234\\
1.8344590768	2.0208531823\\
1.8123608371	2.0042271756\\
1.7906474113	1.9873981487\\
1.7693163772	1.9703673541\\
1.7483650654	1.9531360313\\
1.7277906001	1.9357054586\\
1.7075898766	1.9180768492\\
1.6877595951	1.9002513972\\
1.6682962953	1.8822302798\\
1.6491963378	1.8640146109\\
1.6304559418	1.8456054556\\
1.6120712096	1.8270038344\\
1.5940381228	1.8082106927\\
1.5763525715	1.7892269102\\
1.5590103749	1.7700533029\\
1.5420072851	1.7506906049\\
1.5253390099	1.7311394742\\
1.50900123	1.7114004947\\
1.4929896064	1.6914741656\\
1.4772997984	1.6713609068\\
1.4619274769	1.6510610595\\
1.4468683342	1.6305748831\\
1.4321180961	1.6099025578\\
1.4176725363	1.5890441906\\
1.4035274807	1.5679998098\\
1.3896788188	1.5467693694\\
1.3761225201	1.5253527649\\
1.3628546291	1.5037498158\\
1.3498712796	1.4819602795\\
1.3371687089	1.4599838683\\
1.3247432504	1.4378202304\\
1.3125913487	1.4154689708\\
1.3007095638	1.3929296539\\
1.2890945755	1.3702018083\\
1.2777431886	1.3472849337\\
1.2666523386	1.3241785099\\
1.2558190943	1.3008820014\\
1.245240663	1.2773948663\\
1.234914395	1.2537165666\\
1.2248377862	1.2298465745\\
1.2150084831	1.2057843834\\
1.2054242856	1.1815295178\\
1.1960831518	1.1570815436\\
1.1869832011	1.1324400794\\
1.178122719	1.1076048089\\
1.1695001609	1.0825754922\\
1.1611141571	1.0573519789\\
1.1529635184	1.0319342238\\
1.1450472411	1.0063222977\\
1.1373645142	0.98051640439\\
1.1299147273	0.95451689627\\
1.122697479	0.92832428859\\
1.1157125875	0.90193927654\\
1.1089601041	0.87536275388\\
1.1024403276	0.84859582688\\
1.0961538232	0.82163983414\\
1.0901014458	0.79449636282\\
1.0842843681	0.7671672676\\
1.0787041168	0.73965468773\\
1.073362617	0.71196106325\\
1.0682622502	0.68408915383\\
1.063405928	0.6560420538\\
1.0587971871	0.62782320761\\
1.0544403152	0.5994364246\\
1.0503405162	0.57088589186\\
1.0465041327	0.54217618547\\
1.042938947	0.51331228096\\
1.0396545965	0.48429956127\\
1.0366631548	0.45514382189\\
1.033979961	0.42585127735\\
1.0316248275	0.39642856396\\
1.0296238438	0.36688274496\\
1.0280121447	0.33722131965\\
1.0268383025	0.30745224541\\
1.0261715659	0.27758398917\\
1.0261143537	0.24762565214\\
1.0268250447	0.21758727244\\
1.0285624962	0.18748060435\\
1.0317808465	0.15732130647\\
1.0373553028	0.12713593316\\
1.0472062064	0.096988851589\\
1.066396931	0.067117508885\\
1.1113487692	0.038881651965\\
1.2100764866	0.019658081689\\
1.31390049	0.012567562087\\
1.413859718	0.0091423711547\\
1.5137261053	0.0070791413556\\
1.6136488323	0.0057081758116\\
1.7136053474	0.0047382483938\\
1.8135797394	0.0040215182513\\
1.9135637537	0.0034753773202\\
2.0000032431	0.0031023133474\\
};
\addplot [<-,color=red,solid,line width=4.0pt,forget plot]
  table[row sep=crcr]{%
1.2086076653	-1.2127242489\\
1.1988403043	-1.1885224078\\
1.1893089229	-1.1641279594\\
1.1800110069	-1.1395405081\\
1.1709441484	-1.1147597201\\
1.1621060428	-1.0897853354\\
1.1534944836	-1.0646171804\\
1.1451073584	-1.0392551824\\
1.1369426409	-1.0136993826\\
1.1289983838	-0.98794995051\\
1.1212727096	-0.9620072015\\
1.1137637987	-0.93587160973\\
1.1064698758	-0.90954382591\\
1.0993891942	-0.88302469492\\
1.0925200155	-0.85631527127\\
1.0858605864	-0.82941683784\\
1.0794091099	-0.80233092322\\
1.0731637105	-0.77505931988\\
1.067122392	-0.7476041022\\
1.0612829839	-0.71996764303\\
1.0556430772	-0.69215263302\\
1.0501999416	-0.66416209651\\
1.0449504217	-0.63599940763\\
1.0398908039	-0.60766830659\\
1.0350166431	-0.57917291367\\
1.0303225356	-0.55051774264\\
1.0258018152	-0.52170771295\\
1.0214461414	-0.49274816079\\
1.0172449291	-0.46364484811\\
1.0131845461	-0.43440397316\\
1.0092471546	-0.40503217975\\
1.0054089979	-0.37553656972\\
1.0016377924	-0.34592472368\\
0.99788862615	-0.31620473681\\
0.99409725379	-0.28638529122\\
0.99016863314	-0.25647580595\\
0.98595623307	-0.22648677162\\
0.98122212496	-0.19643054715\\
0.97555334713	-0.16632345752\\
0.96816685946	-0.136192096\\
0.9573856425	-0.10609607997\\
0.93894272015	-0.076234226129\\
0.89941683675	-0.047635415889\\
0.80838945474	-0.026320112905\\
0.70311807578	-0.01804357165\\
0.60298032985	-0.014337731116\\
0.50309140379	-0.012218368249\\
0.40315157715	-0.010898140356\\
0.30317836462	-0.01005154642\\
0.20318820823	-0.0095213113758\\
0.10318989316	-0.0092268294122\\
0.0031883060395	-0.0091288003399\\
-1.8043484279e-06	-0.0091287076645\\
};
\addplot [dashed,color=black,line width=2.0pt,forget plot]
  table[row sep=crcr]{%
1.2086076653	 -2.5\\
1.2086076653	 2.5\\
};
\end{axis}

\end{tikzpicture}%

%% file: tikz/deflationbif2.tikz
%
%
\begin{tikzpicture}

\begin{axis}[%
width=4.520833in,
height=3.565625in,
at={(0.758333in,0.48125in)},
scale only axis,
xmin=0.2,
xmax=2,
ymin=-2.5,
ymax=2.5,
yticklabels={,,},
xticklabels={,,},
xtick = {0},
ytick = {0},
xlabel={\Huge $\lambda$},
ylabel={\Huge $u$}
]
\addplot [color=red, mark size = 5.0pt, mark=*,mark options={solid},forget plot]
  table[row sep=crcr]{%
1.2086076653	-1.2127242489\\
1.2054242856	1.1815295178\\
1.2100764866	0.019658081689\\
};

\addplot [color=black,solid,line width=2.0pt,forget plot]
  table[row sep=crcr]{%
2.0252982505	-2.1511656348\\
1.9999309186	-2.1362543389\\
1.9749563582	-2.1211296905\\
1.9503766066	-2.1057937463\\
1.926191269	-2.0902474329\\
1.9023996495	-2.0744917679\\
1.8790007424	-2.0585278475\\
1.8559932461	-2.0423568163\\
1.8333755614	-2.0259798846\\
1.8111458073	-2.0093983028\\
1.7893018209	-1.9926133265\\
1.7678411841	-1.9756262475\\
1.7467612191	-1.9584383362\\
1.7260590322	-1.9410509061\\
1.7057314924	-1.9234652044\\
1.6857752645	-1.9056824608\\
1.6661868433	-1.8877038888\\
1.6469625361	-1.8695306389\\
1.6280984997	-1.8511638136\\
1.609590766	-1.8326044706\\
1.5914352373	-1.8138535919\\
1.5736277161	-1.7949120938\\
1.5561639264	-1.7757808286\\
1.5390395165	-1.7564605655\\
1.5222500824	-1.7369519971\\
1.5057911853	-1.7172557411\\
1.4896583586	-1.6973723292\\
1.4738471258	-1.6773022126\\
1.4583530139	-1.6570457631\\
1.4431715629	-1.6366032689\\
1.4282983382	-1.6159749386\\
1.4137289427	-1.5951609043\\
1.3994590235	-1.5741612199\\
1.3854842815	-1.5529758632\\
1.3718004856	-1.5316047503\\
1.3584034703	-1.5100477215\\
1.3452891468	-1.488304553\\
1.3324535167	-1.4663749736\\
1.3198926639	-1.4442586469\\
1.3076027684	-1.4219551908\\
1.2955801094	-1.399464182\\
1.2838210673	-1.3767851585\\
1.2723221278	-1.353917628\\
1.2610798849	-1.3308610749\\
1.2500910422	-1.3076149676\\
1.2393524145	-1.2841787652\\
1.22886093	-1.2605519292\\
1.2186136291	-1.2367339279\\
1.2086076653	-1.2127242489\\
1.1988403043	-1.1885224078\\
1.1893089229	-1.1641279594\\
1.1800110069	-1.1395405081\\
1.1709441484	-1.1147597201\\
1.1621060428	-1.0897853354\\
1.1534944836	-1.0646171804\\
1.1451073584	-1.0392551824\\
1.1369426409	-1.0136993826\\
1.1289983838	-0.98794995051\\
1.1212727096	-0.9620072015\\
1.1137637987	-0.93587160973\\
1.1064698758	-0.90954382591\\
1.0993891942	-0.88302469492\\
1.0925200155	-0.85631527127\\
1.0858605864	-0.82941683784\\
1.0794091099	-0.80233092322\\
1.0731637105	-0.77505931988\\
1.067122392	-0.7476041022\\
1.0612829839	-0.71996764303\\
1.0556430772	-0.69215263302\\
1.0501999416	-0.66416209651\\
1.0449504217	-0.63599940763\\
1.0398908039	-0.60766830659\\
1.0350166431	-0.57917291367\\
1.0303225356	-0.55051774264\\
1.0258018152	-0.52170771295\\
1.0214461414	-0.49274816079\\
1.0172449291	-0.46364484811\\
1.0131845461	-0.43440397316\\
1.0092471546	-0.40503217975\\
1.0054089979	-0.37553656972\\
1.0016377924	-0.34592472368\\
0.99788862615	-0.31620473681\\
0.99409725379	-0.28638529122\\
0.99016863314	-0.25647580595\\
0.98595623307	-0.22648677162\\
0.98122212496	-0.19643054715\\
0.97555334713	-0.16632345752\\
0.96816685946	-0.136192096\\
0.9573856425	-0.10609607997\\
0.93894272015	-0.076234226129\\
0.89941683675	-0.047635415889\\
0.80838945474	-0.026320112905\\
0.70311807578	-0.01804357165\\
0.60298032985	-0.014337731116\\
0.50309140379	-0.012218368249\\
0.40315157715	-0.010898140356\\
0.30317836462	-0.01005154642\\
0.20318820823	-0.0095213113758\\
0.10318989316	-0.0092268294122\\
0.0031883060395	-0.0091288003399\\
-1.8043484279e-06	-0.0091287076645\\
};
\addplot [color=black,solid,line width=2.0pt,forget plot]
  table[row sep=crcr]{%
2.0000689783	2.1314508024\\
1.975231865	2.1162785454\\
1.9507888522	2.1008954949\\
1.9267394877	2.085302569\\
1.9030830223	2.0695007707\\
1.8798184032	2.0534911832\\
1.8569442979	2.0372749234\\
1.8344590768	2.0208531823\\
1.8123608371	2.0042271756\\
1.7906474113	1.9873981487\\
1.7693163772	1.9703673541\\
1.7483650654	1.9531360313\\
1.7277906001	1.9357054586\\
1.7075898766	1.9180768492\\
1.6877595951	1.9002513972\\
1.6682962953	1.8822302798\\
1.6491963378	1.8640146109\\
1.6304559418	1.8456054556\\
1.6120712096	1.8270038344\\
1.5940381228	1.8082106927\\
1.5763525715	1.7892269102\\
1.5590103749	1.7700533029\\
1.5420072851	1.7506906049\\
1.5253390099	1.7311394742\\
1.50900123	1.7114004947\\
1.4929896064	1.6914741656\\
1.4772997984	1.6713609068\\
1.4619274769	1.6510610595\\
1.4468683342	1.6305748831\\
1.4321180961	1.6099025578\\
1.4176725363	1.5890441906\\
1.4035274807	1.5679998098\\
1.3896788188	1.5467693694\\
1.3761225201	1.5253527649\\
1.3628546291	1.5037498158\\
1.3498712796	1.4819602795\\
1.3371687089	1.4599838683\\
1.3247432504	1.4378202304\\
1.3125913487	1.4154689708\\
1.3007095638	1.3929296539\\
1.2890945755	1.3702018083\\
1.2777431886	1.3472849337\\
1.2666523386	1.3241785099\\
1.2558190943	1.3008820014\\
1.245240663	1.2773948663\\
1.234914395	1.2537165666\\
1.2248377862	1.2298465745\\
1.2150084831	1.2057843834\\
1.2054242856	1.1815295178\\
1.1960831518	1.1570815436\\
1.1869832011	1.1324400794\\
1.178122719	1.1076048089\\
1.1695001609	1.0825754922\\
1.1611141571	1.0573519789\\
1.1529635184	1.0319342238\\
1.1450472411	1.0063222977\\
1.1373645142	0.98051640439\\
1.1299147273	0.95451689627\\
1.122697479	0.92832428859\\
1.1157125875	0.90193927654\\
1.1089601041	0.87536275388\\
1.1024403276	0.84859582688\\
1.0961538232	0.82163983414\\
1.0901014458	0.79449636282\\
1.0842843681	0.7671672676\\
1.0787041168	0.73965468773\\
1.073362617	0.71196106325\\
1.0682622502	0.68408915383\\
1.063405928	0.6560420538\\
1.0587971871	0.62782320761\\
1.0544403152	0.5994364246\\
1.0503405162	0.57088589186\\
1.0465041327	0.54217618547\\
1.042938947	0.51331228096\\
1.0396545965	0.48429956127\\
1.0366631548	0.45514382189\\
1.033979961	0.42585127735\\
1.0316248275	0.39642856396\\
1.0296238438	0.36688274496\\
1.0280121447	0.33722131965\\
1.0268383025	0.30745224541\\
1.0261715659	0.27758398917\\
1.0261143537	0.24762565214\\
1.0268250447	0.21758727244\\
1.0285624962	0.18748060435\\
1.0317808465	0.15732130647\\
1.0373553028	0.12713593316\\
1.0472062064	0.096988851589\\
1.066396931	0.067117508885\\
1.1113487692	0.038881651965\\
1.2100764866	0.019658081689\\
1.31390049	0.012567562087\\
1.413859718	0.0091423711547\\
1.5137261053	0.0070791413556\\
1.6136488323	0.0057081758116\\
1.7136053474	0.0047382483938\\
1.8135797394	0.0040215182513\\
1.9135637537	0.0034753773202\\
2.0000032431	0.0031023133474\\
};
\addplot [color=red,solid,line width=4.0pt,forget plot]
  table[row sep=crcr]{%
1.2086076653	-1.2127242489\\
1.1988403043	-1.1885224078\\
1.1893089229	-1.1641279594\\
1.1800110069	-1.1395405081\\
1.1709441484	-1.1147597201\\
1.1621060428	-1.0897853354\\
1.1534944836	-1.0646171804\\
1.1451073584	-1.0392551824\\
1.1369426409	-1.0136993826\\
1.1289983838	-0.98794995051\\
1.1212727096	-0.9620072015\\
1.1137637987	-0.93587160973\\
1.1064698758	-0.90954382591\\
1.0993891942	-0.88302469492\\
1.0925200155	-0.85631527127\\
1.0858605864	-0.82941683784\\
1.0794091099	-0.80233092322\\
1.0731637105	-0.77505931988\\
1.067122392	-0.7476041022\\
1.0612829839	-0.71996764303\\
1.0556430772	-0.69215263302\\
1.0501999416	-0.66416209651\\
1.0449504217	-0.63599940763\\
1.0398908039	-0.60766830659\\
1.0350166431	-0.57917291367\\
1.0303225356	-0.55051774264\\
1.0258018152	-0.52170771295\\
1.0214461414	-0.49274816079\\
1.0172449291	-0.46364484811\\
1.0131845461	-0.43440397316\\
1.0092471546	-0.40503217975\\
1.0054089979	-0.37553656972\\
1.0016377924	-0.34592472368\\
0.99788862615	-0.31620473681\\
0.99409725379	-0.28638529122\\
0.99016863314	-0.25647580595\\
0.98595623307	-0.22648677162\\
0.98122212496	-0.19643054715\\
0.97555334713	-0.16632345752\\
0.96816685946	-0.136192096\\
0.9573856425	-0.10609607997\\
0.93894272015	-0.076234226129\\
0.89941683675	-0.047635415889\\
0.80838945474	-0.026320112905\\
0.70311807578	-0.01804357165\\
0.60298032985	-0.014337731116\\
0.50309140379	-0.012218368249\\
0.40315157715	-0.010898140356\\
0.30317836462	-0.01005154642\\
0.20318820823	-0.0095213113758\\
0.10318989316	-0.0092268294122\\
0.0031883060395	-0.0091288003399\\
-1.8043484279e-06	-0.0091287076645\\
};
\addplot [dashed,color=black,line width=2.0pt,forget plot]
  table[row sep=crcr]{%
1.2086076653	 -2.5\\
1.2086076653	 2.5\\
};
\end{axis}

\end{tikzpicture}%

%% file: tikz/deflationbif3.tikz
%
%
\begin{tikzpicture}

\begin{axis}[%
width=4.520833in,
height=3.565625in,
at={(0.758333in,0.48125in)},
scale only axis,
xmin=0.2,
xmax=2,
ymin=-2.5,
ymax=2.5,
yticklabels={,,},
xticklabels={,,},
xtick = {0},
ytick = {0},
xlabel={\Huge $\lambda$},
ylabel={\Huge $u$}
]
\addplot [color=red, mark size = 5.0pt, mark=*,mark options={solid},forget plot]
  table[row sep=crcr]{%
1.2086076653	-1.2127242489\\
1.2054242856	1.1815295178\\
1.2100764866	0.019658081689\\
};

\addplot [color=black,solid,line width=2.0pt,forget plot]
  table[row sep=crcr]{%
2.0252982505	-2.1511656348\\
1.9999309186	-2.1362543389\\
1.9749563582	-2.1211296905\\
1.9503766066	-2.1057937463\\
1.926191269	-2.0902474329\\
1.9023996495	-2.0744917679\\
1.8790007424	-2.0585278475\\
1.8559932461	-2.0423568163\\
1.8333755614	-2.0259798846\\
1.8111458073	-2.0093983028\\
1.7893018209	-1.9926133265\\
1.7678411841	-1.9756262475\\
1.7467612191	-1.9584383362\\
1.7260590322	-1.9410509061\\
1.7057314924	-1.9234652044\\
1.6857752645	-1.9056824608\\
1.6661868433	-1.8877038888\\
1.6469625361	-1.8695306389\\
1.6280984997	-1.8511638136\\
1.609590766	-1.8326044706\\
1.5914352373	-1.8138535919\\
1.5736277161	-1.7949120938\\
1.5561639264	-1.7757808286\\
1.5390395165	-1.7564605655\\
1.5222500824	-1.7369519971\\
1.5057911853	-1.7172557411\\
1.4896583586	-1.6973723292\\
1.4738471258	-1.6773022126\\
1.4583530139	-1.6570457631\\
1.4431715629	-1.6366032689\\
1.4282983382	-1.6159749386\\
1.4137289427	-1.5951609043\\
1.3994590235	-1.5741612199\\
1.3854842815	-1.5529758632\\
1.3718004856	-1.5316047503\\
1.3584034703	-1.5100477215\\
1.3452891468	-1.488304553\\
1.3324535167	-1.4663749736\\
1.3198926639	-1.4442586469\\
1.3076027684	-1.4219551908\\
1.2955801094	-1.399464182\\
1.2838210673	-1.3767851585\\
1.2723221278	-1.353917628\\
1.2610798849	-1.3308610749\\
1.2500910422	-1.3076149676\\
1.2393524145	-1.2841787652\\
1.22886093	-1.2605519292\\
1.2186136291	-1.2367339279\\
1.2086076653	-1.2127242489\\
1.1988403043	-1.1885224078\\
1.1893089229	-1.1641279594\\
1.1800110069	-1.1395405081\\
1.1709441484	-1.1147597201\\
1.1621060428	-1.0897853354\\
1.1534944836	-1.0646171804\\
1.1451073584	-1.0392551824\\
1.1369426409	-1.0136993826\\
1.1289983838	-0.98794995051\\
1.1212727096	-0.9620072015\\
1.1137637987	-0.93587160973\\
1.1064698758	-0.90954382591\\
1.0993891942	-0.88302469492\\
1.0925200155	-0.85631527127\\
1.0858605864	-0.82941683784\\
1.0794091099	-0.80233092322\\
1.0731637105	-0.77505931988\\
1.067122392	-0.7476041022\\
1.0612829839	-0.71996764303\\
1.0556430772	-0.69215263302\\
1.0501999416	-0.66416209651\\
1.0449504217	-0.63599940763\\
1.0398908039	-0.60766830659\\
1.0350166431	-0.57917291367\\
1.0303225356	-0.55051774264\\
1.0258018152	-0.52170771295\\
1.0214461414	-0.49274816079\\
1.0172449291	-0.46364484811\\
1.0131845461	-0.43440397316\\
1.0092471546	-0.40503217975\\
1.0054089979	-0.37553656972\\
1.0016377924	-0.34592472368\\
0.99788862615	-0.31620473681\\
0.99409725379	-0.28638529122\\
0.99016863314	-0.25647580595\\
0.98595623307	-0.22648677162\\
0.98122212496	-0.19643054715\\
0.97555334713	-0.16632345752\\
0.96816685946	-0.136192096\\
0.9573856425	-0.10609607997\\
0.93894272015	-0.076234226129\\
0.89941683675	-0.047635415889\\
0.80838945474	-0.026320112905\\
0.70311807578	-0.01804357165\\
0.60298032985	-0.014337731116\\
0.50309140379	-0.012218368249\\
0.40315157715	-0.010898140356\\
0.30317836462	-0.01005154642\\
0.20318820823	-0.0095213113758\\
0.10318989316	-0.0092268294122\\
0.0031883060395	-0.0091288003399\\
-1.8043484279e-06	-0.0091287076645\\
};
\addplot [color=black,solid,line width=2.0pt,forget plot]
  table[row sep=crcr]{%
2.0000689783	2.1314508024\\
1.975231865	2.1162785454\\
1.9507888522	2.1008954949\\
1.9267394877	2.085302569\\
1.9030830223	2.0695007707\\
1.8798184032	2.0534911832\\
1.8569442979	2.0372749234\\
1.8344590768	2.0208531823\\
1.8123608371	2.0042271756\\
1.7906474113	1.9873981487\\
1.7693163772	1.9703673541\\
1.7483650654	1.9531360313\\
1.7277906001	1.9357054586\\
1.7075898766	1.9180768492\\
1.6877595951	1.9002513972\\
1.6682962953	1.8822302798\\
1.6491963378	1.8640146109\\
1.6304559418	1.8456054556\\
1.6120712096	1.8270038344\\
1.5940381228	1.8082106927\\
1.5763525715	1.7892269102\\
1.5590103749	1.7700533029\\
1.5420072851	1.7506906049\\
1.5253390099	1.7311394742\\
1.50900123	1.7114004947\\
1.4929896064	1.6914741656\\
1.4772997984	1.6713609068\\
1.4619274769	1.6510610595\\
1.4468683342	1.6305748831\\
1.4321180961	1.6099025578\\
1.4176725363	1.5890441906\\
1.4035274807	1.5679998098\\
1.3896788188	1.5467693694\\
1.3761225201	1.5253527649\\
1.3628546291	1.5037498158\\
1.3498712796	1.4819602795\\
1.3371687089	1.4599838683\\
1.3247432504	1.4378202304\\
1.3125913487	1.4154689708\\
1.3007095638	1.3929296539\\
1.2890945755	1.3702018083\\
1.2777431886	1.3472849337\\
1.2666523386	1.3241785099\\
1.2558190943	1.3008820014\\
1.245240663	1.2773948663\\
1.234914395	1.2537165666\\
1.2248377862	1.2298465745\\
1.2150084831	1.2057843834\\
1.2054242856	1.1815295178\\
1.1960831518	1.1570815436\\
1.1869832011	1.1324400794\\
1.178122719	1.1076048089\\
1.1695001609	1.0825754922\\
1.1611141571	1.0573519789\\
1.1529635184	1.0319342238\\
1.1450472411	1.0063222977\\
1.1373645142	0.98051640439\\
1.1299147273	0.95451689627\\
1.122697479	0.92832428859\\
1.1157125875	0.90193927654\\
1.1089601041	0.87536275388\\
1.1024403276	0.84859582688\\
1.0961538232	0.82163983414\\
1.0901014458	0.79449636282\\
1.0842843681	0.7671672676\\
1.0787041168	0.73965468773\\
1.073362617	0.71196106325\\
1.0682622502	0.68408915383\\
1.063405928	0.6560420538\\
1.0587971871	0.62782320761\\
1.0544403152	0.5994364246\\
1.0503405162	0.57088589186\\
1.0465041327	0.54217618547\\
1.042938947	0.51331228096\\
1.0396545965	0.48429956127\\
1.0366631548	0.45514382189\\
1.033979961	0.42585127735\\
1.0316248275	0.39642856396\\
1.0296238438	0.36688274496\\
1.0280121447	0.33722131965\\
1.0268383025	0.30745224541\\
1.0261715659	0.27758398917\\
1.0261143537	0.24762565214\\
1.0268250447	0.21758727244\\
1.0285624962	0.18748060435\\
1.0317808465	0.15732130647\\
1.0373553028	0.12713593316\\
1.0472062064	0.096988851589\\
1.066396931	0.067117508885\\
1.1113487692	0.038881651965\\
1.2100764866	0.019658081689\\
1.31390049	0.012567562087\\
1.413859718	0.0091423711547\\
1.5137261053	0.0070791413556\\
1.6136488323	0.0057081758116\\
1.7136053474	0.0047382483938\\
1.8135797394	0.0040215182513\\
1.9135637537	0.0034753773202\\
2.0000032431	0.0031023133474\\
};
\addplot [<-,color=red,solid,line width=4.0pt,forget plot]
  table[row sep=crcr]{%
1.0472062064	0.096988851589\\
1.066396931	0.067117508885\\
1.1113487692	0.038881651965\\
1.2100764866	0.019658081689\\
};
\addplot [->,color=red,solid,line width=4.0pt,forget plot]
  table[row sep=crcr]{%
1.2100764866	0.019658081689\\
1.31390049	0.012567562087\\
1.413859718	0.0091423711547\\
1.5137261053	0.0070791413556\\
};
\addplot [->,color=red,solid,line width=4.0pt,forget plot]
  table[row sep=crcr]{%
1.2054242856	1.1815295178\\
1.1960831518	1.1570815436\\
1.1869832011	1.1324400794\\
1.178122719	1.1076048089\\
1.1695001609	1.0825754922\\
1.1611141571	1.0573519789\\
1.1529635184	1.0319342238\\
1.1450472411	1.0063222977\\
1.1373645142	0.98051640439\\
1.1299147273	0.95451689627\\
1.122697479	0.92832428859\\
1.1157125875	0.90193927654\\
1.1089601041	0.87536275388\\
1.1024403276	0.84859582688\\
1.0961538232	0.82163983414\\
1.0901014458	0.79449636282\\
1.0842843681	0.7671672676\\
};
\addplot [<-,color=red,solid,line width=4.0pt,forget plot]
  table[row sep=crcr]{%
1.50900123	1.7114004947\\
1.4929896064	1.6914741656\\
1.4772997984	1.6713609068\\
1.4619274769	1.6510610595\\
1.4468683342	1.6305748831\\
1.4321180961	1.6099025578\\
1.4176725363	1.5890441906\\
1.4035274807	1.5679998098\\
1.3896788188	1.5467693694\\
1.3761225201	1.5253527649\\
1.3628546291	1.5037498158\\
1.3498712796	1.4819602795\\
1.3371687089	1.4599838683\\
1.3247432504	1.4378202304\\
1.3125913487	1.4154689708\\
1.3007095638	1.3929296539\\
1.2890945755	1.3702018083\\
1.2777431886	1.3472849337\\
1.2666523386	1.3241785099\\
1.2558190943	1.3008820014\\
1.245240663	1.2773948663\\
1.234914395	1.2537165666\\
1.2248377862	1.2298465745\\
1.2150084831	1.2057843834\\
};
\addplot [<-,color=red,solid,line width=4.0pt,forget plot]
  table[row sep=crcr]{%
  1.5057911853	-1.7172557411\\
1.4896583586	-1.6973723292\\
1.4738471258	-1.6773022126\\
1.4583530139	-1.6570457631\\
1.4431715629	-1.6366032689\\
1.4282983382	-1.6159749386\\
1.4137289427	-1.5951609043\\
1.3994590235	-1.5741612199\\
1.3854842815	-1.5529758632\\
1.3718004856	-1.5316047503\\
1.3584034703	-1.5100477215\\
1.3452891468	-1.488304553\\
1.3324535167	-1.4663749736\\
1.3198926639	-1.4442586469\\
1.3076027684	-1.4219551908\\
1.2955801094	-1.399464182\\
1.2838210673	-1.3767851585\\
1.2723221278	-1.353917628\\
1.2610798849	-1.3308610749\\
1.2500910422	-1.3076149676\\
1.2393524145	-1.2841787652\\
1.22886093	-1.2605519292\\
1.2186136291	-1.2367339279\\
1.2086076653	-1.2127242489\\
1.1988403043	-1.1885224078\\
1.1893089229	-1.1641279594\\
1.1800110069	-1.1395405081\\
1.1709441484	-1.1147597201\\
1.1621060428	-1.0897853354\\
1.1534944836	-1.0646171804\\
1.1451073584	-1.0392551824\\
1.1369426409	-1.0136993826\\
1.1289983838	-0.98794995051\\
1.1212727096	-0.9620072015\\
1.1137637987	-0.93587160973\\
1.1064698758	-0.90954382591\\
1.0993891942	-0.88302469492\\
1.0925200155	-0.85631527127\\
1.0858605864	-0.82941683784\\
1.0794091099	-0.80233092322\\
1.0731637105	-0.77505931988\\
1.067122392	-0.7476041022\\
1.0612829839	-0.71996764303\\
1.0556430772	-0.69215263302\\
1.0501999416	-0.66416209651\\
1.0449504217	-0.63599940763\\
1.0398908039	-0.60766830659\\
1.0350166431	-0.57917291367\\
1.0303225356	-0.55051774264\\
1.0258018152	-0.52170771295\\
1.0214461414	-0.49274816079\\
1.0172449291	-0.46364484811\\
1.0131845461	-0.43440397316\\
1.0092471546	-0.40503217975\\
1.0054089979	-0.37553656972\\
1.0016377924	-0.34592472368\\
0.99788862615	-0.31620473681\\
0.99409725379	-0.28638529122\\
0.99016863314	-0.25647580595\\
0.98595623307	-0.22648677162\\
0.98122212496	-0.19643054715\\
0.97555334713	-0.16632345752\\
0.96816685946	-0.136192096\\
0.9573856425	-0.10609607997\\
0.93894272015	-0.076234226129\\
0.89941683675	-0.047635415889\\
0.80838945474	-0.026320112905\\
0.70311807578	-0.01804357165\\
0.60298032985	-0.014337731116\\
0.50309140379	-0.012218368249\\
0.40315157715	-0.010898140356\\
0.30317836462	-0.01005154642\\
0.20318820823	-0.0095213113758\\
0.10318989316	-0.0092268294122\\
0.0031883060395	-0.0091288003399\\
-1.8043484279e-06	-0.0091287076645\\
};
\addplot [dashed,color=black,line width=2.0pt,forget plot]
  table[row sep=crcr]{%
1.2086076653	 -2.5\\
1.2086076653	 2.5\\
};
\end{axis}

\end{tikzpicture}%

%% file: tikz/sufficientconditions1.tikz
%
%
\begin{tikzpicture}

\begin{axis}[%
width=4.822222in,
height=3.616667in,
at={(0.808889in,0.606667in)},
scale only axis,
separate axis lines,
every outer x axis line/.append style={black},
every x tick label/.append style={font=\color{black}},
xmin=-4,
xmax=6,
xtick={\empty},
every outer y axis line/.append style={black},
every y tick label/.append style={font=\color{black}},
ymin=-4.5,
ymax=3,
ytick={\empty}
]
\addplot [color=red,line width=3.0pt,only marks,mark=asterisk,mark options={solid},forget plot]
  table[row sep=crcr]{%
-2	1\\
2	-1\\
};
\addplot [color=black,line width=3.0pt,only marks,mark=asterisk,mark options={solid},forget plot]
  table[row sep=crcr]{%
-0.3	0.7\\
};
\node at (axis cs:-2,1) [anchor=south east] {\LARGE $u_1^*$};
\node at (axis cs:2,-1) [anchor=south west] {\LARGE $u_2^*$};
\node at (axis cs:-0.2,0.5) [anchor=south west] {\LARGE $u_0$};

\node at (axis cs:-2.3,0.2) [anchor=south east] {\large $\rho_1$};
\node at (axis cs:1.5,-1.5) [anchor=south east] {\large $\rho_2$};
\node at (axis cs:5.2,2) [anchor = south] {\huge $D$};

\addplot [->,color=black,solid,line width=1.0pt,forget plot]
  table[row sep=crcr]{%
-2	1\\
-2.785	-0.698166421207658\\
};
\addplot [->,color=black,solid,line width=1.0pt,forget plot]
  table[row sep=crcr]{%
2	-1\\
0.873991130820399	-2.11\\
};

\addplot [color=blue,dashed,line width=2.0pt,forget plot]
  table[row sep=crcr]{%
-3.8	1.50909090909091\\
-3.8027027027027	1.5\\
-3.82702702702703	1.4\\
-3.84594594594595	1.3\\
-3.85945945945946	1.2\\
-3.86756756756757	1.1\\
-3.87027027027027	1\\
-3.86756756756757	0.9\\
-3.85945945945946	0.8\\
-3.84594594594595	0.7\\
-3.82702702702703	0.6\\
-3.8027027027027	0.5\\
-3.8	0.490909090909091\\
-3.77142857142857	0.4\\
-3.73428571428571	0.3\\
-3.7	0.220000000000001\\
-3.69090909090909	0.2\\
-3.63939393939394	0.0999999999999996\\
-3.6	0.0315789473684211\\
-3.58064516129032	0\\
-3.51290322580645	-0.1\\
-3.5	-0.117391304347826\\
-3.43448275862069	-0.2\\
-3.4	-0.24\\
-3.34444444444444	-0.3\\
-3.3	-0.344444444444445\\
-3.24	-0.4\\
-3.2	-0.434482758620689\\
-3.11739130434783	-0.5\\
-3.1	-0.512903225806452\\
-3	-0.580645161290323\\
-2.96842105263158	-0.6\\
-2.9	-0.639393939393939\\
-2.8	-0.690909090909091\\
-2.78	-0.7\\
-2.7	-0.734285714285714\\
-2.6	-0.771428571428571\\
-2.50909090909091	-0.8\\
-2.5	-0.802702702702703\\
-2.4	-0.827027027027027\\
-2.3	-0.845945945945946\\
-2.2	-0.85945945945946\\
-2.1	-0.867567567567568\\
-2	-0.87027027027027\\
-1.9	-0.867567567567568\\
-1.8	-0.85945945945946\\
-1.7	-0.845945945945946\\
-1.6	-0.827027027027027\\
-1.5	-0.802702702702703\\
-1.49090909090909	-0.8\\
-1.4	-0.771428571428571\\
-1.3	-0.734285714285714\\
-1.22	-0.7\\
-1.2	-0.690909090909091\\
-1.1	-0.639393939393939\\
-1.03157894736842	-0.6\\
-1	-0.580645161290323\\
-0.9	-0.512903225806452\\
-0.882608695652174	-0.5\\
-0.8	-0.434482758620689\\
-0.760000000000001	-0.4\\
-0.7	-0.344444444444444\\
-0.655555555555556	-0.3\\
-0.6	-0.239999999999999\\
-0.565517241379311	-0.2\\
-0.5	-0.117391304347826\\
-0.487096774193548	-0.1\\
-0.419354838709677	0\\
-0.4	0.0315789473684211\\
-0.360606060606061	0.0999999999999996\\
-0.309090909090909	0.2\\
-0.3	0.220000000000001\\
-0.265714285714286	0.3\\
-0.228571428571429	0.4\\
-0.2	0.490909090909092\\
-0.197297297297297	0.5\\
-0.172972972972973	0.6\\
-0.154054054054054	0.7\\
-0.140540540540541	0.8\\
-0.132432432432433	0.9\\
-0.12972972972973	1\\
-0.132432432432433	1.1\\
-0.140540540540541	1.2\\
-0.154054054054054	1.3\\
-0.172972972972973	1.4\\
-0.197297297297297	1.5\\
-0.2	1.50909090909091\\
-0.228571428571428	1.6\\
-0.265714285714286	1.7\\
-0.3	1.78\\
-0.309090909090909	1.8\\
-0.360606060606061	1.9\\
-0.4	1.96842105263158\\
-0.419354838709677	2\\
-0.487096774193548	2.1\\
-0.5	2.11739130434783\\
-0.565517241379311	2.2\\
-0.6	2.24\\
-0.655555555555555	2.3\\
-0.7	2.34444444444444\\
-0.76	2.4\\
-0.8	2.43448275862069\\
-0.882608695652174	2.5\\
-0.9	2.51290322580645\\
-1	2.58064516129032\\
-1.03157894736842	2.6\\
-1.1	2.63939393939394\\
-1.2	2.69090909090909\\
-1.22	2.7\\
-1.3	2.73428571428571\\
-1.4	2.77142857142857\\
-1.49090909090909	2.8\\
-1.5	2.8027027027027\\
-1.6	2.82702702702703\\
-1.7	2.84594594594595\\
-1.8	2.85945945945946\\
-1.9	2.86756756756757\\
-2	2.87027027027027\\
-2.1	2.86756756756757\\
-2.2	2.85945945945946\\
-2.3	2.84594594594595\\
-2.4	2.82702702702703\\
-2.5	2.8027027027027\\
-2.50909090909091	2.8\\
-2.6	2.77142857142857\\
-2.7	2.73428571428571\\
-2.78	2.7\\
-2.8	2.69090909090909\\
-2.9	2.63939393939394\\
-2.96842105263158	2.6\\
-3	2.58064516129032\\
-3.1	2.51290322580645\\
-3.11739130434783	2.5\\
-3.2	2.43448275862069\\
-3.24	2.4\\
-3.3	2.34444444444444\\
-3.34444444444444	2.3\\
-3.4	2.24\\
-3.43448275862069	2.2\\
-3.5	2.11739130434783\\
-3.51290322580645	2.1\\
-3.58064516129032	2\\
-3.6	1.96842105263158\\
-3.63939393939394	1.9\\
-3.69090909090909	1.8\\
-3.7	1.78\\
-3.73428571428571	1.7\\
-3.77142857142857	1.6\\
-3.8	1.50909090909091\\
};
\addplot [color=blue,dashed,line width=2.0pt,forget plot]
  table[row sep=crcr]{%
3.5	-0.5\\
3.46206896551724	-0.4\\
3.41724137931034	-0.3\\
3.4	-0.266666666666667\\
3.36296296296296	-0.2\\
3.3	-0.1\\
3.224	0\\
3.2	0.0285714285714291\\
3.13478260869565	0.0999999999999996\\
3.1	0.134782608695653\\
3.02857142857143	0.2\\
3	0.224\\
2.9	0.3\\
2.8	0.362962962962963\\
2.73333333333333	0.4\\
2.7	0.417241379310345\\
2.6	0.462068965517241\\
2.5	0.5\\
2.4	0.529032258064516\\
2.3	0.551612903225806\\
2.2	0.567741935483871\\
2.1	0.57741935483871\\
2	0.580645161290323\\
1.9	0.577419354838709\\
1.8	0.567741935483871\\
1.7	0.551612903225806\\
1.6	0.529032258064516\\
1.5	0.5\\
1.4	0.462068965517241\\
1.3	0.417241379310344\\
1.26666666666667	0.4\\
1.2	0.362962962962962\\
1.1	0.3\\
1	0.224\\
0.971428571428571	0.2\\
0.9	0.134782608695653\\
0.865217391304347	0.0999999999999996\\
0.800000000000001	0.0285714285714291\\
0.776	0\\
0.7	-0.0999999999999994\\
0.637037037037037	-0.2\\
0.600000000000001	-0.266666666666665\\
0.582758620689655	-0.3\\
0.537931034482758	-0.4\\
0.5	-0.5\\
0.470967741935484	-0.6\\
0.448387096774194	-0.7\\
0.432258064516129	-0.8\\
0.422580645161291	-0.9\\
0.419354838709677	-1\\
0.42258064516129	-1.1\\
0.432258064516129	-1.2\\
0.448387096774194	-1.3\\
0.470967741935484	-1.4\\
0.5	-1.5\\
0.537931034482758	-1.6\\
0.582758620689655	-1.7\\
0.600000000000001	-1.73333333333333\\
0.637037037037037	-1.8\\
0.7	-1.9\\
0.776	-2\\
0.800000000000001	-2.02857142857143\\
0.865217391304347	-2.1\\
0.9	-2.13478260869565\\
0.971428571428571	-2.2\\
1	-2.224\\
1.1	-2.3\\
1.2	-2.36296296296296\\
1.26666666666667	-2.4\\
1.3	-2.41724137931034\\
1.4	-2.46206896551724\\
1.5	-2.5\\
1.6	-2.52903225806452\\
1.7	-2.55161290322581\\
1.8	-2.56774193548387\\
1.9	-2.57741935483871\\
2	-2.58064516129032\\
2.1	-2.57741935483871\\
2.2	-2.56774193548387\\
2.3	-2.55161290322581\\
2.4	-2.52903225806452\\
2.5	-2.5\\
2.6	-2.46206896551724\\
2.7	-2.41724137931035\\
2.73333333333333	-2.4\\
2.8	-2.36296296296296\\
2.9	-2.3\\
3	-2.224\\
3.02857142857143	-2.2\\
3.1	-2.13478260869565\\
3.13478260869565	-2.1\\
3.2	-2.02857142857143\\
3.224	-2\\
3.3	-1.9\\
3.36296296296296	-1.8\\
3.4	-1.73333333333333\\
3.41724137931034	-1.7\\
3.46206896551724	-1.6\\
3.5	-1.5\\
3.52903225806452	-1.4\\
3.55161290322581	-1.3\\
3.56774193548387	-1.2\\
3.57741935483871	-1.1\\
3.58064516129032	-1\\
3.57741935483871	-0.9\\
3.56774193548387	-0.8\\
3.55161290322581	-0.7\\
3.52903225806452	-0.6\\
3.5	-0.5\\
};
\addplot [color=black,dashed,line width=1.0pt,forget plot]
  table[row sep=crcr]{%
-0.3	0.7\\
-0.343589743589744	0.690667251410148\\
-0.387179487179487	0.773815827780037\\
-0.430769230769231	0.673155278567666\\
-0.474358974358974	0.700805937132635\\
-0.517948717948718	0.775823132299339\\
-0.561538461538462	0.696284618619216\\
-0.605128205128205	0.757596629495177\\
-0.648717948717949	0.797613219036216\\
-0.692307692307692	0.797889558261297\\
-0.735897435897436	0.733922634559372\\
-0.779487179487179	0.777162695563987\\
-0.823076923076923	0.774353336997309\\
-0.866666666666667	0.797419699663723\\
-0.91025641025641	0.818180715707199\\
-0.953846153846154	0.772881298810958\\
-0.997435897435897	0.871716595946247\\
-1.04102564102564	0.847592540569097\\
-1.08461538461538	0.798014287034419\\
-1.12820512820513	0.854063117074045\\
-1.17179487179487	0.905640090283417\\
-1.21538461538462	0.895510687580617\\
-1.25897435897436	0.891513225668248\\
-1.3025641025641	0.872869391121006\\
-1.34615384615385	0.855853556095368\\
-1.38974358974359	0.900619872836935\\
-1.43333333333333	0.869852479546569\\
-1.47692307692308	0.886008498674815\\
-1.52051282051282	0.964680883250648\\
-1.56410256410256	0.969303327996252\\
-1.60769230769231	0.966119961669133\\
-1.65128205128205	0.989558715321564\\
-1.69487179487179	0.907614640246951\\
-1.73846153846154	0.95595109163453\\
-1.78205128205128	0.976779099507579\\
-1.82564102564103	1.01881265073888\\
-1.86923076923077	0.996599267383734\\
-1.91282051282051	0.971318523374498\\
-1.95641025641026	1.02110860207763\\
-2	1\\
};

\end{axis}
\end{tikzpicture}%

%% file: tikz/sufficientconditions2.tikz
%
%
\begin{tikzpicture}

\begin{axis}[%
width=4.822222in,
height=3.616667in,
at={(0.808889in,0.606667in)},
scale only axis,
separate axis lines,
every outer x axis line/.append style={black},
every x tick label/.append style={font=\color{black}},
xmin=-4,
xmax=6,
xtick={\empty},
every outer y axis line/.append style={black},
every y tick label/.append style={font=\color{black}},
ymin=-4.5,
ymax=3,
ytick={\empty}
]
\addplot [color=green,line width=2.0pt,only marks,mark=asterisk,mark options={solid},forget plot]
  table[row sep=crcr]{%
-2	1\\
};
\addplot [color=red,line width=2.0pt,only marks,mark=asterisk,mark options={solid},forget plot]
  table[row sep=crcr]{%
2	-1\\
};
\addplot [color=black,line width=2.0pt,only marks,mark=asterisk,mark options={solid},forget plot]
  table[row sep=crcr]{%
-0.3	0.7\\
};

\node[text=green] at (axis cs:-2,1) [anchor=south east] {\LARGE $u_1^*$};
\node at (axis cs:2,-1) [anchor=south west] {\LARGE $u_2^*$};
\node at (axis cs:-0.2,0.5) [anchor=south west] {\LARGE $u_0$};

\node at (axis cs:0.5,-1.5) [anchor=south east] {\large $\rho'_2$};
\node at (axis cs:5.2,2) [anchor = south] {\huge $D$};

\addplot [->,color=black,solid,line width=1.0pt,forget plot]
  table[row sep=crcr]{%
2	-1\\
-0.82842617152962	-2\\
};

\addplot [color=gray!20,dashed,line width=2.0pt,forget plot]
  table[row sep=crcr]{%
-3.8	1.50909090909091\\
-3.8027027027027	1.5\\
-3.82702702702703	1.4\\
-3.84594594594595	1.3\\
-3.85945945945946	1.2\\
-3.86756756756757	1.1\\
-3.87027027027027	1\\
-3.86756756756757	0.9\\
-3.85945945945946	0.8\\
-3.84594594594595	0.7\\
-3.82702702702703	0.6\\
-3.8027027027027	0.5\\
-3.8	0.490909090909091\\
-3.77142857142857	0.4\\
-3.73428571428571	0.3\\
-3.7	0.220000000000001\\
-3.69090909090909	0.2\\
-3.63939393939394	0.0999999999999996\\
-3.6	0.0315789473684211\\
-3.58064516129032	0\\
-3.51290322580645	-0.1\\
-3.5	-0.117391304347826\\
-3.43448275862069	-0.2\\
-3.4	-0.24\\
-3.34444444444444	-0.3\\
-3.3	-0.344444444444445\\
-3.24	-0.4\\
-3.2	-0.434482758620689\\
-3.11739130434783	-0.5\\
-3.1	-0.512903225806452\\
-3	-0.580645161290323\\
-2.96842105263158	-0.6\\
-2.9	-0.639393939393939\\
-2.8	-0.690909090909091\\
-2.78	-0.7\\
-2.7	-0.734285714285714\\
-2.6	-0.771428571428571\\
-2.50909090909091	-0.8\\
-2.5	-0.802702702702703\\
-2.4	-0.827027027027027\\
-2.3	-0.845945945945946\\
-2.2	-0.85945945945946\\
-2.1	-0.867567567567568\\
-2	-0.87027027027027\\
-1.9	-0.867567567567568\\
-1.8	-0.85945945945946\\
-1.7	-0.845945945945946\\
-1.6	-0.827027027027027\\
-1.5	-0.802702702702703\\
-1.49090909090909	-0.8\\
-1.4	-0.771428571428571\\
-1.3	-0.734285714285714\\
-1.22	-0.7\\
-1.2	-0.690909090909091\\
-1.1	-0.639393939393939\\
-1.03157894736842	-0.6\\
-1	-0.580645161290323\\
-0.9	-0.512903225806452\\
-0.882608695652174	-0.5\\
-0.8	-0.434482758620689\\
-0.760000000000001	-0.4\\
-0.7	-0.344444444444444\\
-0.655555555555556	-0.3\\
-0.6	-0.239999999999999\\
-0.565517241379311	-0.2\\
-0.5	-0.117391304347826\\
-0.487096774193548	-0.1\\
-0.419354838709677	0\\
-0.4	0.0315789473684211\\
-0.360606060606061	0.0999999999999996\\
-0.309090909090909	0.2\\
-0.3	0.220000000000001\\
-0.265714285714286	0.3\\
-0.228571428571429	0.4\\
-0.2	0.490909090909092\\
-0.197297297297297	0.5\\
-0.172972972972973	0.6\\
-0.154054054054054	0.7\\
-0.140540540540541	0.8\\
-0.132432432432433	0.9\\
-0.12972972972973	1\\
-0.132432432432433	1.1\\
-0.140540540540541	1.2\\
-0.154054054054054	1.3\\
-0.172972972972973	1.4\\
-0.197297297297297	1.5\\
-0.2	1.50909090909091\\
-0.228571428571428	1.6\\
-0.265714285714286	1.7\\
-0.3	1.78\\
-0.309090909090909	1.8\\
-0.360606060606061	1.9\\
-0.4	1.96842105263158\\
-0.419354838709677	2\\
-0.487096774193548	2.1\\
-0.5	2.11739130434783\\
-0.565517241379311	2.2\\
-0.6	2.24\\
-0.655555555555555	2.3\\
-0.7	2.34444444444444\\
-0.76	2.4\\
-0.8	2.43448275862069\\
-0.882608695652174	2.5\\
-0.9	2.51290322580645\\
-1	2.58064516129032\\
-1.03157894736842	2.6\\
-1.1	2.63939393939394\\
-1.2	2.69090909090909\\
-1.22	2.7\\
-1.3	2.73428571428571\\
-1.4	2.77142857142857\\
-1.49090909090909	2.8\\
-1.5	2.8027027027027\\
-1.6	2.82702702702703\\
-1.7	2.84594594594595\\
-1.8	2.85945945945946\\
-1.9	2.86756756756757\\
-2	2.87027027027027\\
-2.1	2.86756756756757\\
-2.2	2.85945945945946\\
-2.3	2.84594594594595\\
-2.4	2.82702702702703\\
-2.5	2.8027027027027\\
-2.50909090909091	2.8\\
-2.6	2.77142857142857\\
-2.7	2.73428571428571\\
-2.78	2.7\\
-2.8	2.69090909090909\\
-2.9	2.63939393939394\\
-2.96842105263158	2.6\\
-3	2.58064516129032\\
-3.1	2.51290322580645\\
-3.11739130434783	2.5\\
-3.2	2.43448275862069\\
-3.24	2.4\\
-3.3	2.34444444444444\\
-3.34444444444444	2.3\\
-3.4	2.24\\
-3.43448275862069	2.2\\
-3.5	2.11739130434783\\
-3.51290322580645	2.1\\
-3.58064516129032	2\\
-3.6	1.96842105263158\\
-3.63939393939394	1.9\\
-3.69090909090909	1.8\\
-3.7	1.78\\
-3.73428571428571	1.7\\
-3.77142857142857	1.6\\
-3.8	1.50909090909091\\
};
\addplot [color=gray!20,dashed,line width=2.0pt,forget plot]
  table[row sep=crcr]{%
3.5	-0.5\\
3.46206896551724	-0.4\\
3.41724137931034	-0.3\\
3.4	-0.266666666666667\\
3.36296296296296	-0.2\\
3.3	-0.1\\
3.224	0\\
3.2	0.0285714285714291\\
3.13478260869565	0.0999999999999996\\
3.1	0.134782608695653\\
3.02857142857143	0.2\\
3	0.224\\
2.9	0.3\\
2.8	0.362962962962963\\
2.73333333333333	0.4\\
2.7	0.417241379310345\\
2.6	0.462068965517241\\
2.5	0.5\\
2.4	0.529032258064516\\
2.3	0.551612903225806\\
2.2	0.567741935483871\\
2.1	0.57741935483871\\
2	0.580645161290323\\
1.9	0.577419354838709\\
1.8	0.567741935483871\\
1.7	0.551612903225806\\
1.6	0.529032258064516\\
1.5	0.5\\
1.4	0.462068965517241\\
1.3	0.417241379310344\\
1.26666666666667	0.4\\
1.2	0.362962962962962\\
1.1	0.3\\
1	0.224\\
0.971428571428571	0.2\\
0.9	0.134782608695653\\
0.865217391304347	0.0999999999999996\\
0.800000000000001	0.0285714285714291\\
0.776	0\\
0.7	-0.0999999999999994\\
0.637037037037037	-0.2\\
0.600000000000001	-0.266666666666665\\
0.582758620689655	-0.3\\
0.537931034482758	-0.4\\
0.5	-0.5\\
0.470967741935484	-0.6\\
0.448387096774194	-0.7\\
0.432258064516129	-0.8\\
0.422580645161291	-0.9\\
0.419354838709677	-1\\
0.42258064516129	-1.1\\
0.432258064516129	-1.2\\
0.448387096774194	-1.3\\
0.470967741935484	-1.4\\
0.5	-1.5\\
0.537931034482758	-1.6\\
0.582758620689655	-1.7\\
0.600000000000001	-1.73333333333333\\
0.637037037037037	-1.8\\
0.7	-1.9\\
0.776	-2\\
0.800000000000001	-2.02857142857143\\
0.865217391304347	-2.1\\
0.9	-2.13478260869565\\
0.971428571428571	-2.2\\
1	-2.224\\
1.1	-2.3\\
1.2	-2.36296296296296\\
1.26666666666667	-2.4\\
1.3	-2.41724137931034\\
1.4	-2.46206896551724\\
1.5	-2.5\\
1.6	-2.52903225806452\\
1.7	-2.55161290322581\\
1.8	-2.56774193548387\\
1.9	-2.57741935483871\\
2	-2.58064516129032\\
2.1	-2.57741935483871\\
2.2	-2.56774193548387\\
2.3	-2.55161290322581\\
2.4	-2.52903225806452\\
2.5	-2.5\\
2.6	-2.46206896551724\\
2.7	-2.41724137931035\\
2.73333333333333	-2.4\\
2.8	-2.36296296296296\\
2.9	-2.3\\
3	-2.224\\
3.02857142857143	-2.2\\
3.1	-2.13478260869565\\
3.13478260869565	-2.1\\
3.2	-2.02857142857143\\
3.224	-2\\
3.3	-1.9\\
3.36296296296296	-1.8\\
3.4	-1.73333333333333\\
3.41724137931034	-1.7\\
3.46206896551724	-1.6\\
3.5	-1.5\\
3.52903225806452	-1.4\\
3.55161290322581	-1.3\\
3.56774193548387	-1.2\\
3.57741935483871	-1.1\\
3.58064516129032	-1\\
3.57741935483871	-0.9\\
3.56774193548387	-0.8\\
3.55161290322581	-0.7\\
3.52903225806452	-0.6\\
3.5	-0.5\\
};
\addplot [color=blue,dashed,line width=2.0pt,forget plot]
  table[row sep=crcr]{%
4.9	-0.233333333333335\\
4.89122807017544	-0.2\\
4.86140350877193	-0.1\\
4.8280701754386	0\\
4.8	0.0761904761904764\\
4.79090909090909	0.0999999999999996\\
4.74909090909091	0.2\\
4.70363636363636	0.3\\
4.7	0.307407407407407\\
4.65283018867925	0.4\\
4.6	0.496551724137932\\
4.59803921568627	0.5\\
4.53725490196078	0.6\\
4.5	0.657575757575757\\
4.47142857142857	0.7\\
4.4	0.8\\
4.32127659574468	0.9\\
4.3	0.925641025641026\\
4.23555555555556	1\\
4.2	1.0390243902439\\
4.14186046511628	1.1\\
4.1	1.14186046511628\\
4.0390243902439	1.2\\
4	1.23555555555556\\
3.92564102564103	1.3\\
3.9	1.32127659574468\\
3.8	1.4\\
3.7	1.47142857142857\\
3.65757575757576	1.5\\
3.6	1.53725490196078\\
3.5	1.59803921568627\\
3.49655172413793	1.6\\
3.4	1.65283018867925\\
3.30740740740741	1.7\\
3.3	1.70363636363636\\
3.2	1.74909090909091\\
3.1	1.79090909090909\\
3.07619047619048	1.8\\
3	1.8280701754386\\
2.9	1.86140350877193\\
2.8	1.89122807017544\\
2.76666666666667	1.9\\
2.7	1.91694915254237\\
2.6	1.93898305084746\\
2.5	1.95762711864407\\
2.4	1.9728813559322\\
2.3	1.98474576271186\\
2.2	1.99322033898305\\
2.1	1.99830508474576\\
2	2\\
1.9	1.99830508474576\\
1.8	1.99322033898305\\
1.7	1.98474576271186\\
1.6	1.9728813559322\\
1.5	1.95762711864407\\
1.4	1.93898305084746\\
1.3	1.91694915254237\\
1.23333333333333	1.9\\
1.2	1.89122807017544\\
1.1	1.86140350877193\\
1	1.8280701754386\\
0.923809523809524	1.8\\
0.9	1.79090909090909\\
0.800000000000001	1.74909090909091\\
0.7	1.70363636363636\\
0.692592592592593	1.7\\
0.600000000000001	1.65283018867925\\
0.503448275862068	1.6\\
0.5	1.59803921568627\\
0.4	1.53725490196078\\
0.342424242424242	1.5\\
0.3	1.47142857142857\\
0.2	1.4\\
0.100000000000001	1.32127659574468\\
0.0743589743589741	1.3\\
0	1.23555555555556\\
-0.0390243902439023	1.2\\
-0.0999999999999996	1.14186046511628\\
-0.141860465116279	1.1\\
-0.2	1.0390243902439\\
-0.235555555555555	1\\
-0.3	0.925641025641026\\
-0.321276595744681	0.9\\
-0.4	0.8\\
-0.471428571428572	0.7\\
-0.5	0.657575757575757\\
-0.537254901960784	0.6\\
-0.598039215686274	0.5\\
-0.6	0.496551724137932\\
-0.652830188679245	0.4\\
-0.7	0.307407407407408\\
-0.703636363636364	0.3\\
-0.74909090909091	0.2\\
-0.790909090909091	0.0999999999999996\\
-0.8	0.0761904761904764\\
-0.828070175438596	0\\
-0.86140350877193	-0.1\\
-0.891228070175438	-0.2\\
-0.9	-0.233333333333334\\
-0.916949152542373	-0.3\\
-0.938983050847458	-0.4\\
-0.957627118644068	-0.5\\
-0.972881355932203	-0.6\\
-0.984745762711864	-0.7\\
-0.993220338983051	-0.8\\
-0.998305084745763	-0.9\\
-1	-1\\
-0.998305084745763	-1.1\\
-0.993220338983051	-1.2\\
-0.984745762711864	-1.3\\
-0.972881355932203	-1.4\\
-0.957627118644068	-1.5\\
-0.938983050847458	-1.6\\
-0.916949152542373	-1.7\\
-0.9	-1.76666666666667\\
-0.891228070175438	-1.8\\
-0.86140350877193	-1.9\\
-0.828070175438596	-2\\
-0.8	-2.07619047619048\\
-0.790909090909091	-2.1\\
-0.74909090909091	-2.2\\
-0.703636363636364	-2.3\\
-0.7	-2.30740740740741\\
-0.652830188679245	-2.4\\
-0.6	-2.49655172413793\\
-0.598039215686274	-2.5\\
-0.537254901960784	-2.6\\
-0.5	-2.65757575757576\\
-0.471428571428572	-2.7\\
-0.4	-2.8\\
-0.321276595744681	-2.9\\
-0.3	-2.92564102564103\\
-0.235555555555555	-3\\
-0.2	-3.0390243902439\\
-0.141860465116279	-3.1\\
-0.0999999999999996	-3.14186046511628\\
-0.0390243902439023	-3.2\\
0	-3.23555555555556\\
0.0743589743589741	-3.3\\
0.100000000000001	-3.32127659574468\\
0.2	-3.4\\
0.3	-3.47142857142857\\
0.342424242424242	-3.5\\
0.4	-3.53725490196078\\
0.5	-3.59803921568627\\
0.50344827586207	-3.6\\
0.600000000000001	-3.65283018867925\\
0.692592592592593	-3.7\\
0.7	-3.70363636363636\\
0.800000000000001	-3.74909090909091\\
0.9	-3.79090909090909\\
0.923809523809524	-3.8\\
1	-3.8280701754386\\
1.1	-3.86140350877193\\
1.2	-3.89122807017544\\
1.23333333333333	-3.9\\
1.3	-3.91694915254237\\
1.4	-3.93898305084746\\
1.5	-3.95762711864407\\
1.6	-3.9728813559322\\
1.7	-3.98474576271186\\
1.8	-3.99322033898305\\
1.9	-3.99830508474576\\
2	-4\\
2.1	-3.99830508474576\\
2.2	-3.99322033898305\\
2.3	-3.98474576271186\\
2.4	-3.9728813559322\\
2.5	-3.95762711864407\\
2.6	-3.93898305084746\\
2.7	-3.91694915254237\\
2.76666666666667	-3.9\\
2.8	-3.89122807017544\\
2.9	-3.86140350877193\\
3	-3.8280701754386\\
3.07619047619048	-3.8\\
3.1	-3.79090909090909\\
3.2	-3.74909090909091\\
3.3	-3.70363636363636\\
3.30740740740741	-3.7\\
3.4	-3.65283018867925\\
3.49655172413793	-3.6\\
3.5	-3.59803921568627\\
3.6	-3.53725490196078\\
3.65757575757576	-3.5\\
3.7	-3.47142857142857\\
3.8	-3.4\\
3.9	-3.32127659574468\\
3.92564102564103	-3.3\\
4	-3.23555555555556\\
4.0390243902439	-3.2\\
4.1	-3.14186046511628\\
4.14186046511628	-3.1\\
4.2	-3.0390243902439\\
4.23555555555556	-3\\
4.3	-2.92564102564103\\
4.32127659574468	-2.9\\
4.4	-2.8\\
4.47142857142857	-2.7\\
4.5	-2.65757575757576\\
4.53725490196078	-2.6\\
4.59803921568627	-2.5\\
4.6	-2.49655172413793\\
4.65283018867925	-2.4\\
4.7	-2.30740740740741\\
4.70363636363636	-2.3\\
4.74909090909091	-2.2\\
4.79090909090909	-2.1\\
4.8	-2.07619047619048\\
4.8280701754386	-2\\
4.86140350877193	-1.9\\
4.89122807017544	-1.8\\
4.9	-1.76666666666667\\
4.91694915254237	-1.7\\
4.93898305084746	-1.6\\
4.95762711864407	-1.5\\
4.9728813559322	-1.4\\
4.98474576271186	-1.3\\
4.99322033898305	-1.2\\
4.99830508474576	-1.1\\
5	-1\\
4.99830508474576	-0.9\\
4.99322033898305	-0.8\\
4.98474576271186	-0.7\\
4.9728813559322	-0.6\\
4.95762711864407	-0.5\\
4.93898305084746	-0.4\\
4.91694915254237	-0.3\\
4.9	-0.233333333333335\\
};
\addplot [color=black,dashed,line width=1.0pt,forget plot]
  table[row sep=crcr]{%
-0.3	0.7\\
-0.241025641025641	0.705647244405225\\
-0.182051282051282	0.606402605982062\\
-0.123076923076923	0.641396239442344\\
-0.0641025641025641	0.518943428117323\\
-0.00512820512820511	0.481672179475824\\
0.0538461538461538	0.435947669344627\\
0.112820512820513	0.340765549821398\\
0.171794871794872	0.268392497763863\\
0.230769230769231	0.221168970414089\\
0.28974358974359	0.34178075106348\\
0.348717948717949	0.167146357646752\\
0.407692307692308	0.24924222891995\\
0.466666666666667	0.175680351934917\\
0.525641025641026	0.16430618835216\\
0.584615384615385	0.133754247537772\\
0.643589743589744	-0.0694980417263975\\
0.702564102564102	-0.0622456120554403\\
0.761538461538461	0.011497181327626\\
0.82051282051282	-0.099246323008053\\
0.879487179487179	-0.0925129159040852\\
0.938461538461538	-0.218938534297256\\
0.997435897435897	-0.356155743936554\\
1.05641025641026	-0.277988033677144\\
1.11538461538462	-0.399934773531619\\
1.17435897435897	-0.384256797188263\\
1.23333333333333	-0.388334941061886\\
1.29230769230769	-0.45543991869409\\
1.35128205128205	-0.50283953135377\\
1.41025641025641	-0.577415596095952\\
1.46923076923077	-0.658857728990903\\
1.52820512820513	-0.66548998052645\\
1.58717948717949	-0.79283630376744\\
1.64615384615385	-0.716697248652666\\
1.7051282051282	-0.69045624441181\\
1.76410256410256	-0.906551720931106\\
1.82307692307692	-0.96211265655474\\
1.88205128205128	-0.835573490555363\\
1.94102564102564	-1.00702198339876\\
2	-1\\
};
\end{axis}
\end{tikzpicture}%

%% file: tikz/centering1.tikz
%
%
\begin{tikzpicture}

\begin{axis}[%
width=4.520833in,
height=3.390625in,
at={(0.758333in,0.56875in)},
scale only axis,
xmin=-4,
xmax=6,
xtick={\empty},
ymin=-3,
ymax=4.5,
ytick={\empty}
]
\addplot [color=red,line width=3.0pt,only marks,mark=asterisk,mark options={solid},forget plot]
  table[row sep=crcr]{%
-2	1\\
2	-1\\
};
\addplot [color=black,line width=3.0pt,only marks,mark=asterisk,mark options={solid},forget plot]
  table[row sep=crcr]{%
-0.3	0.7\\
};
\node at (axis cs:-2.4,1) [anchor=south east] {\LARGE $u_1^*$};
\node at (axis cs:2,-1) [anchor=south] {\LARGE $u_2^*$};
\node at (axis cs:-0.2,0.5) [anchor=south west] {\LARGE $u_0$};

\node at (axis cs:5.2,3.5) [anchor = south] {\huge $D$};

\node at (axis cs:-1.5,-0.1) [anchor = south] {\large $\rho^+$};

\addplot [->,color=black,solid,line width=1.0pt,forget plot]
  table[row sep=crcr]{%
-0.3 0.7\\
-1.9	-0.5\\
};

\addplot [color=blue,dashed,line width=2.0pt,forget plot]
  table[row sep=crcr]{%
-0.299999999999991	-1.3\\
-0.4	-1.2974358974359\\
-0.5	-1.28974358974359\\
-0.6	-1.27692307692308\\
-0.7	-1.25897435897436\\
-0.8	-1.23589743589744\\
-0.9	-1.20769230769231\\
-0.923076923076923	-1.2\\
-1	-1.17297297297297\\
-1.1	-1.13243243243243\\
-1.17058823529412	-1.1\\
-1.2	-1.08571428571429\\
-1.3	-1.03142857142857\\
-1.35238095238095	-1\\
-1.4	-0.96969696969697\\
-1.5	-0.9\\
-1.6	-0.819354838709678\\
-1.62222222222222	-0.8\\
-1.7	-0.727586206896552\\
-1.72758620689655	-0.7\\
-1.8	-0.622222222222222\\
-1.81935483870968	-0.6\\
-1.9	-0.5\\
-1.96969696969697	-0.4\\
-2	-0.352380952380953\\
-2.03142857142857	-0.3\\
-2.08571428571429	-0.2\\
-2.1	-0.170588235294119\\
-2.13243243243243	-0.0999999999999996\\
-2.17297297297297	0\\
-2.2	0.0769230769230772\\
-2.20769230769231	0.1\\
-2.23589743589744	0.2\\
-2.25897435897436	0.3\\
-2.27692307692308	0.4\\
-2.28974358974359	0.5\\
-2.2974358974359	0.6\\
-2.3	0.699999999999991\\
-2.2974358974359	0.8\\
-2.28974358974359	0.9\\
-2.27692307692308	1\\
-2.25897435897436	1.1\\
-2.23589743589744	1.2\\
-2.20769230769231	1.3\\
-2.2	1.32307692307692\\
-2.17297297297297	1.4\\
-2.13243243243243	1.5\\
-2.1	1.57058823529412\\
-2.08571428571429	1.6\\
-2.03142857142857	1.7\\
-2	1.75238095238095\\
-1.96969696969697	1.8\\
-1.9	1.9\\
-1.81935483870968	2\\
-1.8	2.02222222222222\\
-1.72758620689655	2.1\\
-1.7	2.12758620689655\\
-1.62222222222222	2.2\\
-1.6	2.21935483870968\\
-1.5	2.3\\
-1.4	2.36969696969697\\
-1.35238095238095	2.4\\
-1.3	2.43142857142857\\
-1.2	2.48571428571429\\
-1.17058823529412	2.5\\
-1.1	2.53243243243243\\
-1	2.57297297297297\\
-0.923076923076924	2.6\\
-0.9	2.60769230769231\\
-0.8	2.63589743589744\\
-0.7	2.65897435897436\\
-0.6	2.67692307692308\\
-0.5	2.68974358974359\\
-0.4	2.6974358974359\\
-0.3	2.7\\
-0.2	2.6974358974359\\
-0.0999999999999996	2.68974358974359\\
0	2.67692307692308\\
0.100000000000001	2.65897435897436\\
0.2	2.63589743589744\\
0.3	2.60769230769231\\
0.323076923076924	2.6\\
0.4	2.57297297297297\\
0.5	2.53243243243243\\
0.570588235294117	2.5\\
0.600000000000001	2.48571428571429\\
0.7	2.43142857142857\\
0.752380952380952	2.4\\
0.800000000000001	2.36969696969697\\
0.9	2.3\\
1	2.21935483870968\\
1.02222222222222	2.2\\
1.1	2.12758620689655\\
1.12758620689655	2.1\\
1.2	2.02222222222222\\
1.21935483870968	2\\
1.3	1.9\\
1.36969696969697	1.8\\
1.4	1.75238095238095\\
1.43142857142857	1.7\\
1.48571428571429	1.6\\
1.5	1.57058823529412\\
1.53243243243243	1.5\\
1.57297297297297	1.4\\
1.6	1.32307692307692\\
1.60769230769231	1.3\\
1.63589743589744	1.2\\
1.65897435897436	1.1\\
1.67692307692308	1\\
1.68974358974359	0.9\\
1.6974358974359	0.8\\
1.7	0.7\\
1.6974358974359	0.6\\
1.68974358974359	0.5\\
1.67692307692308	0.4\\
1.65897435897436	0.3\\
1.63589743589744	0.2\\
1.60769230769231	0.1\\
1.6	0.0769230769230761\\
1.57297297297297	0\\
1.53243243243243	-0.0999999999999996\\
1.5	-0.170588235294118\\
1.48571428571429	-0.2\\
1.43142857142857	-0.3\\
1.4	-0.352380952380953\\
1.36969696969697	-0.4\\
1.3	-0.5\\
1.21935483870968	-0.6\\
1.2	-0.622222222222223\\
1.12758620689655	-0.7\\
1.1	-0.727586206896552\\
1.02222222222222	-0.8\\
1	-0.819354838709677\\
0.9	-0.9\\
0.800000000000001	-0.969696969696969\\
0.752380952380952	-1\\
0.7	-1.03142857142857\\
0.600000000000001	-1.08571428571429\\
0.570588235294118	-1.1\\
0.5	-1.13243243243243\\
0.4	-1.17297297297297\\
0.323076923076923	-1.2\\
0.3	-1.20769230769231\\
0.2	-1.23589743589744\\
0.100000000000001	-1.25897435897436\\
0	-1.27692307692308\\
-0.0999999999999996	-1.28974358974359\\
-0.2	-1.2974358974359\\
-0.299999999999991	-1.3\\
};
\addplot [color=black,dashed,line width=1.0pt,forget plot]
  table[row sep=crcr]{%
-0.3	0.7\\
-0.343589743589744	0.655300847775446\\
-0.387179487179487	0.758637482949203\\
-0.430769230769231	0.775205537352269\\
-0.474358974358974	0.788897028238147\\
-0.517948717948718	0.781534196463602\\
-0.561538461538462	0.78042092486565\\
-0.605128205128205	0.755451444076663\\
-0.648717948717949	0.722850756799165\\
-0.692307692307692	0.75706150883907\\
-0.735897435897436	0.732994827041633\\
-0.779487179487179	0.728322130464779\\
-0.823076923076923	0.845004697036038\\
-0.866666666666667	0.776156727750367\\
-0.91025641025641	0.78315636782935\\
-0.953846153846154	0.795336969204956\\
-0.997435897435897	0.819125105520385\\
-1.04102564102564	0.84855303954517\\
-1.08461538461538	0.781488920240703\\
-1.12820512820513	0.887218639644166\\
-1.17179487179487	0.860930059244797\\
-1.21538461538462	0.904030455451275\\
-1.25897435897436	0.850976272550041\\
-1.3025641025641	0.870446274689689\\
-1.34615384615385	0.83112412274832\\
-1.38974358974359	0.853560596362385\\
-1.43333333333333	0.919536967435317\\
-1.47692307692308	0.887391787116704\\
-1.52051282051282	0.963202951924731\\
-1.56410256410256	0.877255546890528\\
-1.60769230769231	0.989379382223428\\
-1.65128205128205	0.94325939034609\\
-1.69487179487179	0.970983936472578\\
-1.73846153846154	1.01378514825788\\
-1.78205128205128	0.936080382916278\\
-1.82564102564103	0.958973473897942\\
-1.86923076923077	0.972703869918093\\
-1.91282051282051	1.01629023403442\\
-1.95641025641026	1.03049217697661\\
-2	1\\
};
\end{axis}
\end{tikzpicture}%

%% file: tikz/centering2.tikz
%
%
\begin{tikzpicture}

\begin{axis}[%
width=4.520833in,
height=3.390625in,
at={(0.758333in,0.56875in)},
scale only axis,
xmin=-4,
xmax=6,
xtick={\empty},
ymin=-3,
ymax=4.5,
ytick={\empty}
]
\addplot [color=red,line width=3.0pt,only marks,mark=asterisk,mark options={solid},forget plot]
  table[row sep=crcr]{%
2	-1\\
};
\addplot [color=green,line width=3.0pt,only marks,mark=asterisk,forget plot]
  table[row sep=crcr]{%
-2	1\\
};
\addplot [color=black,line width=3.0pt,only marks,mark=asterisk,mark options={solid},forget plot]
  table[row sep=crcr]{%
-0.3	0.7\\
};

\node[text=green] at (axis cs:-2.3,1) [anchor=south east] {\LARGE $u_1^*$};
\node at (axis cs:2,-1) [anchor=south] {\LARGE $u_2^*$};
\node at (axis cs:-0.2,0.5) [anchor=south west] {\LARGE $u_0$};

\node at (axis cs:5.2,3.5) [anchor = south] {\huge $D$};

\node at (axis cs:-1.5,-0.1) [anchor = south] {\large $\rho'^+$};

\addplot [color=gray!20,dashed,line width=2.0pt,forget plot]
  table[row sep=crcr]{%
-0.299999999999991	-1.3\\
-0.4	-1.2974358974359\\
-0.5	-1.28974358974359\\
-0.6	-1.27692307692308\\
-0.7	-1.25897435897436\\
-0.8	-1.23589743589744\\
-0.9	-1.20769230769231\\
-0.923076923076923	-1.2\\
-1	-1.17297297297297\\
-1.1	-1.13243243243243\\
-1.17058823529412	-1.1\\
-1.2	-1.08571428571429\\
-1.3	-1.03142857142857\\
-1.35238095238095	-1\\
-1.4	-0.96969696969697\\
-1.5	-0.9\\
-1.6	-0.819354838709678\\
-1.62222222222222	-0.8\\
-1.7	-0.727586206896552\\
-1.72758620689655	-0.7\\
-1.8	-0.622222222222222\\
-1.81935483870968	-0.6\\
-1.9	-0.5\\
-1.96969696969697	-0.4\\
-2	-0.352380952380953\\
-2.03142857142857	-0.3\\
-2.08571428571429	-0.2\\
-2.1	-0.170588235294119\\
-2.13243243243243	-0.0999999999999996\\
-2.17297297297297	0\\
-2.2	0.0769230769230772\\
-2.20769230769231	0.1\\
-2.23589743589744	0.2\\
-2.25897435897436	0.3\\
-2.27692307692308	0.4\\
-2.28974358974359	0.5\\
-2.2974358974359	0.6\\
-2.3	0.699999999999991\\
-2.2974358974359	0.8\\
-2.28974358974359	0.9\\
-2.27692307692308	1\\
-2.25897435897436	1.1\\
-2.23589743589744	1.2\\
-2.20769230769231	1.3\\
-2.2	1.32307692307692\\
-2.17297297297297	1.4\\
-2.13243243243243	1.5\\
-2.1	1.57058823529412\\
-2.08571428571429	1.6\\
-2.03142857142857	1.7\\
-2	1.75238095238095\\
-1.96969696969697	1.8\\
-1.9	1.9\\
-1.81935483870968	2\\
-1.8	2.02222222222222\\
-1.72758620689655	2.1\\
-1.7	2.12758620689655\\
-1.62222222222222	2.2\\
-1.6	2.21935483870968\\
-1.5	2.3\\
-1.4	2.36969696969697\\
-1.35238095238095	2.4\\
-1.3	2.43142857142857\\
-1.2	2.48571428571429\\
-1.17058823529412	2.5\\
-1.1	2.53243243243243\\
-1	2.57297297297297\\
-0.923076923076924	2.6\\
-0.9	2.60769230769231\\
-0.8	2.63589743589744\\
-0.7	2.65897435897436\\
-0.6	2.67692307692308\\
-0.5	2.68974358974359\\
-0.4	2.6974358974359\\
-0.3	2.7\\
-0.2	2.6974358974359\\
-0.0999999999999996	2.68974358974359\\
0	2.67692307692308\\
0.100000000000001	2.65897435897436\\
0.2	2.63589743589744\\
0.3	2.60769230769231\\
0.323076923076924	2.6\\
0.4	2.57297297297297\\
0.5	2.53243243243243\\
0.570588235294117	2.5\\
0.600000000000001	2.48571428571429\\
0.7	2.43142857142857\\
0.752380952380952	2.4\\
0.800000000000001	2.36969696969697\\
0.9	2.3\\
1	2.21935483870968\\
1.02222222222222	2.2\\
1.1	2.12758620689655\\
1.12758620689655	2.1\\
1.2	2.02222222222222\\
1.21935483870968	2\\
1.3	1.9\\
1.36969696969697	1.8\\
1.4	1.75238095238095\\
1.43142857142857	1.7\\
1.48571428571429	1.6\\
1.5	1.57058823529412\\
1.53243243243243	1.5\\
1.57297297297297	1.4\\
1.6	1.32307692307692\\
1.60769230769231	1.3\\
1.63589743589744	1.2\\
1.65897435897436	1.1\\
1.67692307692308	1\\
1.68974358974359	0.9\\
1.6974358974359	0.8\\
1.7	0.7\\
1.6974358974359	0.6\\
1.68974358974359	0.5\\
1.67692307692308	0.4\\
1.65897435897436	0.3\\
1.63589743589744	0.2\\
1.60769230769231	0.1\\
1.6	0.0769230769230761\\
1.57297297297297	0\\
1.53243243243243	-0.0999999999999996\\
1.5	-0.170588235294118\\
1.48571428571429	-0.2\\
1.43142857142857	-0.3\\
1.4	-0.352380952380953\\
1.36969696969697	-0.4\\
1.3	-0.5\\
1.21935483870968	-0.6\\
1.2	-0.622222222222223\\
1.12758620689655	-0.7\\
1.1	-0.727586206896552\\
1.02222222222222	-0.8\\
1	-0.819354838709677\\
0.9	-0.9\\
0.800000000000001	-0.969696969696969\\
0.752380952380952	-1\\
0.7	-1.03142857142857\\
0.600000000000001	-1.08571428571429\\
0.570588235294118	-1.1\\
0.5	-1.13243243243243\\
0.4	-1.17297297297297\\
0.323076923076923	-1.2\\
0.3	-1.20769230769231\\
0.2	-1.23589743589744\\
0.100000000000001	-1.25897435897436\\
0	-1.27692307692308\\
-0.0999999999999996	-1.28974358974359\\
-0.2	-1.2974358974359\\
-0.299999999999991	-1.3\\
};
\addplot [color=blue,dashed,line width=2.0pt,forget plot]
  table[row sep=crcr]{%
0.466666666666665	-2.2\\
0.4	-2.21694915254237\\
0.3	-2.23898305084746\\
0.2	-2.25762711864407\\
0.100000000000001	-2.2728813559322\\
0	-2.28474576271186\\
-0.0999999999999996	-2.29322033898305\\
-0.2	-2.29830508474576\\
-0.3	-2.3\\
-0.4	-2.29830508474576\\
-0.5	-2.29322033898305\\
-0.6	-2.28474576271186\\
-0.7	-2.2728813559322\\
-0.8	-2.25762711864407\\
-0.9	-2.23898305084746\\
-1	-2.21694915254237\\
-1.06666666666667	-2.2\\
-1.1	-2.19122807017544\\
-1.2	-2.16140350877193\\
-1.3	-2.1280701754386\\
-1.37619047619048	-2.1\\
-1.4	-2.09090909090909\\
-1.5	-2.04909090909091\\
-1.6	-2.00363636363636\\
-1.60740740740741	-2\\
-1.7	-1.95283018867925\\
-1.79655172413793	-1.9\\
-1.8	-1.89803921568627\\
-1.9	-1.83725490196078\\
-1.95757575757576	-1.8\\
-2	-1.77142857142857\\
-2.1	-1.7\\
-2.2	-1.62127659574468\\
-2.22564102564103	-1.6\\
-2.3	-1.53555555555556\\
-2.3390243902439	-1.5\\
-2.4	-1.44186046511628\\
-2.44186046511628	-1.4\\
-2.5	-1.3390243902439\\
-2.53555555555556	-1.3\\
-2.6	-1.22564102564103\\
-2.62127659574468	-1.2\\
-2.7	-1.1\\
-2.77142857142857	-1\\
-2.8	-0.957575757575758\\
-2.83725490196078	-0.9\\
-2.89803921568627	-0.8\\
-2.9	-0.796551724137931\\
-2.95283018867925	-0.7\\
-3	-0.607407407407407\\
-3.00363636363636	-0.6\\
-3.04909090909091	-0.5\\
-3.09090909090909	-0.4\\
-3.1	-0.376190476190476\\
-3.1280701754386	-0.3\\
-3.16140350877193	-0.2\\
-3.19122807017544	-0.0999999999999996\\
-3.2	-0.0666666666666656\\
-3.21694915254237	0\\
-3.23898305084746	0.1\\
-3.25762711864407	0.2\\
-3.2728813559322	0.3\\
-3.28474576271186	0.4\\
-3.29322033898305	0.5\\
-3.29830508474576	0.6\\
-3.3	0.7\\
-3.29830508474576	0.8\\
-3.29322033898305	0.9\\
-3.28474576271186	1\\
-3.2728813559322	1.1\\
-3.25762711864407	1.2\\
-3.23898305084746	1.3\\
-3.21694915254237	1.4\\
-3.2	1.46666666666667\\
-3.19122807017544	1.5\\
-3.16140350877193	1.6\\
-3.1280701754386	1.7\\
-3.1	1.77619047619048\\
-3.09090909090909	1.8\\
-3.04909090909091	1.9\\
-3.00363636363636	2\\
-3	2.00740740740741\\
-2.95283018867925	2.1\\
-2.9	2.19655172413793\\
-2.89803921568627	2.2\\
-2.83725490196078	2.3\\
-2.8	2.35757575757576\\
-2.77142857142857	2.4\\
-2.7	2.5\\
-2.62127659574468	2.6\\
-2.6	2.62564102564103\\
-2.53555555555556	2.7\\
-2.5	2.7390243902439\\
-2.44186046511628	2.8\\
-2.4	2.84186046511628\\
-2.3390243902439	2.9\\
-2.3	2.93555555555556\\
-2.22564102564103	3\\
-2.2	3.02127659574468\\
-2.1	3.1\\
-2	3.17142857142857\\
-1.95757575757576	3.2\\
-1.9	3.23725490196078\\
-1.8	3.29803921568627\\
-1.79655172413793	3.3\\
-1.7	3.35283018867925\\
-1.60740740740741	3.4\\
-1.6	3.40363636363636\\
-1.5	3.44909090909091\\
-1.4	3.49090909090909\\
-1.37619047619048	3.5\\
-1.3	3.5280701754386\\
-1.2	3.56140350877193\\
-1.1	3.59122807017544\\
-1.06666666666667	3.6\\
-1	3.61694915254237\\
-0.9	3.63898305084746\\
-0.8	3.65762711864407\\
-0.7	3.6728813559322\\
-0.6	3.68474576271186\\
-0.5	3.69322033898305\\
-0.4	3.69830508474576\\
-0.3	3.7\\
-0.2	3.69830508474576\\
-0.0999999999999996	3.69322033898305\\
0	3.68474576271186\\
0.100000000000001	3.6728813559322\\
0.2	3.65762711864407\\
0.3	3.63898305084746\\
0.4	3.61694915254237\\
0.466666666666665	3.6\\
0.5	3.59122807017544\\
0.600000000000001	3.56140350877193\\
0.7	3.5280701754386\\
0.776190476190477	3.5\\
0.800000000000001	3.49090909090909\\
0.9	3.44909090909091\\
1	3.40363636363636\\
1.00740740740741	3.4\\
1.1	3.35283018867925\\
1.19655172413793	3.3\\
1.2	3.29803921568628\\
1.3	3.23725490196078\\
1.35757575757576	3.2\\
1.4	3.17142857142857\\
1.5	3.1\\
1.6	3.02127659574468\\
1.62564102564103	3\\
1.7	2.93555555555556\\
1.7390243902439	2.9\\
1.8	2.84186046511628\\
1.84186046511628	2.8\\
1.9	2.7390243902439\\
1.93555555555556	2.7\\
2	2.62564102564103\\
2.02127659574468	2.6\\
2.1	2.5\\
2.17142857142857	2.4\\
2.2	2.35757575757576\\
2.23725490196078	2.3\\
2.29803921568627	2.2\\
2.3	2.19655172413793\\
2.35283018867925	2.1\\
2.4	2.00740740740741\\
2.40363636363636	2\\
2.44909090909091	1.9\\
2.49090909090909	1.8\\
2.5	1.77619047619048\\
2.5280701754386	1.7\\
2.56140350877193	1.6\\
2.59122807017544	1.5\\
2.6	1.46666666666667\\
2.61694915254237	1.4\\
2.63898305084746	1.3\\
2.65762711864407	1.2\\
2.6728813559322	1.1\\
2.68474576271186	1\\
2.69322033898305	0.9\\
2.69830508474576	0.8\\
2.7	0.700000000000018\\
2.69830508474576	0.6\\
2.69322033898305	0.5\\
2.68474576271186	0.4\\
2.6728813559322	0.3\\
2.65762711864407	0.2\\
2.63898305084746	0.1\\
2.61694915254237	0\\
2.6	-0.0666666666666692\\
2.59122807017544	-0.0999999999999996\\
2.56140350877193	-0.2\\
2.5280701754386	-0.3\\
2.5	-0.376190476190477\\
2.49090909090909	-0.4\\
2.44909090909091	-0.5\\
2.40363636363636	-0.6\\
2.4	-0.607407407407408\\
2.35283018867925	-0.7\\
2.3	-0.796551724137932\\
2.29803921568627	-0.8\\
2.23725490196078	-0.9\\
2.2	-0.957575757575758\\
2.17142857142857	-1\\
2.1	-1.1\\
2.02127659574468	-1.2\\
2	-1.22564102564103\\
1.93555555555556	-1.3\\
1.9	-1.3390243902439\\
1.84186046511628	-1.4\\
1.8	-1.44186046511628\\
1.7390243902439	-1.5\\
1.7	-1.53555555555556\\
1.62564102564103	-1.6\\
1.6	-1.62127659574468\\
1.5	-1.7\\
1.4	-1.77142857142857\\
1.35757575757576	-1.8\\
1.3	-1.83725490196078\\
1.2	-1.89803921568628\\
1.19655172413793	-1.9\\
1.1	-1.95283018867925\\
1.00740740740741	-2\\
1	-2.00363636363636\\
0.9	-2.04909090909091\\
0.800000000000001	-2.09090909090909\\
0.776190476190477	-2.1\\
0.7	-2.1280701754386\\
0.600000000000001	-2.16140350877193\\
0.5	-2.19122807017544\\
0.466666666666665	-2.2\\
};
\addplot [color=black,dashed,line width=1.0pt,forget plot]
  table[row sep=crcr]{%
-0.3	0.7\\
-0.241025641025641	0.628337239106672\\
-0.182051282051282	0.524161450634171\\
-0.123076923076923	0.573607903963025\\
-0.0641025641025641	0.49281082057641\\
-0.00512820512820511	0.417185087986414\\
0.0538461538461538	0.380250873260165\\
0.112820512820513	0.475902506672688\\
0.171794871794872	0.386360286749301\\
0.230769230769231	0.301385947673107\\
0.28974358974359	0.346529058950489\\
0.348717948717949	0.141315135468696\\
0.407692307692308	0.22603229166342\\
0.466666666666667	0.180586824452661\\
0.525641025641026	0.102115874799917\\
0.584615384615385	-0.0170073343430486\\
0.643589743589744	0.0220063726316734\\
0.702564102564102	-0.0810382430076832\\
0.761538461538461	-0.157790798049648\\
0.82051282051282	-0.18568482153336\\
0.879487179487179	-0.0928065367067091\\
0.938461538461538	-0.301094052827206\\
0.997435897435897	-0.310477047187015\\
1.05641025641026	-0.39181322412236\\
1.11538461538462	-0.357809434740961\\
1.17435897435897	-0.487086949650139\\
1.23333333333333	-0.353895063138619\\
1.29230769230769	-0.537591438649551\\
1.35128205128205	-0.601838717161802\\
1.41025641025641	-0.60262918418498\\
1.46923076923077	-0.616480774323559\\
1.52820512820513	-0.7309481725571\\
1.58717948717949	-0.595793849340777\\
1.64615384615385	-0.772042971771039\\
1.7051282051282	-0.822581918873698\\
1.76410256410256	-0.913231981377099\\
1.82307692307692	-0.909581974853217\\
1.88205128205128	-1.00355025984088\\
1.94102564102564	-0.955324627918716\\
2	-1\\
};

\addplot [->,color=black,solid,line width=1.0pt,forget plot]
  table[row sep=crcr]{%
-0.3 0.7\\
-2.7	-1.1\\
};
\end{axis}
\end{tikzpicture}%

%% file: tikz/eulerbeam-asym.tikz
%
%
%
\begin{tikzpicture}

\begin{axis}[
title={\LARGE{Buckling of an Euler elastica, $\mu = 1/2$}},
x label style={at={(axis description cs:0.5,-0.17)},anchor=south},
every x tick label/.append style={font=\color{black},font=\LARGE},
y label style={at={(axis description cs:-0.08,.5)},anchor=south},
every y tick label/.append style={font=\color{black},font=\LARGE},
xlabel={\LARGE{$\lambda$}},
ylabel={\LARGE{$\textrm{sign}(\theta'(0))\|\theta\|_2$}},
xmin=0, xmax=12.5663706143592,
ymin=-3, ymax=3,
axis on top=false,
width=11cm,
xmajorgrids,
ymajorgrids
]
\addplot [thick, black]
table {%
0 -0.0456435458031004
0.1 -0.0456897566215017
0.2 -0.0458289519464837
0.3 -0.0460628414125962
0.4 -0.0463943240035325
0.5 -0.0468275988152333
0.6 -0.0473682799954924
0.7 -0.0480235944900098
0.8 -0.0488026315384207
0.9 -0.0497166755298753
1 -0.0507796440175419
1.1 -0.0520086645970824
1.2 -0.0534248385326
1.3 -0.0550542599520089
1.4 -0.0569293909722448
1.5 -0.0590909416068735
1.6 -0.0615904795245506
1.7 -0.0644941174097412
1.8 -0.0678878283996072
1.9 -0.0718852849166987
2 -0.0766397221729893
2.1 -0.0823624040791957
2.2 -0.0893524545963223
2.3 -0.0980466780450128
2.4 -0.109106765944403
2.5 -0.123579196889196
2.6 -0.14320345917943
2.7 -0.171029990441205
2.8 -0.212638545883699
2.9 -0.277919650933159
3 -0.37916190595772
3.1 -0.513752137954878
3.2 -0.656227861525535
3.3 -0.787858760210531
3.4 -0.904652222229563
3.5 -1.00799988893724
3.6 -1.10012051234409
3.7 -1.18296961277066
3.8 -1.25810246809262
3.9 -1.32673355821612
4 -1.38981643918479
4.1 -1.44810929709234
4.2 -1.50222362592475
4.3 -1.55265951772005
4.4 -1.5998312650692
4.5 -1.64408611967987
4.6 -1.68571821946219
4.7 -1.72497907566649
4.8 -1.76208558067849
4.9 -1.79722620450789
5 -1.83056584980578
5.1 -1.86224970010705
5.2 -1.89240630289073
5.3 -1.92115006412944
5.4 -1.94858328515208
5.5 -1.97479783985625
5.6 -1.99987656656803
5.7 -2.02389443145658
5.8 -2.0469195075254
5.9 -2.06901380355461
6 -2.09023397007017
6.1 -2.11063190384295
6.2 -2.13025526812664
6.3 -2.149147942507
6.4 -2.16735041362112
6.5 -2.18490011594323
6.6 -2.2018317301937
6.7 -2.21817744561703
6.8 -2.23396719131723
6.9 -2.2492288409841
7 -2.26398839464618
7.1 -2.27827014051526
7.2 -2.2920967995167
7.3 -2.30548965471116
7.4 -2.31846866748925
7.5 -2.33105258215117
7.6 -2.34325902126802
7.7 -2.35510456685528
7.8 -2.366604844051
7.9 -2.37777458696862
8 -2.38862770370159
8.1 -2.39917733456517
8.2 -2.40943590529249
8.3 -2.41941517570743
8.4 -2.42912628433551
8.5 -2.43857978935786
8.6 -2.44778570626786
8.7 -2.45675354254763
8.8 -2.46549232964698
8.9 -2.47401065251554
9 -2.48231667691211
9.1 -2.49041817469091
9.2 -2.49832254724388
9.3 -2.50603684725915
9.4 -2.51356779894065
9.5 -2.5209218168181
9.6 -2.5281050232572
9.7 -2.53512326479331
9.8 -2.54198212736858
9.9 -2.54868695056216
10 -2.55524284089405
10.1 -2.56165468427354
10.2 -2.56792715765672
10.3 -2.57406473997242
10.4 -2.58007172237042
10.5 -2.58595221784127
10.6 -2.5917101702529
10.7 -2.59734936284534
10.8 -2.60287342622147
10.9 -2.6082858458688
11 -2.61358996924407
11.1 -2.61878901245061
11.2 -2.62388606653507
11.3 -2.62888410342932
11.4 -2.63378598156014
11.5 -2.63859445114839
11.6 -2.64331215921763
11.7 -2.64794165433028
11.8 -2.65248539106866
11.9 -2.65694573427653
12 -2.661324963076
12.1 -2.6656252746734
12.2 -2.6698487879668
12.3 -2.6739975469672
12.4 -2.67807352404414
12.5 -2.68207862300642
};
\addplot [thick, black]
table {%
3.4 0.324461834820978
3.5 0.199376215170977
3.6 0.149196856407026
3.7 0.119395249466346
3.8 0.0993037892868468
3.9 0.0847620269164794
4 0.0737289638735825
4.1 0.065067212247592
4.2 0.0580873330919445
4.3 0.0523452210842044
4.4 0.0475411802600391
4.5 0.0434653425701361
4.6 0.0399662278143075
4.7 0.0369316998995077
4.8 0.0342769463423688
4.9 0.0319366242451632
5 0.0298595753025338
5.1 0.0280051779446604
5.2 0.0263407729551105
5.3 0.0248398107802834
5.4 0.0234804948760914
5.5 0.0222447727929365
5.6 0.0211175754126523
5.7 0.0200862361597951
5.8 0.0191400437118788
5.9 0.0182698888401576
6 0.0174679947822292
6.1 0.0167276971125697
6.2 0.0160432713650142
6.3 -0.0154097951330435
6.4 -0.0148230381033071
6.5 -0.0142793748017768
6.6 -0.0137757162196247
6.7 -0.0133094575759297
6.8 -0.0128784403728215
6.9 -0.0124809276819331
7 -0.0121155923406553
7.1 -0.0117815184990134
7.2 -0.0114782178189225
7.3 -0.0112056626817188
7.4 -0.0109643401377913
7.5 -0.0107553322281691
7.6 -0.0105804310215133
7.7 -0.0104423007240712
7.8 -0.0103447053394418
7.9 -0.0102928299595867
8 -0.0102937393217438
8.1 -0.0103570433050016
8.2 -0.0104958833694615
8.3 -0.0107284370151256
8.4 -0.0110802869421705
8.5 -0.0115883052798377
8.6 -0.0123073330770949
8.7 -0.0133223479669319
8.8 -0.0147722490641476
8.9 -0.0169006586386786
9 -0.0201777244162837
9.1 -0.0256423017621069
9.2 -0.036113277849873
9.3 -0.0622660422508576
9.4 -0.152895490133912
9.5 -0.293957235041669
9.6 -0.404081112553464
9.7 -0.492101837454608
9.8 -0.566551091050174
9.9 -0.631774592433902
10 -0.690231563682853
10.1 -0.743456331982058
10.2 -0.792478271901261
10.3 -0.838026581707949
10.4 -0.880639921269754
10.5 -0.920729568058835
10.6 -0.958618055628813
10.7 -0.994563957277074
10.8 -1.02877843028734
10.9 -1.06143663371259
11 -1.09268583437059
11.1 -1.12265130346653
11.2 -1.15144069775193
11.3 -1.17914737566633
11.4 -1.20585294886782
11.5 -1.23162927432997
11.6 -1.25654003014508
11.7 -1.28064197680859
11.8 -1.30398597759536
11.9 -1.32661783209877
12 -1.34857896321371
12.1 -1.36990698795898
12.2 -1.39063619534774
12.3 -1.41079794921991
12.4 -1.43042103000263
12.5 -1.44953192638894
};
\addplot [blue, mark=square*, mark size=3, mark options={draw=black}, only marks]
table {%
3.4 0.324461834820978
};
\addplot [thick, black]
table {%
3.4 0.562356369837922
3.5 0.791220226041219
3.6 0.933938954991128
3.7 1.0469970399247
3.8 1.14261776082105
3.9 1.22617643738268
4 1.30066797870894
4.1 1.36798833027502
4.2 1.4294387619907
4.3 1.4859637990196
4.4 1.53827774741837
4.5 1.58693803521758
4.6 1.63239058088626
4.7 1.67499932552941
4.8 1.71506626248365
4.9 1.75284548913169
5 1.78855334823464
5.1 1.82237592539701
5.2 1.85447470794619
5.3 1.88499093371144
5.4 1.91404898621609
5.5 1.94175908265962
5.6 1.9682194286226
5.7 1.9935179646377
5.8 2.01773379621
5.9 2.04093837534452
6 2.06319648486302
6.1 2.08456706464251
6.2 2.10510390998278
6.3 2.1248562656689
6.4 2.14386933429288
6.5 2.16218471359113
6.6 2.17984077462441
6.7 2.19687299035397
6.8 2.21331422238521
6.9 2.22919497224414
7 2.24454360243251
7.1 2.25938653161018
7.2 2.27374840752959
7.3 2.28765226075964
7.4 2.30111964175623
7.5 2.31417074344242
7.6 2.3268245111357
7.7 2.33909874138992
7.8 2.35101017109401
7.9 2.36257455798125
8 2.37380675354407
8.1 2.38472077008186
8.2 2.39532983735838
8.3 2.40564646288977
8.4 2.41568247833653
8.5 2.42544908628069
8.6 2.43495690222291
8.7 2.44421599313421
8.8 2.45323591290205
8.9 2.46202573497161
9 2.47059408244878
9.1 2.47894915590205
9.2 2.48709875907407
9.3 2.49505032269169
9.4 2.50281092654255
9.5 2.51038731996955
9.6 2.51778594091857
9.7 2.52501293366145
9.8 2.53207416530409
9.9 2.53897524117873
10 2.54572151921037
10.1 2.55231812333817
10.2 2.55876995606591
10.3 2.56508171020804
10.4 2.57125787989257
10.5 2.57730277087609
10.6 2.58322051022174
10.7 2.58901505538619
10.8 2.59469020275829
10.9 2.60024959568798
11 2.60569673204111
11.1 2.61103497131302
11.2 2.61626754133074
11.3 2.62139754457177
11.4 2.62642796412482
11.5 2.63136166931604
11.6 2.63620142102273
11.7 2.64094987669441
11.8 2.64560959510006
11.9 2.6501830408188
12 2.65467258848992
12.1 2.65908052683721
12.2 2.66340906248144
12.3 2.66766032355377
12.4 2.67183636312202
12.5 2.67593916244116
};
\addplot [blue, mark=square*, mark size=3, mark options={draw=black}, only marks]
table {%
6.3 0.145621366655801
};
\addplot [thick, black]
table {%
6.3 -0.145621366655801
6.4 -0.382787845162594
6.5 -0.518586928949479
6.6 -0.623317340989801
6.7 -0.710914390208753
6.8 -0.787177567333629
6.9 -0.855195759296455
7 -0.916851430132238
7.1 -0.973395319924362
7.2 -1.02571029206748
7.3 -1.07444810166726
7.4 -1.12010684624185
7.5 -1.16307781987191
7.6 -1.20367538691384
7.7 -1.24215685321788
7.8 -1.2787361549547
7.9 -1.31359357040338
8 -1.34688278552453
8.1 -1.37873614714959
8.2 -1.40926864336925
8.3 -1.43858097023096
8.4 -1.46676192970517
8.5 -1.49389032969797
8.6 -1.52003650750764
8.7 -1.54526356453947
8.8 -1.56962837680849
8.9 -1.59318242933156
9 -1.61597251073418
9.1 -1.63804129583361
9.2 -1.65942783764815
9.3 -1.68016798557454
9.4 -1.70029474292295
9.5 -1.71983857429114
9.6 -1.73882767117487
9.7 -1.75728818259224
9.8 -1.77524441623139
9.9 -1.79271901462975
10 -1.80973311009607
10.1 -1.82630646144916
10.2 -1.84245757513238
10.3 -1.85820381284632
10.4 -1.87356148750117
10.5 -1.88854594901153
10.6 -1.90317166122581
10.7 -1.91745227109155
10.8 -1.93140067099928
10.9 -1.94502905511401
11 -1.95834897039248
11.1 -1.97137136288954
11.2 -1.98410661987746
11.3 -1.99656460823419
11.4 -2.00875470949862
11.5 -2.02068585194164
11.6 -2.032366539959
11.7 -2.04380488105557
11.8 -2.05500861065885
11.9 -2.06598511497232
12 -2.07674145205526
12.1 -2.08728437129489
12.2 -2.09762033141885
12.3 -2.1077555171798
12.4 -2.11769585483027
12.5 -2.12744702649345
};
\addplot [thick, black]
table {%
6.3 0.14562136668957
6.4 0.382787845330578
6.5 0.518586928949847
6.6 0.623317341033611
6.7 0.710914390208894
6.8 0.787177567333629
6.9 0.85519575929645
7 0.91685143029653
7.1 0.973395319965584
7.2 1.02571029207953
7.3 1.07444810167123
7.4 1.1201068462433
7.5 1.16307781987248
7.6 1.20367538691407
7.7 1.24215685321799
7.8 1.27873615495475
7.9 1.3135935704034
8 1.34688278552454
8.1 1.37873614714959
8.2 1.40926864336925
8.3 1.43858097023096
8.4 1.46676192970517
8.5 1.49389032969797
8.6 1.52003650750764
8.7 1.54526356453946
8.8 1.56962837680849
8.9 1.59318242933156
9 1.61597251073418
9.1 1.63804129583361
9.2 1.65942783764815
9.3 1.68016798557454
9.4 1.70029474292295
9.5 1.71983857429114
9.6 1.73882767117487
9.7 1.75728818259223
9.8 1.77524441623139
9.9 1.79271901462975
10 1.80973311009607
10.1 1.82630646144915
10.2 1.84245757513238
10.3 1.85820381284632
10.4 1.87356148750118
10.5 1.88854594901154
10.6 1.9031716612258
10.7 1.91745227109156
10.8 1.93140067099929
10.9 1.94502905511401
11 1.95834897039248
11.1 1.97137136288954
11.2 1.98410661987747
11.3 1.99656460823419
11.4 2.00875470949862
11.5 2.02068585194164
11.6 2.032366539959
11.7 2.04380488105557
11.8 2.05500861065885
11.9 2.06598511497232
12 2.07674145205526
12.1 2.08728437129489
12.2 2.09762033141885
12.3 2.1077555171798
12.4 2.11769585483027
12.5 2.12744702649346
};
\addplot [thick, black]
table {%
9.6 0.0460215039653755
9.7 0.0291213738392661
9.8 0.021489458994538
9.9 0.0171439495411652
10 0.0143488473014991
10.1 0.0124067224965454
10.2 0.0109824309545797
10.3 0.00989498281695972
10.4 0.00903819748790175
10.5 0.00834579727957645
10.6 0.00777436713835718
10.7 0.00729434767498683
10.8 0.00688496526991467
10.9 0.00653122882834346
11 0.00622207395487723
11.1 0.00594917466265032
11.2 0.00570615909737603
11.3 0.00548807683009909
11.4 0.00529103649720944
11.5 0.00511193606570299
11.6 0.00494828227486511
11.7 0.0047980519023039
11.8 0.00465959035950625
11.9 0.00453153566088819
12 0.00441276082789893
12.1 0.00430232988735916
12.2 0.00419946404172505
12.3 0.00410351556546207
12.4 0.00401394766585427
12.5 0.00393031903484466
};
\addplot [blue, mark=square*, mark size=3, mark options={draw=black}, only marks]
table {%
9.6 0.0460215039653755
};
\addplot [thick, black]
table {%
9.6 0.357701019939148
9.7 0.462795549011396
9.8 0.54503406478917
9.9 0.614744877805898
10 0.676124167349181
10.1 0.731404738927892
10.2 0.781952197185754
10.3 0.828677827021054
10.4 0.872227517516033
10.5 0.913079833520997
10.6 0.951601666780246
10.7 0.988082023424975
10.8 1.0227536361488
10.9 1.0558073863817
11 1.08740226928451
11.1 1.11767248306703
11.2 1.14673259948877
11.3 1.17468141665815
11.4 1.20160488413578
11.5 1.22757836055627
11.6 1.25266838169337
11.7 1.27693406329137
11.8 1.30042822722414
11.9 1.32319831516993
12 1.34528713705581
12.1 1.36673348955478
12.2 1.38757267132225
12.3 1.4078369153959
12.4 1.42755575456026
12.5 1.44675633202358
};
\addplot [thick, black]
table {%
3.4 0.562356369837922
3.39754485428713 0.552662998487628
3.39524398604855 0.542931837236986
3.39310315906615 0.533164204652564
3.39112848563688 0.523361616671099
3.3893264478801 0.513525815056988
3.3877039187931 0.503658799931324
3.38626818244274 0.493762866851163
3.38502695247608 0.483840648924106
3.38398838788065 0.473895164422189
3.38316110461957 0.46392987029976
3.38255418140148 0.453948721902559
3.38217715742191 0.44395623895927
3.3820400194344 0.433957577643994
3.38215317499668 0.423958608056633
3.38252740821883 0.413965995850992
3.38317381387568 0.403987285918667
3.3841037054165 0.394030984976944
3.38532849233488 0.38410663861643
3.38685952270587 0.374224896862764
3.38870788764604 0.364397560695837
3.39088418620423 0.354637600420324
3.39339825191524 0.344959135572308
3.39625884602382 0.335377365557907
3.3994733271117 0.325908440906436
3.40304731217533 0.316569267338784
};
\addplot [thick, black]
table {%
6.4 0.382787845330579
6.39478970146545 0.37425903358788
6.38967073904906 0.365675511618697
6.3846451535076 0.357037440123509
6.37971501437258 0.348345030156596
6.37488241828036 0.339598546070412
6.370149486999 0.330798308611283
6.3655183652184 0.32194469812875
6.36099121809169 0.3130381579031
6.35657022851981 0.304079197606105
6.35225759416966 0.295068396909465
6.34805552421901 0.286006409266173
6.34396623582221 0.276893965899072
6.33999195029421 0.26773188004796
6.33613488901018 0.258521051542245
6.3323972690222 0.249262471796258
6.32878129839712 0.239957229362604
6.32528917128209 0.2306065162309
6.32192306270735 0.221211635132594
6.31868512314131 0.211774008224566
6.31557747281364 0.202295187668364
6.31260219582853 0.192776868852039
6.30976133409351 0.183220907329626
6.30705688109273 0.173629341048669
6.30449077553781 0.164004420199433
6.30206489493443 0.154348648221186
6.29978104910444 0.144664839416022
6.29764097370876 0.134956201796479
6.29564632381818 0.125226459176435
6.29379866758184 0.115480035962132
6.29209948004557 0.105722345335636
6.29055013717363 0.0959602542996929
6.28915191012709 0.0862028645552712
6.28790595985345 0.0764628867880692
6.28681333203914 0.0667591994259354
6.28587495247645 0.0571219501509794
6.28509162289231 0.0476036192125789
6.28446401728321 0.038305630442342
6.28399267879585 0.0294507713256636
6.28367801718808 0.0216049338618031
6.28352030689926 0.0163092417413446
6.28351968574392 -0.0162849794744437
6.2836761542681 -0.0215499772283152
6.28398957573759 -0.0293836311876273
6.28445967678948 -0.0382334481338072
6.28508604872926 -0.0475290499473746
6.28586814945576 -0.0570461297902862
6.28680530599565 -0.0666826854263882
6.28789671761832 -0.0763859855775879
6.2891414594968 -0.0861257609952284
6.29053848687499 -0.0958830689362791
6.2920866396976 -0.105645162143015
6.29378464765416 -0.115402916821924
6.29563113558712 -0.125149451912283
6.29762462921089 -0.134879344861371
6.29976356108848 -0.144588164748103
6.30204627681186 -0.154272183048075
6.30447104133211 -0.163928188212273
6.30703604538843 -0.173553363191512
6.30973941198492 -0.183145202341035
6.31257920286849 -0.192701453646153
6.31555342496318 -0.202220077610752
6.3186600367197 -0.211699217333821
6.32189695434338 -0.221137176232759
6.32526205786694 -0.23053240106978
6.32875319703869 -0.239883468704936
6.33236819700066 -0.249189075497893
6.33610486373627 -0.258448028613865
6.33996098926891 -0.267659238707958
6.3439343565989 -0.276821713620403
6.34802274436779 -0.285934552818094
6.35222393124369 -0.294996942394946
6.35653570002392 -0.304008150496427
6.36095584145366 -0.312967523071069
6.36548215776161 -0.321874479878365
6.37011246591822 -0.330728510707775
6.37484460061862 -0.339529171767037
6.3796764170014 -0.34827608222275
6.38460579310835 -0.356968920870181
6.38963063209637 -0.365607422923276
6.39474886421224 -0.374191376917798
6.39995844854084 -0.382720621721132
6.40525737453817 -0.391195043645672
};
\addplot [thick, black]
table {%
9.6 0.0460215039653755
9.59040246922723 0.0488123390316122
9.58091374590165 0.0519512485718483
9.57157091561691 0.055498846311223
9.56242538298472 0.0595258047474339
9.55354797042977 0.0641116532073212
9.54503439763804 0.0693401816215167
9.53700897734376 0.0752892221821569
9.5296221065008 0.0820134337278602
9.52303586703565 0.0895227770007787
9.51739594238554 0.0977663323490338
9.51279890919138 0.106634126039385
9.50927252961425 0.115980180788257
9.50677971831128 0.12565431065914
9.50523987696233 0.13552612257918
9.50455318679613 0.145494716743941
9.50461839932311 0.155487674557904
9.50534244269589 0.165455429030468
9.50664405821142 0.175365063656899
9.5084541193141 0.185195189303677
9.51071447577507 0.194932195495964
9.51337633531845 0.204567663935496
9.51639865482687 0.21409663338796
9.51974672332439 0.223516448995287
9.52339098147428 0.232826000728529
9.52730606475174 0.242025217681912
9.53147003792678 0.251114730077388
9.53586378559091 0.260095641504261
9.54047052695445 0.268969374035749
9.54527542847459 0.277737561909767
9.55026529314331 0.286401977853447
9.55542830982065 0.294964481556244
9.56075384970259 0.303426983316126
9.56623229993435 0.311791418187955
9.57185492663713 0.320059727468421
9.57761376135667 0.328233845359804
9.58350150627473 0.336315689328162
9.58951145454315 0.344307153121033
9.59563742288308 0.352210101720196
9.60187369419796 0.360026367719625
};
\addplot [thick, red, mark=*, mark size=3, mark options={draw=black}, only marks]
table {%
0.5 -0.0468275988152333
3.5 -1.00799988893724
6.5 -2.18490011594323
9.5 -2.5209218168181
12.5 -2.68207862300642
};
\end{axis}

\end{tikzpicture}